\documentclass[12pt]{amsart}
\usepackage{amssymb}
\usepackage{amsmath}
\usepackage{amsbsy}
\usepackage{amscd}
\usepackage[mathscr]{eucal}
\usepackage[all]{xy}
\advance\textheight by \topskip
\makeatletter
\def\cal{\mathcal}
\def\Bbb{\mathbb}
\def\frak{\mathfrak}

\newenvironment{pf*}[1]{\proof[#1]}{\endproof}
\renewcommand{\thesubsection}{\thesection(\@roman\c@subsection)}
\makeatother
\newtheorem{Theorem}[equation]{Theorem}
\newtheorem{Corollary}[equation]{Corollary}
\newtheorem{Lemma}[equation]{Lemma}
\newtheorem{Proposition}[equation]{Proposition}

\theoremstyle{definition}
\newtheorem{Definition}[equation]{Definition}
\newtheorem{Example}[equation]{Example}

\makeatletter
\renewcommand\section{\@startsection{section}{1}%
  {\z@}{.7\linespacing\@plus\linespacing}{.5\linespacing}%
  {\reset@font\normalfont\bfseries\centering}}
\makeatother

\theoremstyle{remark}
\newtheorem{Remark}[equation]{Remark}
\newtheorem*{Claim}{Claim}

\newtheorem*{Acknowledgements}{Acknowledgements}
\numberwithin{equation}{section}
\numberwithin{figure}{section}

\newcommand{\thmref}[1]{Theorem~\ref{#1}}
\newcommand{\secref}[1]{\S\ref{#1}}
\newcommand{\lemref}[1]{Lemma~\ref{#1}}
\newcommand{\propref}[1]{Proposition~\ref{#1}}
\newcommand{\remref}[1]{Remark~\ref{#1}}
\newcommand{\corref}[1]{Corollary~\ref{#1}}
\newcommand{\subsecref}[1]{\S\ref{#1}}

\newcommand{\romnum}[1]{\romannumeral #1} 
\newcommand{\Romnum}[1]{\expandafter\uppercase\expandafter{\romannumeral #1}} 
\newcommand{\C}{{\Bbb C}}
\newcommand{\Z}{{\Bbb Z}}
\newcommand{\Q}{{\Bbb Q}}
\newcommand{\R}{{\Bbb R}}
\newcommand{\HH}{{\Bbb H}}

\newcommand{\CP}{\operatorname{\C P}}
\newcommand{\RP}{\operatorname{\R P}}

\newcommand{\SO}{\operatorname{\rm SO}}
\newcommand{\U}{\operatorname{\rm U}}
\newcommand{\Pin}{\operatorname{\rm Pin}}
\newcommand{\Spin}{\operatorname{\rm Spin}}
\newcommand{\Spinc}{\Spin^{c}}
\newcommand{\Spincm}{\Spin^{c_-}}
\newcommand{\SP}{\operatorname{\rm Sp}}
\newcommand{\su}{\operatorname{\frak su}}

\newcommand{\Map}{\operatorname{Map}}

\newcommand{\Hom}{\operatorname{Hom}}

\newcommand{\Ker}{\operatorname{Ker}}
\newcommand{\Coker}{\operatorname{Coker}}
\newcommand{\Ima}{\operatorname{Im}}

\newcommand{\sign}{\operatorname{sign}}
\newcommand{\supp}{\operatorname{supp}}%

\newcommand{\ind}{\mathop{\text{\rm ind}}\nolimits}
\newcommand{\tr}{\operatorname{tr}}

\newcommand{\M}{{\cal M}} 
\newcommand{\G}{{\cal G}} 

\newcommand{\A}{{\cal A}} 
\newcommand{\CC}{{\cal C}} 
\newcommand{\B}{{\cal B}} 
\newcommand{\OO}{{\operatorname{O}}} %
\newcommand{\Sc}{{\cal S}}%
\newcommand{\I}{{\cal I}}%
\newcommand{\D}{{\cal D}}%
\newcommand{\K}{{\cal K}}%
\newcommand{\UU}{{\cal U}}%
\newcommand{\E}{{\cal E}}%

\newcommand{\Lie}{\operatorname{Lie}}

\newcommand{\id}{\operatorname{id}}

\newcommand{\SW}{\operatorname{SW}}
\newcommand{\X}{\tilde{X}}

\newcommand{\deux}{\{\pm 1\}}
\begin{document}
\title[$\Pin^-(2)$-monopole invariants]{$\Pin^-(2)$-monopole invariants}
%
%
%
%
\author{Nobuhiro Nakamura}
\address{Department of Mathematics, Gakushuin University, 1-5-1, Mejiro, Toshima-ku, Tokyo, 171-8588, Japan}
\email{nobuhiro@math.gakushuin.ac.jp}
%
%
\begin{abstract}
We introduce a diffeomorphism invariant of $4$-manifolds, the $\Pin^-(2)$-mono-pole invariant,  defined 
by using the $\Pin^-(2)$-monopole equations. 
We compute the invariants of several $4$-manifolds, and prove gluing formulae.
By using the invariants, we construct exotic smooth structures on the connected sum of an elliptic surface $E(n)$ with arbitrary number of the $4$-manifolds of the form of $S^2\times\Sigma$ or $S^1\times Y$ where $\Sigma$ is a compact Riemann surface with positive genus and $Y$ is a closed $3$-manifold. 
As another application, we give an estimate of the genus of surfaces embedded in a $4$-manifold $X$ representing a class $\alpha\in H_2(X;l)$, where $l$ is a local coefficient on $X$.
\end{abstract}
%
%
%
%
%
%
\maketitle
%
%
%
%
\section{Introduction}\label{sec:intro}
%
%
In the paper \cite{Pin2}, we introduced the $\Pin^-(2)$-monopole equations which are a twisted or a real version of the Seiberg-Witten equations, and obtained several constraints on the intersection forms with local coefficients of $4$-manifolds by analyzing the moduli spaces.
In this article, we investigate diffeomorphism invariants defined by using the $\Pin^-(2)$-monopole equations,  which we will call {\it $\Pin^-(2)$-monopole invariants}. 
We compute the invariants of several $4$-manifolds, and prove connected-sum formulae.
We give two applications. 
The first application is to construct exotic smooth structures on $E(n) \#(\#_{i=1}^k(S^2\times\Sigma_i))\#(\#_{j=1}^l(S^1\times Y_j))$ where $\Sigma_i$ are compact Riemann surfaces with positive genus and $Y_j$ are closed $3$-manifolds. 
The second application is an estimate of the genus of surfaces embedded in a $4$-manifold $X$ representing a class $\alpha\in H_2(X;l)$, where $l$ is a local coefficient on $X$, which can be considered as a local coefficient analogue of the adjunction inequalities in the Seiberg-Witten theory \cite{KM, FS, MST, OS}.
%
%
\subsection{Exotic smooth structures}\label{subsec:exotic}
%
%
%
%
We state the first application:
\begin{Theorem}\label{thm:exotic}
For any positive integer $n$, there exists a set ${\cal S}_n$ of infinitely many distinct smooth structures on the elliptic surface $E(n)$ which have the following significance{\rm :}
For $\sigma\in{\cal S}_n$, let $E(n)_\sigma$ be the manifold with the smooth structure $\sigma$ homeomorphic to $E(n)$. 
Let $Z$ be a connected sum of arbitrary positive number of $4$-manifolds, each of which is  $S^2\times\Sigma$ or $S^1\times Y$ where $\Sigma$ is a compact Riemann surface with positive genus and $Y$ is a closed $3$-manifold. 
Then,   $E(n)_\sigma \# Z$ for different $\sigma$ are mutually non-diffeomorphic.
\end{Theorem}
\begin{Remark}
A famous result due to C.~T.~C.~Wall tells us that any pair of simply-connected smooth $4$-manifolds $M_1$ and $M_2$ which have isomorphic  intersection forms are {\it stably diffeomorphic} for stabilization by taking connected sums with $k(S^2\times S^2)$ for sufficiently large $k$. (See e.g. \cite{Kirby}.)
\thmref{thm:exotic} says that there exist infinitely many exotic structures on $E(n)$ which can not be stabilized by $S^2\times\Sigma$ with positive $g(\Sigma)$ or $S^1\times Y^3$.
\end{Remark}

%
%
\subsection{$\Pin^-(2)$-monopole invariants}\label{subsec:Pin2-inv}
%
%
To prove the theorem above,  the $\Pin^-(2)$-monopole invariant will be defined and used.
We remark that the $\Pin^-(2)$-monopole equations are defined on a {\it $\Spin^{c_-}$-structure} (\subsecref{subsec:spinc-} and \cite{Pin2}, Section 3), which is a $\Pin^-(2)$-analogue of $\Spin^c$-structure.
One of the special features of the $\Pin^-(2)$-monopole theory is that the moduli spaces may be nonorientable.
Hence, in general, $\Z_2$-valued invariants will be defined.
Only when the moduli space is orientable, $\Z$-valued invariants can be defined.
Here, we state several nonvanishing results on the $\Pin^-(2)$-monopole invariants.

A $\Spincm$-structure is an object on a double covering $\X\to X$ of a $4$-manifold $X$ rather than on $X$ itself.
For a $\Spin^{c_-}$-structure on $\X\to X$, an $\OO(2)$-bundle $E$  called the characteristic $\OO(2)$-bundle is associated (\subsecref{subsec:spinc-}).   
Let $l$ be the $\Z$-bundle associated to the double covering $\X\to X$, i.e., $l=\X\times_{\deux}\Z$.
The $l$-coefficient Euler class of $E$ in $H^2(X;l)$ is denoted by $\tilde{c}_1(E)$. 
More precisely, we need to fix an $l$-coefficient orientation of $E$ to define the Euler class $\tilde{c}_1(E)$. 
(See \subsecref{subsec:spinc-}.)

An Enriques surface $N_0$ has a double covering $\pi\colon K_0\to N_0$ with $K_0$  a $K3$ surface.
More generally, a smooth $4$-manifold $N$ which is homotopy equivalent to an Enriques surface is known to be homeomorphic to the standard Enriques surface \cite{Okonek}, and  has a double covering $\pi\colon K\to N$ such that $K$ is a homotopy $K3$ surface.  
Let $l_K= K\times_{\deux}\Z$.
\begin{Theorem}\label{thm:Enriques}
There exists a $\Spin^{c_-}$-structure $c$ on $\pi\colon K\to N$ which satisfies the following:
\begin{itemize}
\item $\pi^*\tilde{c}_1(E)=0$, where $E$ is the characteristic $\OO(2)$-bundle and $\pi^*\colon H^2(N;l_K)\to H^2(K;\Z)$ is the induced homomorphism.  
\item the $\Z_2$-valued $\Pin^-(2)$-monopole invariant of $(N,c)$ is nontrivial.
\end{itemize}
\end{Theorem}
\begin{Remark}
The virtual dimension of the moduli space of $(N,c)$ is $0$.
\end{Remark}
\begin{Remark}
\thmref{thm:Enriques} is proved by \thmref{thm:mod2} which relates the $\Pin^-(2)$-monopole invariants of $N$ with the Seiberg-Witten invariants of the double covering $K$, together with the non-vanishing result due to J.~Morgan and Z.~Szab\'{o} \cite{MS} for homotopy $K3$ surfaces.
\end{Remark}
Next we state a connected-sum formula for $\Pin^-(2)$-monopole invariants.
Before that, we note the following remarks.
In general, an ordinary $\Spin^c$-structure can be seen as a reduction of an {\it untwisted} $\Spin^{c_-}$-structure defined on a {\it trivial} double cover $\X\to X$ (\subsecref{subsec:spinc-}). 
Furthermore,  the Seiberg-Witten ($\U(1)$-monopole) equations on a $\Spin^c$-structure can be identified with the $\Pin^-(2)$-monopole equations on the corresponding untwisted $\Spin^{c_-}$-structure (\subsecref{subsec:untwisted}).
Often, we will not distinguish an untwisted $\Spin^{c_-}$-structure and the $\Spin^c$-structure which is its reduction, and use the same symbol.
In the following, we consider the gluing of $\Pin^-(2)$-monopoles and ordinary Seiberg-Witten $\U(1)$-monopoles. 

Let $X_1$ be a $4$-manifold with an ordinary $\Spin^c$-(or untwisted $\Spincm$-)structure $c_1$. 
Let $X_2$ be the manifold $Z$ in \thmref{thm:exotic} whose connected-summands are of the form of $S^2\times\Sigma$ or $S^1\times Y$. 
To define a $\Z$-bundle on $X_2$, consider a $2$-torus $T^2$ with a nontrivial $\Z$-bundle $l_T$.
An oriented Riemann surface $\Sigma$ with positive genus $g$ can be considered as a connected sum of $g$ tori: $\Sigma=T^2\#\cdots\# T^2$. 
Let $l_\Sigma$ be the $\Z$-bundle over $\Sigma$ which is given by the connected sum of $l_T$: $l_\Sigma=l_T\#\cdots \#l_T$.
For  a Riemann surface $\Sigma$ with positive genus, consider the product $S^2\times\Sigma$ with the $\Z$-bundle $l$ which is the pull-back 
$$
l=\pi^*l_{\Sigma},
$$ 
where $\pi\colon S^2\times\Sigma\to\Sigma$ is the projection.
We also consider $S^1\times Y$ with the $\Z$-bundle $l^\prime$ which is the pullback of a nontrivial $\Z$-bundle $l_{S^1}$ over $S^1$.
\begin{Remark}
For $(X;l)$, let $b_k^l=b_k(X;l)=\dim H^k(X;l\otimes\Q)$.
For $(X,l)=(S^2\times\Sigma_g;l)$, $b_0^l=b_2^l=b_4^l=0$ and $b_1^l=b_3^l=2g-2$.
For $(X,l)=(S^1\times Y;l)$, $b_k^l=0$ for all $k$.
\end{Remark}
Recall  $X_2$ is a connected-sum of $4$-manifolds of the form of $S^2\times\Sigma$ or $S^1\times Y$. 
Equip  each component of the form of $S^2\times \Sigma$ (resp. $S^1\times Y$) with the $\Z$-bundles $l$ (resp. $l^\prime$) as above, and define the $\Z$-bundle $l_{X_2}$ on $X_2$ as their connected sum. 
If we write the cardinality of $H^2(X_2;l_{X_2})$ as $n$, there are $n$ distinct isomorphism classes of $\Spincm$-structures for $\tilde{X}_2\to X_2$, where $\tilde{X}_2$ is the double covering associated to $l_{X_2}$. 
(See \propref{prop:Spincm}.)
Each of these $\Spin^{c_-}$-structures has a characteristic $\OO(2)$-bundle $E$ with torsion  $\tilde{c}_1(E)$. 
Let $c_2$ be such a $\Spin^{c_-}$-structure on $X_2$.
We consider the connected sum $X_1\# X_2$ with the $\Spin^{c_-}$-structure $c_1\#c_2$ which is the connected sum of the $\Spin^{c_-}$-structures $c_1$ and $c_2$.
(Here we assume $c_1$ is an untwisted $\Spincm$-structure.)
Then, the following holds:
\begin{Theorem}\label{thm:GF1}
Let $X_1$ be a closed oriented connected $4$-manifolds with a $\Spinc${\rm (}untwisted $\Spincm${\rm )}-structure such that
\begin{itemize}
\item $b_+(X_1)\geq 2$,
\item the virtual dimension of the Seiberg-Witten moduli space for $(X_1,c_1)$ is zero,
\item the Seiberg-Witten invariant for $(X_1,c_1)$ is odd.
\end{itemize}
Let $X_2$ and $l_{X_2}$ be as above.
Then, for any $\Spincm$-structure $c_2$ on $\X_2\to X_2$,  the $\Pin^-(2)$-monopole invariant of $(X_1\#X_2,c_1\#c_2)$ is nonzero. 
\end{Theorem}
\begin{Remark}
The virtual dimension $d$ of the moduli space of $(X_1\#X_2,c_1\#c_2)$ is positive: 
For instance, if $X_2=\#_{i=1}^k(S^2\times\Sigma_i)\#\#_{j=1}^m(S^1\times Y_j)$, then 
$$
d= \sum_{i=1}^k (2g(\Sigma_i)-2) + (k+m)=2\sum_{i=1}^kg(\Sigma_i)-k+m \geq k+m.
$$
\end{Remark}
\begin{Remark}
This non-vanishing result would be interesting because of the following two points:
First, although the dimension of the moduli space is positive, the (co)homological (not cohomotopical) invariant is nontrivial.
Second, if $X_2$ contains a component of the form of $S^2\times\Sigma$, all of the Seiberg-Witten invariants and the cohomotopy refinement \cite{BF} of $X_1\# X_2$ are $0$ because $S^2\times\Sigma$ admits a positive scalar curvature metric and $b_+(S^2\times\Sigma)>0$.
\end{Remark}
\begin{Remark}
It is worth to notice that $b_+(X_2;l)=0$. 
In fact, \thmref{thm:GF1} can be considered as a $\Pin^-(2)$-monopole analogue of the Seiberg-Witten gluing formulae for connected sums $X_1\# X_2$ when $X_1$ is a $4$-manifold with positive $b_+(X_1)$ and $X_2$ is one of the following:
\begin{enumerate}
\item  $X_2$ is a $4$-manifold with $b_1(X_2)=b_+(X_2)=0$, (Froyshov \cite{FC}, Chapter 14 for general cases; Fintushel-Stern \cite{FS}, Theorem 1.4 and Nicolaescu \cite{Nicolaescu}, \S4.6.2 for $\overline{\CP}^2$; Kotschick-Morgan-Taubes \cite{KMT}, Proposition 2 for rational homology $4$-spheres),   
\item $X_2=S^1\times S^3$,  (Ozsv\'{a}th-Szab\'{o} \cite{OS2}) or
\item $X_2$ is a connected sum of several manifolds in (1) or (2) above.
\end{enumerate}
\end{Remark}
\begin{Remark}
\thmref{thm:GF1} is a special case of \thmref{thm:GF}.
\end{Remark}
As mentioned above, the $\Pin^-(2)$-monopole invariants are defined as $\Z_2$-valued invariants.
But in some exceptional cases, we can define $\Z$-valued invariants.
For instance, the non-vanishing result for homotopy Enriques surfaces (\thmref{thm:Enriques}) is refined as follows:
\begin{Theorem}\label{thm:Enriques-ref}
The $\Z$-valued $\Pin^-(2)$-monopole invariant for $(N,c)$ in \thmref{thm:Enriques} is odd.
\end{Theorem}
Furthermore, the following holds for connected sums of homotopy Enriques surfaces.
\begin{Theorem}\label{thm:Enriques-sum}
For any integer $n\geq 2$, let $X_n=N_1\# N_2\#\cdots\# N_n$ where each $N_i$ is a homotopy Enriques surface.
Then $X_n$ has a $\Spincm$-structure $c_n$ such that 
\begin{itemize}
\item the $\Z_2$-valued $\Pin^-(2)$-monopole invariant is $0$, but
\item the $\Z$-valued invariant is nontrivial.
\end{itemize}
\end{Theorem}
\begin{Remark}
Since $b_+(N_i)\geq 1$, the Seiberg-Witten invariants and Donaldson invariants of $X_n$ are $0$.
\end{Remark}
%
%
%
%
\subsection{The genus of embedded surfaces}\label{subsec:genus}
%
%
We state the second application of the $\Pin^-(2)$-monopole invariants, which is an estimate of the genus of embedded surfaces representing a local-coefficient class. 
Let $X$ be a closed oriented connected $4$-manifold and suppose a nontrivial double covering $\X\to X$ is given, and let $l=\X\times_{\deux}\Z$.
Then a homology class $\alpha\in H_2(X;l)$ is represented by an embedded surface $\Sigma$ as follows:
\begin{itemize}
\item $\Sigma$ is a connected surface embedded in $X$. Let $i\colon \Sigma\to X$ be the embedding map.
\item The orientation system of $\Sigma$ is identified with the pull-back $i^*l$ of $l$ by $i$.
\item If $i_*\colon H_2(\Sigma;i^*l)\to H_2(X;l)$ is the induced homomorphism and $[\Sigma]\in H_2(\Sigma;i^*l)$ is the fundamental class, then $\alpha=i_*[\Sigma]$.
\end{itemize}
Conversely, a connected embedded surface $\Sigma$ whose orientation system is the restriction of $l$ has its fundamental class $[\Sigma]$ in $H_2(X;l)$.

For such embedded surfaces, the following adjunction inequality holds.
\begin{Theorem}\label{thm:genus}
Let $c$ be a $\Spin^{c_-}$-structure on $\X\to X$, and $\tilde{c}$ be the $\Spinc$-structure on $\X$ induced from $c$ {\rm (}see {\rm \secref{sec:Pin2})}.
Suppose at least one of the following occurs{\rm :}
\begin{itemize}
\item $b_+(X;l)\geq 2$ and the $\Pin^-(2)$-monopole invariant of $(X,c)$ is nontrivial.
\item $b_+(\X)\geq 2$ and the ordinary Seiberg-Witten invariant of $(\X,\tilde{c})$ is nontrivial.
\end{itemize}
Suppose a class $\alpha\in H_2(X;l)$ is represented by a connected embedded surface as above.
If $\alpha$ has infinite order and $\alpha\cdot\alpha\geq 0$, then
$$
-\chi(\Sigma)\geq |\tilde{c}_1(E)\cdot\alpha| + \alpha\cdot\alpha,
$$
where $\chi(\Sigma)$ is the Euler characteristic of $\Sigma$.
\end{Theorem}
Combining \thmref{thm:genus} with the non-vanishing 
results in \subsecref{subsec:Pin2-inv}, we obtain the 
following estimates for several concrete $4$-manifolds.
\begin{Theorem}\label{thm:genus-ap}
Suppose a pair $(X,l)$ of $4$-manifold $X$, and a $\Z$-bundle $l$ over $X$ is one of the following:
\begin{itemize}
\item $(N_1\# N_2\#\cdots\# N_n,l_1\#\cdots \#l_n)$, 
where each $N_i$ is a homotopy Enriques surface, and 
$l_i$ is a nontrivial $\Z$-bundle, or
\item $(E(2)\#Z,l)$ as in \thmref{thm:exotic}. 
\end{itemize}
Let $\Sigma$ be a connected embedded surface as above representing a class $\alpha\in H_2(X;l)$. 
If $\alpha$ has infinite order and $\alpha\cdot\alpha\geq 0$, then
$$
-\chi(\Sigma)\geq \alpha\cdot\alpha.
$$ 
\end{Theorem}
\begin{Remark}
The number $\alpha\cdot\alpha$ is the normal Euler number of the embedding $\Sigma\subset X$.
\end{Remark}
From this, we can also obtain some kind of equivariant adjunction inequality on the double coverings: 
\begin{Corollary}\label{cor:genus-cover}
Let $\X\to X$ be the double covering associated with $(X,l)$ in \thmref{thm:genus-ap}, and $\iota\colon\X\to\X$ be the covering transformation.
Suppose an oriented connected surface $\Sigma$ embedded in $\X$ satisfies the property that $[\Sigma]-\iota_*[\Sigma]$ has infinite order in $H_2(\X;\Z)$ and $[\Sigma]\cdot[\Sigma]\geq 0$. 
If $\Sigma\cap\iota(\Sigma)=\emptyset$, then 
\begin{equation}\label{eq:adjunction}
-\chi(\Sigma)\geq [\Sigma]\cdot[\Sigma].
\end{equation}
\end{Corollary}
\begin{Example}
Let us examine \corref{cor:genus-cover} for a simple example.
Let $X=K3\#(T^2\times S^2)$. 
Consider the double cover  $\X\to X$ which is associated to a nontrivial double cover $T^2\times S^2\to T^2\times S^2$.
Then   $\X=K_1\# (T^2\times S^2)\#K_2$, where $K_i$ are copies of $K3$.
Let $\sigma=[pt\times S^2]$ and $\tau=[T^2\times pt]$ in $H_2(T^2\times S^2;\Z)$.
Take a $2$-sphere $S$ representing $\sigma$ embedded in the $T^2\times S^2$-component, and oriented connected surfaces $\Sigma_i$ ($i=1,2$) embedded in the $K_i$-components so that $[\Sigma_i]\neq 0$, $[\Sigma_i]^2\geq 0$, $\iota(\Sigma_1)\cap \Sigma_2=\emptyset$ and $\iota_*[\Sigma_1]\neq [\Sigma_2]$. 
Then we can arrange to take a connected sum $\Sigma=\Sigma_1\#S\#\Sigma_2$ in $\X$ such that $\Sigma\cap\iota(\Sigma)=\emptyset$.
Such a $\Sigma$ certainly satisfies \eqref{eq:adjunction} because of the adjunction inequality for $K3$.
On the other hand, we can construct oriented connected surfaces $\Sigma$ embedded in $\X$ with $\Sigma\cap \iota(\Sigma)\neq \emptyset$ which violate \eqref{eq:adjunction} as follows.
Let $g_1$  be the genus of $\Sigma_1$ above.
We can take an embedded $2$-torus $T$ representing $\tau + n\sigma$ so that $2n>2g_1-[\Sigma_1]^2$.
Then take a connected sum $\Sigma=\Sigma_1\#T$ in $\X$.
Since $[\Sigma]\cdot\iota_*[\Sigma] = (\sigma+n\tau)^2=2n>0$, we have $\Sigma\cap \iota(\Sigma)\neq \emptyset$.
\end{Example}
The organization of the paper is as follows.
In Section 2, we introduce $\Pin^-(2)$-monopole invariants, and discuss the relation with the Seiberg-Witten invariants on the double covering, and prove \thmref{thm:Enriques} and \thmref{thm:Enriques-ref}.
In Section 3, several versions of gluing formulae are stated, and assuming these, we prove \thmref{thm:exotic} and \thmref{thm:Enriques-sum}. 
Sections 4-6 are devoted to the proof of the gluing theorems stated in \S2. 
Section 4 describes the $\Pin^-(2)$-monopole theory on $3$-manifolds.
Section 5 deals with finite energy $\Pin^-(2)$-monopoles on $4$-manifolds with tubular ends.
In Section 6, we give the proofs of  the gluing theorems. 
In Section 7, the proof of the genus estimate (\thmref{thm:genus}) is given. 
The Appendix provides some analytic detail of the gluing construction. 

\begin{Acknowledgements}
The author would like to thank M.~Furuta for helpful discussions and invaluable suggestions at various stages of this work which enabled the author to enrich the contents of this paper.
It is also pleasure to thank the referees for their detailed and valuable comments and pointing out mistakes in the previous  version of the paper.
\end{Acknowledgements}
%
%
\section{$\Pin^-(2)$-monopole invariants}\label{sec:Pin2}
%
%
\subsection{$\Spin^{c_-}$-structures}\label{subsec:spinc-}
%
%
The $\Pin^-(2)$-monopole equations are defined on $\Spin^{c_-}$-structures, which are a $\Pin^-(2)$-version of the $\Spin^c$-structures. 
While a $\Spin^c$-structure is given as a $\Spin^c(4)=\Spin(4)\times_{\{\pm 1\}}\U(1)$-lift of the frame bundle, a $\Spin^{c_-}$-structure is given by a $\Spin(4)\times_{\{\pm 1\}}\Pin^-(2)$-lift of it.
The precise definition is given as follows. (See also \cite{Pin2}, Section 3.)
The group $\Spin(4)\times_{\{\pm 1\}}\Pin^-(2)$ is denoted by $\Spin^{c_-}(4)$.
Let $X$ be a closed oriented connected  Riemannian $4$-manifold with double covering $\X\to X$. 
The $\SO(4)$-frame bundle on $X$ is denoted by $Fr(X)$. 
Since $\Pin^-(2)=\U(1)\cup j\U(1)$, $\Spin^c(4)$ is the identity component of $\Spin^{c_-}(4)$, and $\Spincm(4)/\Spinc(4)=\deux$. 
Also we have $\Spincm(4)/\Pin^-(2)=\SO(4)$ and $\Spincm(4)/\Spin(4)=\OO(2)$.
\begin{Definition}
A $\Spincm$-structure on $\X\to X$ is a triple $(P, \sigma, \tau)$ where 
\begin{itemize}
\item $P$ is a $\Spincm(4)$-bundle over $X$,
\item $\sigma$ is an isomorphism between the $\Z/2$-bundles $P/\Spinc(4)$ and $\X$,
\item $\tau$ is an isomorphism between the $\SO(4)$-bundles $P/\Pin^-(2)$ and $Fr(X)$.
\end{itemize}
\end{Definition}

Instead of the determinant $\U(1)$-bundle for a $\Spin^c$-structure, an $\OO(2)$-bundle $E=P/\Spin(4)$ is associated to a $\Spin^{c_-}$-structure. 
We call this $E$ the {\it characteristic $\OO(2)$-bundle}. 
Let $l$ be the $\Z$-bundle $\X\times_{\deux}\Z$ over $X$.
Then $l$ is related to $E$ by $\det E=l\otimes \R$.
The $l$-coefficient orientation of $E$ (and hence $\tilde{c}_1(E)\in H^2(X;l)$) is determined via the isomorphism $\sigma\colon P/\Spinc(4)\overset{\cong}{\to}\X$ as follows.
As described in \cite{Pin2}, \S 3.3, the $\Spinc(4)$-bundle $P\to P/\Spinc(4)\cong\X$ defines a $\Spinc$-structure on $\X$.
Let $L$ be its determinant line bundle, and $D(L)$, $S(L)$ be its disk and sphere bundles.
Let $E_\R$ be the $\R^2$-bundle associated to $E$, and $D(E_\R)$, $S(E_\R)$ be similar objects.
Then choose the $l$-coefficient orientation of $E$ so that the Thom classes $\tilde{u}\in H^2(D(L),S(L);\Z)$ of $L$ and $u\in H^2(D(E_\R),S(E_\R);l)$ of $E_\R$ satisfy the relation 
\begin{equation}\label{eq:u}
\pi^*u = \tilde{u},
\end{equation}
where $\pi^*$ is the homomorphism induced from the projection $\pi\colon\X\to X$.
Then we also have the relation $\pi^*\tilde{c}_1(E)=c_1(L)$.

The basic fact on $\Spincm$-structures on $\X\to X$ is as follows:
\begin{Proposition}\label{prop:Spincm}
\rm{(1)} For an $\OO(2)$-bundle $E$ over $X$ with $\det E = l\otimes\R$ as above, there exists a $\Spincm$-structure on $\X\to X$ whose characteristic bundle is isomorphic to $E$ if and only if $w_2(X) = w_2(E) + w_1(l\otimes\R)^2$.\\
\rm{(2)} If a $\Spincm$-structure on $\X\to X$ is given, there is a bijective correspondence between the set of isomorphism classes of $\Spincm$-structures on $\X\to X$ and $H^2(X;l)$.
\end{Proposition}
\proof
The assertion (1) is proved in \cite{Pin2}.
To prove the assertion (2), let us consider the exact sequence,
\begin{equation}\label{eq:exact}
1\to S^1\to \Spincm(4)\to \SO(4)\times\deux\to 1.
\end{equation}
From this, we have a fibration,
\begin{equation}\label{eq:fibration}
BS^1\to B\Spincm(4) \to B(\SO(4)\times\deux).
\end{equation}
In \eqref{eq:exact}, $\deux$ gives rise to an automorphism of $S^1$ of complex conjugation. 
If we identify $BS^1$ with $\CP^\infty$,  the action of  $\pi_1(B(\deux)) \cong \Z_2$ on a fiber of \eqref{eq:fibration} can be homotopically identified with complex conjugation on  $\CP^\infty$.
Then $\Spincm$-structures on $\X\to X$ are classified by 
$$
H^2(X;\tilde{\pi}_2(BS^1))\cong H^2(X;l),
$$
where $\tilde{\pi}_2$ is the local coefficient with respect to the $\pi_1(B(\deux))$-action on fibers.
\endproof

Usually, we will assume the covering $\X\to X$ is nontrivial.
But in the case when $\X\to X$ is trivial, the $\Spin^{c_-}(4)$-bundle of a $\Spin^{c_-}$-structure on $X$  has a $\Spin^c(4)$-reduction, and in fact, this reduction induces a $\Spin^c$-structure on $X$.
We will refer to a $\Spin^{c_-}$-structure with trivial $\X$ as an {\it untwisted $\Spin^{c_-}$-structure}.
%
%
%
\subsection{Definition of $\Pin^-(2)$-monopole invariants}
%
%
In this subsection, we introduce $\Pin^-(2)$-monopole invariants.
Let $X$ be an oriented closed connected $4$-manifold with double covering $\X\to X$, and suppose a $\Spin^{c_-}$-structure $c$ on $\X\to X$ is given.
Let $l=\X\times_{\deux}\Z$, $\lambda=l\otimes\R$, and $E$ be the characteristic $\OO(2)$-bundle.  
Then we have $\lambda=\det E$.
Let $\A$ be the space of $\OO(2)$-connections on $E$, $\CC$  the configuration space $\CC=\A\times\Gamma(S^+)$, and
$\CC^*$ the space of irreducible configurations, $\CC^*=\A\times(\Gamma(S^+)\setminus 0)$.
Fix $k\geq 3$ and take $L^2_k$-completion of $\CC$ and $\CC^*$.
The gauge transformation group $\G$ is the $L^2_{k+1}$-completion of $\Gamma(\X\times_{\{\pm 1\}}\U(1))$, where $\{\pm 1\}$ acts on $\U(1)$ by complex conjugation.
We use the same symbols for the completed spaces.
Let $\B^*=\CC^*/\G$.

The (perturbed) $\Pin^-(2)$-monopole equations for $(A,\Phi)\in \CC$ are given as follows:
\begin{equation}\label{eq:monopole}
\left\{
\begin{aligned}
D_A\Phi =& 0,\\
\frac12 F_A^+ =& q(\Phi) + \mu,
\end{aligned}
\right.
\end{equation}
where $D_A$ is the Dirac operator, $q$ is a quadratic form and $\mu\in\Omega^+(i\lambda)$.
(See Section 4 of \cite{Pin2} for the precise meaning and definition of each term of the equations.)
\begin{Remark}
Here we adopt the convention according to \cite{KM2}, slightly different from \cite{Pin2}, with $\frac12$ on the curvature term $F_A^+$.
Of course, this set of the equations is essentially same with that in \cite{Pin2}, because they coincide after an appropriate rescaling. 
\end{Remark}

The moduli space $\M(X,c)=\M_{\Pin^-(2)}(X,c)$ is defined as the space of solutions modulo gauge transformations. 
(The perturbed moduli space is usually denoted by the same symbol.)
\begin{Remark}
When the $\Spincm$-structure is untwisted, since $\X\to X$ is trivial, we have $\G=\Gamma(\X\times_{\{\pm 1\}}\U(1))\cong\Map(X,\U(1))$. 
While the stabilizer of the $\Pin^-(2)$-monopole reducible on a twisted $\Spincm$-structure is $\deux$, that in the untwisted case is $\U(1)$. 
(See also \subsecref{subsec:untwisted}.)
\end{Remark}
For the time being, we suppose the $\Spincm$-structure is twisted.
Suppose $b_+(X;l)\geq 1$.
Then, as in the case of the ordinary Seiberg-Witten theory, by a generic choice of $\mu$, the moduli space $\M(X,c)$ has no reducible and is a compact manifold whose dimension is given by
\begin{equation}\label{eq:dc}
d(c)=\frac14(\tilde{c}_1(E)^2 - \sign(X)) - (b_0(X;l) - b_1(X;l) + b_+(X;l)).
\end{equation}
Note that the index of the Dirac operator $D_A$ is given by $\frac14(\tilde{c}_1(E)^2 - \sign(X))$ and $b_0(X;l)=0$ if $l$ is nontrivial.

In a sense, the $\Pin^-(2)$-monopole invariant of $(X,c)$ is  defined as the fundamental class of the moduli space $[\M(X,c)]\in H_{d(c)}(\B^*)$.
We can  obtain a numerical invariant by evaluating $[\M(X,c)]$ by a cohomology class in $H^{d(c)}(\B^*)$. 
If $\X\to X$ is nontrivial, $\B^*$ has the homotopy type of the classifying space of the group $\Z/2\times \Z^{b_1(X;l)}$.
This fact is stated in \cite{Pin2}, Proposition 25.
However, the proof of Lemma 27 in \cite{Pin2} which is used in the proof of Proposition 25 is incomplete in that it is not proved there that the identity component of $\G$ is contractible.
Here we complement it.
\begin{Lemma}
The gauge transformation group $\G$ is homotopy equivalent to $(\Z/2)\times \Z^{b_1(X;l)}$.
\end{Lemma}
\proof 
Let $\tilde{\G} = \Map (\X,\U(1))$. 
Define the involution $I$ on $\tilde{\G}$ by $u\mapsto \overline{\iota^*u}$ where $\iota\colon\X\to\X$ is the covering transformation and `` $\bar{\cdot}$ " means the complex conjugation.
Then $\G$ is identified with the $I$-fixed point set $\tilde{\G}^I$.
Let $h\colon\tilde{\G}\to[\X,S^1]\cong H^1(\X ;\Z)\cong \Z^{b_1(\X)}$ be the map which sends each element of $\tilde{\G}$ to its homotopy class. 
Put $\tilde{{\cal K}} =\ker h$.
Consider the following diagram:
$$
\begin{CD}
1@>>>\tilde{\cal K}@>>> \tilde{\G} @>>> [X,S^1] @>{h}>>1\\
@. @AAA  @AAA @AA{j}A\\
1@>>>\tilde{\cal K} \cap \tilde{\G}^I @>>>\tilde{\G}^I @>>> h(\tilde{\G}^I)@>>>1 
\end{CD}
$$
The vertical map $j$ is injective since the first and second vertical maps are inclusions.
It is proved that $\pi_0\tilde{\G}^I = \pi_0\G = \Z_2\oplus \Z^{b_1(X;l)}$ in the proof of Lemma 27 in \cite{Pin2}. 
Now it suffices to see that  $\tilde{\cal K} \cap \tilde{\G}^I$ is homotopy equivalent to $\{\pm 1\}$.
Each element  $u\in\tilde{\cal K}$ can be written as $u=\exp(2\pi\sqrt{-1} f)$ for some function $f\colon\X\to \R$.
If $u=\exp(2\pi\sqrt{-1} f)$ is in $\tilde{\cal K} \cap \tilde{\G}^I$, then there is an integer $m$ so that $f(\iota x)= m- f(x)$ for every $x\in\X$.
If we fix a base point $x_0\in\X$ and choose $f$ so that $f(x_0)\in [0,1)$, then such an $m$ is uniquely determined.
Then the homotopy $f_t = tf +(1-t)m/2$ gives the homotopy between $u$ and $\pm 1$.
\endproof

In contrast to the ordinary Seiberg-Witten theory, the moduli space $\M(X,c)$ may be non-orientable.
(A necessary condition for $\M(X,c)$ to be orientable will be given in \subsecref{subsec:orientability}.)
In general, we can define the following  $\Z/2$-valued version of the $\Pin^-(2)$-monopole invariants.
\begin{Definition}
The $\Pin^-(2)$-monopole invariant of $(X,c)$ is defined as a map 
$$
\SW^{\Pin}(X,c)\colon H^{d(c)}(\B^*;\Z/2) \to \Z/2,
$$
given by
$$
\SW^{\Pin}(X,c)(\xi):= \langle\xi, [\M(X,c)]\rangle.
$$ 
\end{Definition}
If $b_+(X;l)\geq 2$, then $\SW^{\Pin}(X,c)$ is a diffeomorphism invariant. 
If $b_+(X;l)=1$, then  $\SW^{\Pin}(X,c)$  depends on the chamber structure of the space of metrics and perturbations.
\begin{Remark}
We give a geometric description of the cohomology classes of $\B^*$ in \subsecref{subsec:mu-map1} and \subsecref{subsec:mu-map2}.
\end{Remark}
\begin{Remark} 
The compactness of $\M(X,c)$ enables us to develop the Bauer-Furuta theory \cite{BF} for the $\Pin^-(2)$-monopole equations.
In fact, we can define a stable cohomotopy refinement of the $\Pin^-(2)$-monopole invariants.
This will be discussed elsewhere.
\end{Remark}
%
%
\subsection{Orientability of the moduli spaces}\label{subsec:orientability}
%
%
The purpose of this subsection is to discuss the orientability of the moduli spaces.
Let us consider the family of Dirac operators $\tilde{\delta}_{Dirac}=\{ D_A\}_{A\in\A}$. 
In \cite{Pin2}, \S4, we introduced a subgroup $\K_\gamma$ in $\G$, which has the properties:
\begin{itemize}
\item $\G/\K_\gamma=\{\pm 1\}$. 
\item $\K_\gamma$ acts on $\A$ freely, and $\A/\K_\gamma$ has the same homotopy type of  $H^1(X;\lambda)/H^1(X;l)$.
\end{itemize}
\begin{Remark}\label{rem:Kgamma}
Here $\gamma$ is a circle embedded in $X$ on which $l$ is nontrivial.
The subgroup $\K_\gamma$ is defined as the set of gauge transformations whose restrictions to $\gamma$ are homotopic to $1$. 
\end{Remark}
Dividing $\tilde{\delta}_{Dirac}$ by $\K_\gamma$, we obtain the family $\delta_{Dirac}=\tilde{\delta}_{Dirac}/\K_\gamma$ over $\A/\K_\gamma$.
\begin{Proposition}\label{prop:orientable}
If the index of the Dirac operator is even and $\det\ind\delta_{Dirac}$ is trivial, then the moduli space is orientable.
\end{Proposition}
\proof 
For a configuration $(A,\Phi)$, let us consider the sequence,
$$
\begin{CD}
0@>>>\Omega^0(i\lambda) @>{\I_\Phi}>>\Omega^1(i\lambda)\oplus \Gamma(S^+)@>{\D_{(A,\Phi)}}>> \Omega^+(i\lambda)\oplus \Gamma(S^-)@>>>0,
\end{CD}
$$
where $\I_\Phi(f)=(-2df,f\Phi)$ and $\D_{(A,\Phi)} (a,\phi)=d^+a-Dq_{\Phi}(\phi),D_A\phi+\frac12\rho(a)\Phi)$, which are the linearizations of the gauge group action and the monopole map.
Let $V=\Omega^1(i\lambda)\oplus \Gamma(S^+)$, and $W=(\Omega^0\oplus\Omega^+)(i\lambda)\oplus\Gamma(S^-)$ and define $\delta_{(A,\Phi)}\colon V\to W$ by,
$$
\delta_{(A,\Phi)} = \I_\Phi^*\oplus \D_{(A,\Phi)}.
$$
Then the family $\tilde{\delta}=\{\delta_{(A,\Phi)}\}_{(A,\Phi)\in\CC}$ defines a bundle homomorphism between the bundles over $\CC$, 
\begin{equation*}
\tilde{\delta}\colon \CC\times V\to\CC\times W.
\end{equation*}
Restricting $\tilde{\delta}$ to $\CC^*$ and dividing by $\G$, we obtain a bundle homomorphism over $\B^*=\CC^*/\G$,
$$
\delta\colon  \CC^*\times_{\G} V\to\CC^*\times_{\G} W.
$$
The moduli space is orientable if $\det\ind\delta$ is trivial.
By deforming $\delta_{(A,\Phi)}$ by $\delta_{(A,t\phi)}$ $(0\leq t\leq 1$), we may assume 
%
$\tilde{\delta}=\{(d^*\oplus d^+)\oplus D_A\}_{(A,\Phi)\in\CC}$.
%
Since $(d^*\oplus d^+)$ does not depend on $(A,\Phi)$,  $\det\ind (d^*\oplus d^+)$ is trivial.
Therefore it suffices to consider the Dirac family 
\begin{equation}\label{eq:td}
\tilde{\delta}^\prime=\{ D_A\}_{(A,\Phi)\in\CC}\colon \CC\times\Gamma(S^+)\to \CC\times\Gamma(S^-).
\end{equation}
Then \eqref{eq:td} can be identified with the pull-back of $\tilde{\delta}_{Dirac}$, via the projection $p\colon \CC\to\A$ with $p(A,\Phi)=A$.
Dividing \eqref{eq:td} by $\K_\gamma$, we obtain $\tilde{\delta}^\prime/\K\colon\CC\times_{\K_\gamma} \Gamma(S^+) \to \CC\times_{\K_\gamma}\Gamma(S^-)$.
Note that $\CC/\K_\gamma$ is homotopic to $\A/\K_\gamma$. 
Thus $\ind(\tilde{\delta}^\prime/\K)$ is identified with $p^*\ind(\delta_{Dirac})$, which is trivial by the assumption.
Hence $\det\ind\delta$ is trivial if and only if $\det \left((p^*\ind(\delta_{Dirac}))|_{\CC^*})/\deux\right)$ over $\CC^*/\G$ is trivial.
Note that $\CC^*/\G\simeq \RP^\infty\times T^{b_1(X;l)}$. 
Let $\eta\to \CC^*/\G$ be the nontrivial real line bundle which represents the generator of $H^1(\RP^\infty;\Z_2)$.  
Then by the assumptions, we see that $\det \left((p^*\ind(\delta_{Dirac})|_{\CC^*})/\deux\right)\cong \eta^{\otimes\ind D}$.
Thus the proposition is proved.
\endproof
\begin{Remark}
For instance, if $b_1(X;l)=0$ and the Dirac index is even, then the moduli space is orientable.
\end{Remark}
Note that $H^*(\B^*;\Z)/\mathrm{Tor}\cong H^*(T^{b_1(X;l)};\Z)$.
Suppose the moduli space $\M(X)$ is orientable. 
Fixing an orientation, we can define $\Z$-valued $\Pin^-(2)$-monopole invariants $\SW^{\Pin}_{\Z}$ by evaluating the fundamental class $[\M(X)]$ by infinite-order classes $\xi$ in $H^*(\B^*;\Z)$:
$$
\SW^{\Pin}_{\Z}(X,c)(\xi) = \langle\xi,[\M(X)]\rangle.
$$ 
%
%
\subsection{$\Pin^-(2)$-monopoles on untwisted $\Spincm$-structures}\label{subsec:untwisted}
%
%
Let us consider an untwisted $\Spincm$-structure $c=(P,\sigma,\tau)$ on a (trivial) double covering $\X\to X$.
The two connected components of $\X$ will be denoted by $X_+$ and $X_-$ according to the rule described below.
Consider the $\Spinc$-structure on $\X$ which is defined by the projection $P\to P/\Spinc(4)\cong \X$.
Its restrictions to the components $X_+$ and $X_-$ of $\X$ are mutually complex conjugate $\Spinc$-structures $c_+$ and $c_-$ (see \cite{Pin2}, \S2(\romnum{3})).
Let $i_\pm\colon X_\pm\to \X$ be the inclusion maps.
Let $L_\pm$ be the determinant line bundles of $c_\pm$, and their Thom classes be $u_\pm$.
Then $X_+$ is chosen to satisfy
$$
u_+=i_+^*(\tilde{u})=i_+^*\circ\pi^*(u),
$$
where $u$ and $\tilde{u}$ are the Thom classes as in \eqref{eq:u}.
We call the $\Spinc$-structure $c_+$ the {\it canonical reduction}.

\begin{Remark}
When a $\Spinc$-structure $c_0$ with $\Spinc(4)$-bundle $P_c\to X$ is given, the $\Spincm(4)$-bundle $P=P_c\times_{\Spinc(4)}\Spincm(4)$  defines an untwisted $\Spincm$-structure $c$ on $\X=P/\Spinc(4)\to X$.  
Then $c_0$ is the canonical reduction of $c$.
\end{Remark}
 
As {\it real} vector bundles, we have identifications among spinor bundles for $c$, $c_+$ and $c_-$,
$$
S^\pm_c \cong S^\pm_{c_+}\cong S^\pm_{c_-}.
$$
Also as {\it real} vector bundles,  we have identifications among the $\R^2$-vector bundle associated to the characteristic $\OO(2)$-bundle $E$ of $c$ and the determinant line bundles $L_{\pm}$.
If an $\OO(2)$-connection $A$ on $E$ is given, we have $\U(1)$-connections $A_\pm$  on $L_\pm$ induced from $A$ by reduction.
As {\it real} operators, the covariant derivatives of $A$ and $A_\pm$  can be identified, and therefore the Dirac operators induced from $A$ and  $A_\pm$ can also be identified as real operators.
Furthermore, it can be seen that the $\Pin^-(2)$-monopole solutions on $c$ can be identified with the Seiberg-Witten solutions on $c_\pm$ via the identifications above:
\begin{Proposition}\label{prop:untwisted}
Let $c$ be an untwisted $\Spincm$-structure, and $c_\pm$ the $\Spinc$-structures which are its reductions as above.
Then there are identifications among the set of $\Pin^-(2)$-monopole solutions on $c$ and the sets of Seiberg-Witten solutions on $c_\pm$.
Moreover, at the level of moduli spaces, we have
$$
\M_{\Pin^-(2)}(X,c)\cong{\M}_{\U(1)}(X,c_+)\cong{\M}_{\U(1)}(X,c_-),
$$
where ${\M}_{\U(1)}$ means the ordinary Seiberg-Witten {\rm (}$\U(1)$-monopole{\rm )} moduli spaces.
\end{Proposition}
In what follows,  when we use a  phrase like ``a  $\Spinc$( untwisted $\Spincm$)-structure $c$", it means an untwisted $\Spincm$-structure and  its canonical reduction.
We consider them to be an equivalent object, and use them alternatively according to situations. 
%
%
\subsection{Relation with the Seiberg-Witten invariants of the double coverings}\label{subsec:relation}
%
%
Let us consider a twisted $\Spincm$-structure $c$ on a (nontrivial) covering $\pi\colon\X\to X$.
If we pull-back the $\Spincm$-structure $c$ to $\X$, the pulled-back $\Spincm$-structure $\tilde{c}$ on $\X$ is {\it untwisted}.
If $P$ is the $\Spincm(4)$-bundle for $c$, the projection $P\to P/\Spinc\cong \X$ can be considered as a $\Spinc(4)$-bundle over $\X$ which defines a $\Spinc$-structure $\tilde{c}_+$ over $\X$ which is, in fact, the canonical reduction of $\tilde{c}$.
Then $\pi^*P$ is identified with $P\times_{\Spinc(4)}\Spincm(4)$.  
The covering transformation $\iota\colon\X\to\X$ has a natural lift $\tilde{\iota}$ on $\tilde{c}$ which is given by a $\Spincm(4)$-bundle morphism of $P\times_{\Spinc(4)}\Spincm(4)$ defined by $\tilde{\iota}([p,g]) = [pJ,J^{-1}g]$ for $[p,g]\in P\times_{\Spinc(4)}\Spincm(4)$, where $J=[1,j^{-1}]\in \Spincm(4)=\Spin(4)\times_{\deux}\Pin^-(2)$.
Then there is a bijective correspondence between the configuration space of $c$ and the space of $\tilde{\iota}$-invariant configurations on $\tilde{c}$.
If we interpret the objects on $\tilde{c}$ in terms of the $\Spinc$-structure $\tilde{c}_+$,  the $\tilde{\iota}$-action is identified with the antilinear involution $I$ defined in \cite{Pin2}, \S 4(v).
Thus we can identify configurations on $(X,c)$ with $I$-invariant configurations on $(\X,\tilde{c}_+)$. 
In particular, we have,
\begin{Proposition}[\cite{Pin2}, Proposition 4.11]\label{prop:rel}
There is a bijective correspondence between the set of $\Pin^-(2)$-monopole solutions on $(X,c)$ and the set of $I$-invariant Seiberg-Witten solutions on $(\X,\tilde{c}_+)$. 
Moreover we have
\begin{equation}\label{eq:rel}
\M_{\Pin^-(2)}(X,c)\cong{\M}_{\U(1)}(\X,\tilde{c}_+)^I.
\end{equation}
\end{Proposition}
Let us discuss the relation of the $\Pin^-(2)$-monopole invariants of $X$ and the Seiberg-Witten invariants of $\X$.
Mimicking the arguments in \cite{RW} or \cite{Nfree}, 
we can prove a formula which relates the $\Pin^-(2)$-monopole invariants of $(X,c)$ with the Seiberg-Witten invariants of $(\X,\tilde{c}_+)$ as follows.
\begin{Theorem}\label{thm:mod2}
If $d(c)=0$ and $b_1(\X)=0$, then 
\begin{equation}\label{eq:mod2}
\SW^{\U(1)}(\X,\tilde{c}_+)\equiv \sum_{c_\sigma} \SW^{\Pin}(X,c_\sigma) \mod 2
\end{equation}
where $\SW^{\U(1)}(\X,\tilde{c}_+)$ is the Seiberg-Witten invariant of $(\X,\tilde{c}_+)$, and $c_\sigma$ runs through all $\Spin^{c_-}$-structures on $X$ whose pull-back on $\X$ are isomorphic to $\tilde{c}_+$.
\end{Theorem}
\begin{Remark}
Since the $I$-action is free and $d(c)=0$, the virtual dimension of the Seiberg-Witten moduli for $(\X,\tilde{c}_+)$ is also zero.
\end{Remark}
\begin{Remark}
The set of $c_\sigma$'s as above is identified with 
$$
\{c+ a\,|\, a\in\ker(\pi^*\colon H^*(X;l)\to H^*(\X;\pi^*l))\}.
$$ 
\end{Remark}
\proof[Proof of \thmref{thm:mod2}]
In the $I$-equivariant setting, the moduli space $\M_{\U(1)}(\X,\tilde{c}_+)$ is decomposed into the $I$-invariant part and the free part.
The $I$-invariant part is identified with $\M_{\Pin^-(2)}(X,c)$ as in \eqref{eq:rel}.
On the other hand, if the free part is a $0$-dimensional manifold, then the number of elements in the free part is even, because $\Z/2$  acts freely. 
Now, the theorem follows if the equivariant transversality can be achieved by an equivariant perturbation.
This issue is discussed in \cite{Nfree}. ({\it Cf.} \cite{RW}.)
It is easy to achieve the transversality on the free part.
For the $I$-invariant part, on each point $\xi\in\M_{\U(1)}(\X,\tilde{c}_+)^I$,  consider the Kuranishi model $f_\xi\colon H_1\to H_2$, where $H_1$ and $H_2$ are finite dimensional $I$-linear vector spaces.
Since the $I$-action on the base space $\X$ is free, the Lefschetz formula tells us that $H_1$ and $H_2$ are  isomorphic as the $I$-spaces.
Then fixing an $I$-linear isomorphism $L_\xi\colon H_1\to H_2$, we can perturb the equations $I$-equivariantly by using $L_\xi$ to achieve the transversality around $\xi$. 
 \endproof
%
%
Now, we can prove \thmref{thm:Enriques} and \thmref{thm:Enriques-ref}.
\proof[Proof of \thmref{thm:Enriques} and \thmref{thm:Enriques-ref}]
There exists a $\Spin^{c_-}$-structure $c$ on $N$ whose associated $\OO(2)$-bundle is isomorphic to $\underline{\R}\oplus (l_K\otimes\R)$.
Then the associated $\Spin^c$-structure $\tilde{c}$ on the double cover $K$ has a trivial determinant line bundle.
Then $\SW^{\U(1)}(K,\tilde{c})$ is congruent to one modulo $2$ by Morgan-Szab\'{o} \cite{MS}.
On the other hand, since $b_1(N;l)=0$, the Dirac index is even and $d(c)=0$ for the $\Spincm$-structure $c$, the moduli space is orientable, and by fixing an orientation, the $\Z$-valued invariant is defined.
Then, by \thmref{thm:mod2}, there is a $\Spin^{c_-}$-structure $c^\prime$ such that $\SW_\Z^{\Pin}(N, c^\prime)$ is odd.
\endproof
\begin{Remark}
At present, the author does not know the exact value of $\SW_\Z^{\Pin}(N, c^\prime)$ for any homotopy Enriques surface $N$.
\end{Remark}
%
%
\section{Gluing formulae}\label{sec:gluingformulae}
%
%
In this section, we state several versions of gluing formulae for the $\Pin^-(2)$-monopole invariants, and prove \thmref{thm:exotic} and \thmref{thm:Enriques-sum}. 
Before that, we introduce two kinds of $\mu$-maps in order to represent various cohomology classes of $\B^*$.
%
%
\subsection{$\mu$-map (1)}\label{subsec:mu-map1}
%
%
In this subsection, we define the first $\mu$-map, $\mu_{\mathcal{E}}$.
The isomorphism class of a double cover $\X\to X$ is determined by a homomorphism $\rho\colon\pi_1(X)\to \deux$.
Let $H=\pi_1(\X)$. 
When the double cover $\X\to X$ is nontrivial, we have the exact sequence
$$
1\to H \to \pi_1 (X)\overset{\rho}{\to}\deux\to 1.
$$
Let $\iota_*$ be the involution on the rational cohomology group $H_1(\X;\Q)$ induced from the covering transformation $\iota\colon\X\to\X$. 
If we write its $(+1)$(resp. $(-1)$)-eigenspace as $H^+_1$ (resp. $H^-_1$), we have the identifications $H_1^+\cong H_1(X;\Q)$ and $H_1^-\cong H_1(X;l\otimes\Q)$, where $l=\X\times_{\deux}\Z$. 
On the other hand, $H_1(\X;\Q)$ is identified with $(H/[H,H])\otimes\Q$.
Then we can choose loops $\gamma_1,\ldots,\gamma_b$ in $X$, where $b=b_1(X;l)$, such that 
\begin{enumerate}
\item[(C1)] the homotopy class of each $\gamma_i$ is in $\ker\rho$, and 
\item[(C2)] the homology classes of $\gamma_1,\ldots,\gamma_b$ generate $H_1(X;l)/\mathrm{Tor}$.
\end{enumerate}
Note that the restriction of $l$ to $\gamma_i$ is a trivial $\Z$-bundle and the restriction $E_{\gamma_i}$ of $E$ to $\gamma_i$ has a unique $\U(1)$-reduction according the $l$-orientation of $E$.

Let $E$ be the characteristic $\OO(2)$-bundle of a $\Spincm$-structure on a nontrivial double covering $\X\to X$, and $\pi\colon X\times \CC^*\to X$ be the projection. 
We define the universal characteristic $\OO(2)$-bundle $\cal E$ over $X\times \B^*$ as  ${\cal E} = \pi^*E/\G$.
Then we have its characteristic classes 
$$
\tilde{c}_1({\cal E}) \in H^2(X\times\B^*; l\hat{\otimes} \Z),\quad w_2({\cal E})\in  H^2(X\times \B^*;\Z_2),
$$
where $\hat{\otimes}$ denotes the exterior tensor product of local coefficients.
Now let us define the $\mu$-maps 
$$
\hat{\mu}_{\cal E} \colon H_1(X;l)\to H^1(\B^*;\Z), \quad \mu_{\cal E} \colon H_1(X;\Z_2)\to H^1(\B^*;\Z_2),
$$
by the formula
$$
\hat{\mu}_{\cal E} (\alpha) = \tilde{c}_1({\cal E})/\alpha, \quad \mu_{\cal E}(\alpha) = w_2({\cal E})/\alpha.
$$
Since the restriction of $\det E=l\otimes \R$ to $\gamma_i$ is a trivial $\R$-bundle over $\gamma_i$, for any $\OO(2)$-connection $A$ on $E$, the holonomy $\mathrm{Hol}_{\gamma_i}(A)$ around $\gamma_i$  is contained in $\SO(2)\subset\OO(2)$.
Let $\hat{\theta}\in H^1(\SO(2);\Z)$ and $\theta\in H^1(\SO(2);\Z_2)$ be the generators.
\begin{Proposition}\label{prop:hol}
$\hat{\mu}_{\cal E} (\gamma_i) =\mathrm{Hol}_{\gamma_i}^*\hat{\theta}$, ${\mu}_{\cal E} (\gamma_i) =\mathrm{Hol}_{\gamma_i}^*\theta$.
\end{Proposition}
\begin{Remark}
As in the proposition above, we sometimes abuse the symbol for a loop to denote its homotopy class or homology class. 
\end{Remark}
\proof
(The proof is parallel to the ordinary Seiberg-Witten case. {\it Cf.} \cite{OS2}, \S9.)
For a loop $\beta\colon S^1\to \B^*$, the restriction ${\cal E}|_{\gamma_i\times\beta}$ has a $\U(1)$-reduction associated to the $\U(1)$-reduction of  $E|_{\gamma_i}$. 
Then
$$
\langle\tilde{c}_1({\cal E})/\gamma_i,\beta\rangle=\langle c_1({\cal E}|_{\gamma_i\times\beta}), \gamma_i\times\beta\rangle = \deg(\mathrm{Hol}_{\gamma_i}\circ\beta).
$$
\endproof
Since $\B^*\simeq \RP^\infty\times T^b$, $H_1(\B^*;\Z_2)$ and $H_1(\B^*;\Z)$ have decompositions
$$
H_1(\B^*;\Z_2)=H_P\oplus H_T,\quad H_1(\B^*;\Z)=\hat{H}_P\oplus \hat{H}_T,
$$
where $H_P$ is a subgroup isomorphic to $H_1(\RP^\infty;\Z_2)\cong \Z_2$, $\hat{H}_P\cong H_1(\RP^\infty;\Z)\cong\Z_2$, $H_T\cong H_1(T^b;\Z_2)\cong\Z_2^b$ and $\hat{H}_T\cong H_1(T^b;\Z)\cong\Z^b$.
Let $\eta_1$ (resp. $\hat{\eta}_1$) be the generator of $H_P$ (resp. $\hat{H}_P$).
\begin{Corollary}\label{cor:mu1}
There exist basis $\tau_1,\ldots,\tau_b$ for $H_T$ and $\hat{\tau}_1,\ldots,\hat{\tau}_b$ for $\hat{H}_T$ such that
\begin{itemize}
\item $\langle\mu_{\cal E}(\gamma_i),\tau_j\rangle=\delta_{ij},\quad \langle\mu_{\cal E}(\gamma_i),\eta_1\rangle=0$, 
\item $\langle\hat{\mu}_{\cal E}(\gamma_i),\hat{\tau}_j\rangle=\delta_{ij},\quad \langle\hat{\mu}_{\cal E}(\gamma_i),\hat{\eta}_1\rangle=0$.
\end{itemize}
\end{Corollary} 
\proof
The assertions for $\tau_i$ and $\hat{\tau}_i$ are obvious from \propref{prop:hol}.
On the other hand, the class $\eta_1$ is represented by a path $\tilde{\eta}_1=\{(A_t,\Phi_t)\}_{t\in[0,1]}$ in $\CC^*$ such that $A_t=A_0$ and $\Phi_1=-\Phi_0$, and therefore $(A_1,\Phi_1)$ is gauge equivalent to $(A_0,\Phi_0)$ by the constant gauge transformation $-1$.
\endproof
\begin{Remark}\label{rem:jacobian}
For each $\gamma_i$ as above, the holonomy map $\mathrm{Hol}_{\gamma_i}\colon \A/\G\to S^1$ represents a cohomology class $\bar{\gamma}_i$ in $H^1(\A/\G;\Z)\cong [\A/\G,S^1]$.
In fact, $(\bar{\gamma}_1,\ldots,\bar{\gamma}_b)$ gives a basis for $H^1(\A/\G;\Z)$.
\end{Remark}
%
%
\subsection{$\mu$-map (2)}\label{subsec:mu-map2}
%
%
We define the second $\mu$-map $\mu_{\mathcal{F}}$. 
When we define the involution $I$ on $\X\times \C$ by $I(x,v)=(\iota x,\bar{v})$. we have an $\R^2$-bundle $E_0=(\X\times\C)/I$ over $X$ which is identified with $\underline{\R}\oplus(l\otimes\sqrt{-1}\R)$.
Then $\G=\Gamma(\X\times_{\deux}\U(1))$ naturally acts on $E_0$ by $(x,v)\mapsto (x,u(x)v)$.
Let $\pi\colon X\times\CC^*\to X$ be the projection, and define the $\R^2$-bundle $\mathcal{F}$ over $X\times \B^*$ by $\mathcal{F}=\pi^*E_0/\G$.
By using the Stiefel-Whitney class $w_2(\mathcal{F})\in H^2(X\times\B^*;\Z_2)$, define the $\mu$-map  $\mu_{\mathcal{F}}$ for $k=0,1$ as follows:
$$
\mu_{\mathcal{F}}\colon H_k(X;\Z_2)\to H^{2-k}(\B^*;\Z_2),\quad \mu_{\mathcal{F}}(\alpha) = w_2(\mathcal{F})/\alpha.
$$   

Let us consider the case when $\alpha\in H_1(X;\Z_2)$.
By the universal coefficient theorem, we have a split exact sequence
$$
0\to H_1(X;l)\otimes\Z_2\to H_1(X;\Z_2)\to \mathrm{Tor}(H_0(X;l),\Z_2)\to 0.
$$
Then there is a loop $\nu$ in $X$ such that 
\begin{itemize}
\item[(N)] the homology class of $\nu$ corresponds to the generator of $\mathrm{Tor}(H_0(X;l),\Z_2)\cong \Z_2$.
\end{itemize}
Let $\eta_1$ and $\tau_1,\ldots,\tau_b$ be the basis for $H_1(\B^*;\Z_2)=H_P\oplus H_T$ as in \subsecref{subsec:mu-map1}.
\begin{Proposition}\label{prop:muF}
$\langle\mu_{\mathcal{F}}(\nu),\eta_1\rangle =1$, and $\langle\mu_{\mathcal{F}}(\nu),\tau_i\rangle =0$ for any $i$.
\end{Proposition}
\proof
As in the proof of \corref{cor:mu1}, the class $\eta_1$ is represented by a path $\tilde{\eta}_1=\{(A_t,\Phi_t)\}_{t\in[0,1]}$ in $\CC^*$ such that $A_t=A_0$ and $(A_1,\Phi_1)=(-1)(A_0,\Phi_0)$.
Then $\mathcal{F}|_{\nu\times \eta_1}$ is identified with $[0,1]\times[0,1]\times\C/\sim$, where
$$
(0,y,v)\sim (1,y,\bar{v}),\quad (x,0,v)\sim(x,1,-v).
$$
In other words, when $\pi_i\colon S^1\times S^1\to S^1$ is the $i$-th projection and $\varepsilon\to S^1$ is a nontrivial $\R$-bundle over $S^1$,
$$
\mathcal{F}|_{\nu\times \eta_1}\cong \pi_2^*\varepsilon\oplus (\pi_1^*\varepsilon\otimes\pi^*_2\varepsilon).
$$ 
Then the first assertion follows because $w_2(\mathcal{F}|_{\nu\times \eta_1})=w_1(\pi^*_2\varepsilon)w_1(\pi_1^*\varepsilon\otimes\pi^*_2\varepsilon)$ is the generator of $H^2(\nu\times \eta;\Z_2)$.

Recall that $\pi_0\G\cong H^1(X;l)\oplus\Z_2$.
For the dual basis $\check{\gamma}_i\in  H^1(X;l)$ of $\gamma_i\in H_1(X;l)$, we can take $u_i\in\G$ representing $\check{\gamma}_i$.
Then $u_i|_\nu\simeq 1$, and we may assume $u_i|_\nu=1$.
The homology class $\tau_i\in H_1(X;\Z_2)$ is represented by a path $\tilde{\tau}_i = \{(A_t,\Phi_t)\}_{t\in[0,1]}$ such that $\Phi_1=u_i\Phi_0$ and $A_t = A_0 + t(2u_i^{-1}du_i)$.
Then $\mathcal{F}|_{\nu\times\tau_i}$ can be identified with $[0,1]\times[0,1]\times\C/\sim$, where
$$
(0,y,v)\sim (1,y,\bar{v}),\quad (x,0,v)\sim(x,1,v).
$$
Hence  $w_2(\mathcal{F}|_{\nu\times\tau_i})$ is $0$.
\endproof
\begin{Corollary}
$H^*(\B^*;\Z_2)$ is generated by $\mu_{\mathcal{F}}(\nu)$ and $\mu_{\mathcal{E}}(\gamma_i)$ for $i=1,\ldots,b$.
\end{Corollary}
Next we consider $\mu_{\cal F}(x_0)$ for a generator $x_0$ of $H_0(X;\Z_2)$.
\begin{Proposition}
$\mu_{\cal F}(x_0) = \mu_{\cal F}(\nu)\cup \mu_{\cal F}(\nu)$.
\end{Proposition}
\proof
Since $\B^*\simeq \RP^\infty\times T^b$, $H_2(\B^*;\Z_2)$ is generated by 
\begin{itemize}
\item $\eta_2$ corresponding to the generator of $H_2(\RP^\infty;\Z_2)$, 
\item $\eta_1\otimes \tau_j$, where $\eta_1$ and $\tau_j$ are as in \subsecref{subsec:mu-map1}, and 
\item $\tau_i\times \tau_j$ ($i\neq j$).
\end{itemize}

First we prove that $\langle\mu_{\cal F}(x), \eta_2\rangle\neq 0$.
Fix an $\OO(2)$-connection $A_0$ on $E$, and choose $\phi_0,\phi_1,\phi_2\in\Gamma(S^+)$ which are linearly independent.
Let $S$ be the $2$-sphere in $\CC^*$ defined as 
$$
S=\{A_0\}\times\{p\phi_0+q\phi_1+r\phi_2\,|\, p,q,r \in\R,\quad p^2 + q^2+ r^2=1\,\}.
$$
Then the class $\eta_2$ is represented by $[S/\deux]$.
Let $\varepsilon\to \RP^2$ be the canonical line bundle.
We see that ${\cal F}|_{\{x_0\}\times S/\deux}$ is isomorphic to $\varepsilon\oplus \varepsilon$.

Next we prove that $\langle\mu_{\cal F}(x), \eta_1\otimes \tau_i\rangle=\langle\mu_{\cal F}(x), \tau_i\times \tau_j\rangle=0$.
As in the proof of \propref{prop:muF}, we can choose $u_i\in\G$ representing $\check{\gamma}_i$. 
We may assume $u_i(x_0)=1$. 
The homology class $\tau_i\in H_1(X;\Z_2)$ is represented by a path $\tilde{\tau}_i = \{(A_t,\Phi_t)\}_{t\in[0,1]}$ such that $\Phi_1=u_i\Phi_0$ and $A_t = A_0 + t(2u_i^{-1}du_i)$.
Then we can see that 
$$
w_2({\cal F}|_{\{x_0\}\times (\eta_1\times \tau_i)})=w_2({\cal F}|_{\{x_0\}\times (\tau_i\times \tau_j)}) = 0.
$$
\endproof
For cohomology classes of $\B^*$, let
$$
\nu^*=\mu_{\mathcal{F}}(\nu),\quad \gamma_i^* = \mu_{\mathcal{E}}(\gamma_i),\quad \hat{\gamma}_i^* = \hat{\mu}_{\mathcal{E}}(\gamma_i).
$$
Then, for example, $H^*(\B;\Z_2)$ can be written as
$$
H^*(\B;\Z_2) = \Z_2[\nu^*]\otimes\bigwedge (\Z_2\gamma^*_1\oplus\cdots\oplus\Z_2\gamma^*_b),
$$
and a cohomology class $\xi\in H^*(\B;\Z_2)$ can be written as 
$$
\xi = (\nu^*)^a\prod_{i\in I}\gamma_i^*,
$$
where $a$ is a non-negative integer and $I$ is a subset of $\{1,\ldots,b\}$.

For a $\Spinc$(untwisted $\Spincm$)-structure, we have the $\mu$-map of ordinary Seiberg-Witten theory (\cite{OS2}, \S9):
$$
\mu_{0}\colon H_k(X;\Z)\to H^{2-k}(\B^*;\Z)\quad (k=0,1).
$$
For $x\in H_0(X;\Z)$ and $\gamma\in H_1(X;\Z)$, let
$x^* = \mu_{0}(x)$, $\gamma^*=\mu_0(\gamma).$
%
%
\subsection{Cutting down the moduli spaces}\label{subsec:cutdown}
%
%
The purpose of this subsection is to construct the submanifolds in the moduli spaces which are dual to the classes $\mu_{\mathcal{F}}(\nu)$, $\mu_{\mathcal{F}}(x_0)$ and $\mu_{\mathcal{E}}(\gamma_i)$. 
({\it Cf}. \cite{DK}, \S5.2 and \cite{OS}, \S9.)
For a loop $\nu$ in $X$ as in \subsecref{subsec:mu-map2}, fix a tubular neighborhood $n(\nu)$ of $\nu$ which is a smooth open submanifold with smooth boundary in $X$.
Let $\CC_{n(\nu)}^*$ be the space of irreducible configurations on $n(\nu)$, $\G_{n(\nu)}$ be the gauge transformation group and $\B_{n(\nu)}^*=\CC_{n(\nu)}^*/\G_{n(\nu)}$.
Note that $\pi_0\G_{n(\nu)}=\deux$.
Let $\G_{n(\nu)}$ act on $\R$ via the projection $\G_{n(\nu)}\to \pi_0\G_{n(\nu)}=\deux$ and the multiplication of $\deux$. 
Dividing by the diagonal action, we obtain a real line bundle 
$$
{\cal L}_\nu = \CC_{n(\nu)}^*\times_{\G_{n(\nu)}}\R \to \B^*_{n(\nu)}.
$$
Suppose thet the moduli space  $\M(X)$ contains no reducibles and is  perturbed to be a smooth manifold.
Let $M$ be $\M(X)$ itself or its smooth submanifold.
Since the restriction of an irreducible solution on $X$ to an open subset of $X$ is also irreducible by the unique continuation property of the Dirac operator, we have a well-defined restriction map
$$
r_\nu\colon M\to \B_{n(\nu)}^*.
$$
We can choose a section $s$ of $\mathcal{L}_\nu$ so that the pull-back $r^*_\nu s$ is transverse to the zero-section of $r_\nu^*\mathcal {L}_\nu$ (\cite{DK}, 5.2.2).
Then the zero-set of $r_\nu^* s$ is a codimension-one submanifold of $M$ which is dual to the class $\mu_{\mathcal{F}}(\nu)$ in $M$, and is denoted by 
$$
M\cap V_\nu.
$$
Similarly, for the class $\mu_{\mathcal {F}}(x_0)$, we can construct a codimension-two submanifold of $M$ which is dual to $\mu_{\mathcal {F}}(x_0)$ in $M$, and is denoted by
$$
M\cap V_{x_0}.
$$
For the loops $\gamma_i$ chosen in \subsecref{subsec:mu-map1}, let $\mathrm{Hol}_{\gamma_i}\colon M\to S^1$ be the smooth map defined by the holonomy around $\gamma_i$.
When we take a regular value $\theta\in S^1$ of $\mathrm{Hol}_{\gamma_i}$, the inverse image $\mathrm{Hol}_{\gamma_i}^{-1}(\theta)$ is a codimension-one submanifold of $M$ which is dual to $\mu_{\mathcal{E}}(\gamma_i)$ in $M$, and is denoted by 
$$
M\cap V_{\gamma_i}.
$$

%
%
\subsection{Gluing theorems}\label{subsec:GF}
%
%
In this subsection, we state several gluing formulae for $\Pin^-(2)$-monopole invariants, which will be proved in later sections.
The formulae have different forms depending on whether the $\Spincm$-structures are twisted or untwisted, and the moduli spaces contain reducibles or not.
For local coefficients $l_1$ and $l_2$ over $X_1$ and $X_2$, if both of $l_i$ are nontrivial, then we have $b_1(X_1\#X_2;l_1\#l_2)=b_1(X_1;l_1) + b_1(X_2;l_2) +1$ by the Meyer-Vietoris sequence. 
Hence there is an extra generator of $H_1(X_1\#X_2)$ which does not come from $X_1$ and $X_2$.
On the other hand, if one of $l_i$ is trivial, then $b_1(X_1\#X_2;l_1\#l_2) = b_1(X_1;l_1) +  b_1(X_2;l_2)$.
Choose loops $\alpha_1,\ldots,\alpha_{b_1(l_1)}$ in $X_1$, and $\beta_1,\ldots,\beta_{b_1(l_2)}$ in $X_2$, where $b_1(l_i)=b_1(X_i;l_i)$ for $i=1,2$,  and $\delta$ in $X_1\#X_2$ representing an extra generator if both of $l_1$ and $l_2$ are nontrivial, such that
\begin{itemize}
\item $\alpha_1,\ldots,\alpha_{b_1(l_1)}$ and $\beta_1,\ldots,\beta_{b_1(l_2)}$ satisfy the conditions (C1) and (C2) in \subsecref{subsec:mu-map1} for $(X_1,l_1)$ and $(X_2,l_2)$, respectively, and
\item $\alpha_1,\ldots,\alpha_{b_1(l_1)}, \beta_1,\ldots,\beta_{b_1(l_2)}$ and $\delta$ (if exists) satisfy the conditions (C1) and (C2) for $(X_1\#X_2,l_1\#l_2)$.
(We assume $\alpha_i$ and $\beta_j$ are also contained in $X_1\#X_2$.)
\end{itemize}
For each $i=1,2$, if $l_i$ is nontrivial, then choose another loop $\nu_i$ in $X_i$ satisfying the condition (N) before \propref{prop:muF}.
We also assume that $\nu_i$ is contained in $X_1\#X_2$.

The first gluing formula is on the gluing of $\U(1)$-irreducible monopoles and $\Pin^-(2)$-reducible monopoles.
\begin{Theorem}\label{thm:GF}
Let $X_1$ be a closed oriented connected $4$-manifold with $b_+(X_1)\geq 2$ and a $\Spinc${\rm (}untwisted $\Spincm${\rm )}-structure $c_1$. 
Let $X_2$ be a closed oriented connected $4$-manifold which satisfies the following{\rm :}
\begin{itemize}
\item There exists a nontrivial double covering ${\X}_2\to X_2$ with $b_+(X_2;l_2)=0$ where $l_2={\X}_2\times_{\deux}\Z$.
\item There exists a $\Spincm$-structure $c_2$ on ${\X}_2\to X_2$ such that $\tilde{c}_1(E)^2=\sign(X_2)$ {\rm (}and hence the Dirac index is $0$ and $d(c_2)=b_1(X_2;l_2)${\rm )}.
\end{itemize}
For a cohomology class $\xi\in H^*(\B^*(X_1,c_1);\Z_2)$ of the form $\xi=\prod_{i\in I}\mu_0(\alpha_i)$ where $I \subset \{1,\ldots,b_1(l_1)\}$, let $\xi^\prime= \prod_{i\in I}\mu_{\cal E}(\alpha_i)\in H^*(\B^*(X_1\# X_2,c_1\#c_2);\Z_2)$. 
Then we have
$$
\SW^{\Pin}(X_1\#X_2,c_1\#c_2)(\xi^\prime(\nu_2^*)^{2a+1}\beta_1^*\cdots\beta_{b_1(l_2)}^*) \equiv \SW^{\U(1)}(X_1,c_1)(\xi (x^*)^a)\mod 2.
$$
\end{Theorem}
\thmref{thm:GF1} is a corollary of \thmref{thm:GF}.

The second one is a generalized blow-up formula by the gluing of $\Pin^-(2)$-irreducibles and $\U(1)$-reducibles.
\begin{Theorem}[{\it Cf.} \cite{FS,Nicolaescu, FC}]\label{thm:blowup}
Let $X_1$ be a closed oriented connected $4$-manifold with a $\Spincm$-structure $c_1$ with $b_+(X_1;l_1)\geq 2$.
Let $X_2$ be a closed oriented connected $4$-manifold with  a $\Spinc${\rm (}untwisted $\Spincm${\rm)}-structure $c_2$ such that $b_1(X_2)=b_+(X_2)=0$ and $d(c_2) =-1$.
For any $\xi=(\nu_1^*)^a\prod_{i\in I}\alpha^*_{i}$ where $I \subset \{1,\ldots,b_1(l_1)\}$, 
$$
\SW^{\Pin}(X_1\# X_2, c_1\#c_2)(\xi) = \SW^{\Pin}(X_1,c_1)(\xi).
$$ 
\end{Theorem}
\begin{Remark}
In \thmref{thm:blowup}, $\xi$ is assumed to represent  both of the cohomology  classes  of $\B^*(X_1,c_1)$ and $\B^*(X_1\#X_2,c_1\#c_2)$.
The similar remark is valid for the following theorems.
\end{Remark}
The third one is on the gluing of $\Pin^-(2)$-irreducibles and $\Pin^-(2)$-reducibles.
\begin{Theorem}\label{thm:twistedGF}
Let $X_1$ be a closed oriented connected $4$-manifold with a twisted $\Spincm$-structure $c_1$ with $b_+(X;l_1)\geq 2$, and $X_2$ be a manifold with a $\Spincm$-structure $c_2$ as in \thmref{thm:GF}.
Then, for any $\xi=(\nu_1^*)^a\prod_{i\in I}\alpha^*_{i}$ where $I \subset \{1,\ldots,b_1(l_1)\}$, 
$$
\SW^{\Pin}(X_1\# X_2,c_1\#c_2)(\xi\delta^*\beta_1^*\cdots\beta_{b_1(l_2)}^*) = \SW^{\Pin}(X_1,c_1)(\xi).
$$
\end{Theorem}

If $4$-manifolds $X_1$ and $X_2$ have positive $b_+$, then the Seiberg-Witten invariants of $X_1\# X_2$ are always $0$. 
Likewise, the $\Z_2$-valued $\Pin^-(2)$-monopole invariants have a similar property.
\begin{Theorem}\label{thm:vanishing}
Let $X_1$ be a closed oriented connected $4$-manifold with a twisted $\Spincm$-structure $c_1$ with $b_+(X_1;l_1)\geq 1$.
Let $X_2$ be a closed oriented connected $4$-manifold with a {\rm (}twisted or untwisted{\rm )} $\Spincm$-structure $c_2$, and suppose one of the following:
\begin{enumerate}
\item[(\romnum1)] $b_+(X_2)\geq 1$ and $c_2$ is an untwisted $\Spincm$-structure on $X_2$.
\item[(\romnum2)] $c_2$ is a twisted $\Spincm$-structure on $X_2$ with $b_+(X_2;l_2)\geq  1$.
\end{enumerate}
Then $\SW^{\Pin} (X_1\# X_2,c_1\# c_2)(\xi) = 0$ for any class $\xi\in H^*(\B;\Z_2)$.
\end{Theorem}

On the other hand, the $\Z$-valued invariants can be nontrivial for a connected sum $X_1\#X_2$ even when both of $b_+(X_1;l_1)$ and $b_+(X_2;l_2)$ are positive.
Consider $(X_i,l_i)$ $(i=0,1,\ldots, n)$ with nontrivial $l_i$. 
We assume $b_1(X_i,l_i)=0$ for every $i$.
As noticed above, each time we take a connected sum of these, we have an extra generator in the first homology of the connected sum.
Choose loops $\delta_1,\ldots,\delta_n$ in $X_0\#\cdots\# X_n$ representing such extra generators in $H_1( X_0\#\cdots \# X_n; l_1\#\cdots\# l_n)$ satisfying the conditions (C1), (C2). 
\begin{Theorem}\label{thm:ZGF}
Let $n$ be any positive integer.
For $i=0,1,\ldots,n$, let $X_i$ be a closed oriented connected $4$-manifold with a twisted $\Spincm$-structure $c_i$ satisfying 
\begin{itemize}
\item $b_1(X_i;l_i)=0$, $b_+(X_i;l_i)\geq 2$.
\item $d(c_i)=0$, and
\item the index of the Dirac operator is positive and even.
\end{itemize}
Note that in this situation, the moduli space $\M(X_i,c_i)$ is orientable, and the $\Z$-valued invariant $\SW_\Z^{\Pin}(X_i,c_i)(1)$ is defined for a choice of orientation.
Then the glued moduli space $\M(X_0\#\cdots\# X_n, c_1\#\cdots\# c_n)$ is orientable, and
$$
\SW_\Z^{\Pin}(X_0\#\cdots\# X_n, c_1\#\cdots\# c_n)(\hat{\delta}_{1}^*\cdots\hat{\delta}_n^*) = 2^n\prod_{i=0}^n \SW_\Z^{\Pin}(X_i,c_i)(1),
$$ 
for a choice of orientation.
\end{Theorem}
%
%
%
%
\subsection{Proofs of \thmref{thm:exotic} and \thmref{thm:Enriques-sum}}
%
%
In this subsection, we prove \thmref{thm:exotic} and \thmref{thm:Enriques-sum} by assuming \thmref{thm:GF} and \thmref{thm:ZGF}.

\proof[Proof of \thmref{thm:exotic}]
Let $(X_2,l_{X_2})$ be as in \thmref{thm:GF1}.
Then this satisfies the conditions for $X_2$ in \thmref{thm:GF}.
 
For given $n$, required exotic structures on $E(n)$ can be constructed by either logarithmic transformation (see e.g., \cite{GS}) or Fintushel-Stern's knot surgery \cite{FS2}.

First, we discuss on the case of logarithmic transformation.
Let $E(n)_{p.q}$ be the log transformed $E(n)$ with two multiple fibers of multiplicities $p$ and $q$.
For odd $n$, all of $E(n)_{p.q}$ with $\gcd(p,q)=1$ is homeomorphic to $E(n)$.
On the other hand, for even $n$, $E(n)_{p.q}$ is homeomorphic to $E(n)$ if and only if $\gcd(p,q)=1$ and $pq$ is odd.
Let $f\in H^2(E(n)_{p,q})$ be the Poincar\'{e} dual of the homology class of a regular fiber.   
Then there is a primitive class $f_0$ with $f=pqf_0$, and the Poincar\'{e} duals $f_p$ and $f_q$ of the multiple fibers of $p$ and $q$ are given by $f_p=qf_0$ and $f_q = pf_0$.
If we put
$$
D(a,b,c)= a f + bf_p + cf_q, 
$$
then, for $n\geq 2$, the canonical class $K$ is given as $K=D(n-2,p-1,q-1)$.
The Seiberg-Witten basic classes are given by $K-2D(a,b,c)$, where $0\leq a\leq n-2$, $0\leq b \leq p-1$, $0\leq c\leq q-1$, and the value the Seiberg-Witten invariant for the class $K-2D(a,b,c)$ is 
$$
\SW^{\U(1)}(E(n)_{p,q}, K-2D(a,b,c))=(-1)^a
\begin{pmatrix}
n-2\\
a
\end{pmatrix},
$$ 
which is independent of $b$ and $c$.
Similar facts hold for the case when $n=1$.
In general, the number of basic classes whose Seiberg-Witten invariants are odd is changed if $p$ and $q$ are varied. 
By using these facts together with \thmref{thm:GF}, we can find infinitely many $\{p,q\}$ such that $E(n)_{p,q}\#X_2$ have different numbers of basic classes for $\Pin^-(2)$-monopole invariants.

For a knot $K$, let $E(n)_K$ be the manifold obtained by the knot surgery on a regular fiber $T$ with $K$.
If we consider the Seiberg-Witten invariant as a symmetric Laurent polynomial as in \cite{FS2}, the invariant of $E(n)$ is related to that of  $E(n)_K$  by 
$$
\SW^{\U(1)}_{E(n)_K} = \SW^{\U(1)}_{E(n)}\cdot\Delta_K(t), 
$$
where $t=\exp(2[T])$ and $\Delta_K(t)$ is the (symmetrized) Alexander polynomial of $K$.
Now, let $X_K=E(n)_K$, and let us fix a $\Spincm$-structure $c_2$ on $X_2$ as in \thmref{thm:GF}, and consider a function of $\Pin^-(2)$-monopole invariants of $X_K\# X_2$, 
$$
\SW^{\Pin}_{X_K\# (X_2, c_2)} \colon \{ h\in H^2(X_K;\Z)\, |\, h\equiv w_2(X)\mod 2 \} \to \Z_2,
$$
which is defined as 
$$
\SW^{\Pin}_{X_K\# (X_2, c_2)} (h) =\SW^{\Pin} (X_K\# X_2,c(h) \# c_2)(\nu_2^*\beta_1^*\cdots\beta_{b_1(l_2)}^*),
$$
where $c(h)$ is the $\Spinc$-structure on $X_K$ with $c_1=h$.
If we assume $\SW^{\Pin}_{X_K\# (X_2, c_2)} $ as a $\Z_2$-coefficient polynomial, then \thmref{thm:GF1} implies that $\SW^{\Pin}_{X_K\# (X_2, c_2)} $ is the $\Z_2$-reduction of the  $\Z$-coefficient polynomial $\SW^{\U(1)}_{E(n)_K}$. 
Then we can find infinitely many $K$ so that $\SW^{\Pin}_{X_K\# (X_2, c_2)} $ are different.
\endproof

\proof[Proof of \thmref{thm:Enriques-sum}]
For each $(N_i,l_i)$, we have $b_1(X_i;l_i)=0$ and $b_+(X_i;l_i)=2$. 
By \thmref{thm:Enriques-ref}, there is a twisted $\Spincm$-structure $c_i$ such that $d(c_i)=0$, the Dirac index is $2$ and
$\SW_\Z^{\Pin}(X_i,c_i)$ is odd. 
Then the theorem follows from \thmref{thm:ZGF}.
\endproof

%
%
\section{$\Pin^-(2)$-monopole theory on $3$-manifolds}
%
%
Sections 4--6 are devoted to the proof of the gluing theorems in \subsecref{subsec:GF}, and this preparatory section is on the $\Pin^-(2)$-monopole theory on $3$-manifolds.
We refer to \cite{KM2, FC} for the Seiberg-Witten counterpart of the topics in this section.
%
%
\subsection{$\Spin^{c_-}$-structures on $3$-manifolds}
%
%
Define the group $\Spin^{c_-}(3)$ by
$$
\Spin^{c_-}(3) = \Spin(3)\times_{\{\pm 1\}}\Pin^-(2) = \SP(1)\times_{\{\pm 1\}}\Pin^-(2).
$$
Let $Y$ be an oriented closed connected Riemannian $3$-manifold, and $Fr(Y)$ its $\SO(3)$-frame bundle.
Suppose a double covering $\tilde{Y}\to Y$ is given. 
A $\Spin^{c_-}$-structure on $\tilde{Y}\to Y$ consists of  a principal $\Spin^{c_-}(3)$-bundle $P$ and isomorphisms  $\sigma\colon P/\Spinc(3)\to\tilde{Y}$ and $\tau\colon P/\Pin^-(2)\to Fr(Y)$.
The characteristic $\OO(2)$-bundle $E$ is defined as $E=P/\Spin(3)$.
\begin{Remark}
As in the $4$-dimensional case, if $\tilde{Y}\to Y$ is trivial, then a $\Spincm$-structure on $\tilde{Y}\to Y$ can be reduced to a $\Spinc$-structure on $Y$, and is called {\it untwisted}.
\end{Remark}
Define the action of $\Spin^{c_-}(3)$ on $\Ima \HH$ by 
$$
[q,u]\cdot v =qvq^{-1},  
$$
for $[q,u]\in \Spin^{c_-}(3)$ and $v\in\Ima \HH$. 
Then the associated bundle $P\times_{\Spincm(3)}\Ima\HH$ is identified with the tangent bundle $TY$.
Define the $\Spin^{c_-}(3)$-action on $\HH$ by 
$$
[q,u]\cdot\psi = q\psi u^{-1},
$$
for $[q,u]\in \Spin^{c_-}(3)$ and $\psi\in\HH$.
Then we obtain the associated bundle $S=P\times_{\Spincm(3)}\HH$ which is the spinor bundle for the $\Spin^{c_-}$-structure.

The Clifford multiplication is defined as follows.
The identity component of $\Spin^{c_-}(3)$ is a $\Spin^{c}(3)$, and the quotient group $\Spin^{c_-}(3)/\Spin^{c}(3)$ is isomorphic to $\{\pm 1\}$.
Let $\C_-$ be a copy of $\C$ with the $\{\pm 1\}$-action by complex conjugation.
Then $\Spin^{c_-}(3)$ acts on $\C_-$ via the projection $\Spin^{c_-}(3)\to\Spin^{c_-}(3)/\Spin^{c}(3)=\{\pm 1\}$.
If we define 
$$
\rho_0\colon (\Ima \HH)\otimes_\R\C_-\times\HH\to \HH
$$
by $\rho_0(v\otimes a,\psi)=\bar{v}\psi\bar{a}$, then $\rho_0$ is $\Spin^{c_-}(3)$-equivariant.
Let $K=\tilde{Y}\times_{\{\pm 1\}}\C_-$. 
Then we can define the Clifford multiplication
$$
\rho\colon T^*Y\otimes_\R K \to \Hom (S,S),
$$
which induces
$$
\rho\colon\Omega^1(Y;K)\times\Gamma(S)\to\Gamma(S).
$$
Note that $K=\underline{\R}\oplus i\lambda$, and so $\Omega^1(Y;K)=\Omega^1(Y;\underline{\R})\oplus\Omega^1(Y;i\lambda)$.
Although the spinor bundle $S$ does not have an ordinary hermitian inner product, the pointwise {\it twisted} hermitian product 
\begin{equation}\label{eq:hermitien}
\langle\cdot,\cdot\rangle_{K,x} \colon S_x\times S_x\to K_x
\end{equation}
is defined. 
For $\alpha\otimes 1\in T^*Y\otimes K$, the image $\rho(\alpha\otimes 1)$ is a traceless endomorphism which is {\it skew-adjoint} with respect to the inner product \eqref{eq:hermitien}.
The whole image of $T^*Y$ by $\rho$ forms the subbundle of $\Hom (S,S)$, which we write as $\tilde{\su}(S)$, equipped with the inner product $\frac12\tr(a^*b)$.
When $\{e_1$, $e_2$, $e_3\}$ is an oriented frame on $\Lambda^1(Y)$, we assume the orientation convention
$$
\rho(e_1)\rho(e_2)\rho(e_3)=1.
$$ 
We extends $\rho$ to forms by the rule,
$$
\rho(\alpha\wedge\beta)=\frac12 (\rho(\alpha)\rho(\beta)+(-1)^{\deg\alpha\deg\beta}\rho(\beta)\rho(\alpha)).
$$
The orientation convention implies $\rho(*\alpha)=-\rho(\alpha)$ for $1$-forms.
%
%
\subsection{$\Pin^{-}(2)$-monopole equations on $3$-manifolds}
%
%
An $\OO(2)$-connection $B$ on $E$ together with the Levi-Civita connection defines a $\Spin^{c_-}(3)$-connection on $P$.
Then we have the Dirac operator $D_B\colon\Gamma(S)\to\Gamma(S)$ associated to $B$.

The bundle $\Lambda^1(Y)\otimes_\R i\lambda$ is also associated with $P$ as follows.
Let $\varepsilon\colon\Pin^-(2)\to\Pin^-(2)/\U(1)\cong\{\pm 1\}$ be the projection, and let $\Spin^{c_-}(3)$ act on $\Ima \HH$ by 
$$
v\in\Ima \HH \to \varepsilon(u)qvq^{-1}\quad\text{for }[q,u]\in\Spin^{c_-}(3).
$$
Then $\Lambda^1(Y)\otimes_\R i\lambda$ is identified with $P\times_{\Spin^{c_-}(3)}\Ima \HH$.
For $\psi\in\HH$, $\psi i\bar{\psi}$ is in $\Ima \HH$. 
Then the map $\psi\in\HH \to \psi i\bar{\psi}\in\Ima \HH$ is $\Spin^{c_-}(3)$-equivariant, and induces a quadratic map 
$$
q\colon\Gamma(S)\to \Omega^1(Y;i\lambda).
$$
For a closed $2$-form $\eta\in\Omega^2(i\lambda)$, the perturbed $\Pin^-(2)$-monopole equations on $Y$ are defined as 
\begin{equation}\label{eq:Pin2eq3}
\left\{
\begin{aligned}
D_B\Psi&=0,\\
-\frac12(*(F_B&+\eta))=q(\Psi),
\end{aligned}\right.
\end{equation}
for  $\OO(2)$-connections $B$ on $E$ and $\Psi\in \Gamma(S)$. 
The gauge transformation group is given by 
$$
\G_Y = \Gamma (\tilde{Y}\times_{\{\pm 1\}}\U(1)),
$$
where $\{\pm 1\}$ acts on $\U(1)$ by complex conjugation.
\begin{Remark}
If the $\Spincm$-structure is untwisted, then the $3$-dimensional $\Pin^-(2)$-monopole equations are also identified with the $3$-dimensional Seiberg-Witten equations.
\end{Remark}
%
%
\subsection{$\Pin^{-}(2)$-Chern-Simons-Dirac functional}
%
%
Choose a reference $\OO(2)$-connection $B_0$ on $E$.
Let $\A(E)$ be the space of $\OO(2)$-connections on $E$, and $\CC=\A(E)\times\Gamma(S)$.
\begin{Definition} Let $\eta$ be a closed $2$-form in $\Omega^2(\lambda)$.
The (perturbed) $\Pin^{-}(2)$-Chern-Simons-Dirac functional $\vartheta\colon\CC\to\R$ is defined by
\begin{equation}\label{eq:CSD}
\vartheta(B,\Psi)=-\frac18\int_Y(B-B_0)\wedge (F_B+F_{B_0}+i\eta)+\frac12\int_Y\langle D_B\Psi,\Psi\rangle_\R{\rm dvol}_Y.
\end{equation}
\end{Definition}
A few comments on the definition.
For $\alpha\in\Omega^1(i\lambda)$ and $\beta\in\Omega^2(i\lambda)$, $\alpha\wedge \beta$ is in $ \Omega^3(Y;\R)$ since $\lambda^{\otimes 2}$ is trivial.
The inner product $\langle\cdot,\cdot\rangle_\R$ is the real part of \eqref{eq:hermitien}.

The tangent space of $\CC$ at $(B,\Psi)$ is $T_{(B,\Psi)}\CC = \Omega^1(i\lambda)\oplus \Gamma(S)$. 
We equip the tangent space with an $L^2$ metric.
Then the gradient of $\vartheta$ with respect to the $L^2$-metric is given by
$$
\nabla\vartheta = \left(\frac12(*(F_B+i\eta)) + q(\Psi), D_B\Psi\right).
$$
Hence the critical points of $\vartheta$ are the solutions of the $\Pin^-(2)$-monopole equations on $Y$. 

For a critical point $(B,\Psi)$ of $\vartheta$, let $\mathscr{H}_{(B,\Psi)}\colon \Omega^1(i\lambda)\oplus\Gamma(S)\to\Omega^1(i\lambda)\oplus\Gamma(S)$ be the derivative of $\nabla\vartheta$ at $(B,\Psi)$ given as
$$
\mathscr{H}_{(B,\Psi)}(b,\psi)=\left(\frac12*db-Dq_\Psi(\psi),-D_B\psi-\frac12b\Psi\right),
$$
where $Dq_\Psi$ is the linearization of $q$.
A critical point $(B,\Psi)$ is called {\it non-degenerate} if the middle cohomology group of the following complex is $0$:
$$
\begin{CD}
\Omega^0(i\lambda)@>{\mathcal{I}_\Psi}>>\Omega^1(i\lambda)\oplus\Gamma(S)@>{\mathscr{H}_{(B,\Psi)}}>>\Omega^1(i\lambda)\oplus\Gamma(S),
\end{CD}
$$
where $\mathcal{I}_\Psi$ is defined by $\mathcal{I}_\Psi(f)=(-2df,f\Psi)$.

For $g\in\G_Y$, $g^{-1}dg$ is an $i\lambda$-valued $1$-form, and the $\lambda$-valued $1$-form $\frac1{2\pi i}g^{-1}dg$ represents an integral class $[g]\in H^1(Y;l)/\mathrm{Tor}$.
\begin{Proposition}
For $(B,\Psi)\in\CC$ and $g\in\G_Y$, 
$$
\vartheta(g(B,\Psi))-\vartheta(B,\Psi) = 2\pi([g]\cup (\pi \tilde{c}_1(E)-[\eta])[Y],
$$
where $[\eta]\in H^2(Y;\lambda)$ is the de Rham cohomology class of $\eta$. 
\end{Proposition}
%
%
\subsection{Non-degenerate critical point on $S^3$}
%
%
Here, we suppose $Y=S^3$ with a positive scalar curvature metric. 
Since $S^3$ is simply-connected, every $\Spincm$-structure is untwisted.
This is unique up to isomorphism and identified with a unique $\Spinc$-structure.
For a positive scalar curvature metric, every monopole solution is a reducible one, say $(\theta,0)$, which is unique up to gauge. 
Furthermore, the kernel of the Dirac operator $D_\theta$ is trivial.
Since the index of $D_\theta$ is $0$, the cokernel is also trivial, and this implies $(\theta,0)$ is nondegerate.
The stabilizer of $(\theta,0)$ of the gauge group action is denoted by $\Gamma_\theta$:
$$
\Gamma_\theta =\{g\in \Map(S^3;\U(1))\,|\, g(\theta,0)=(\theta,0)\}.
$$
Note that $\Gamma_\theta\cong S^1$.
%
%
\section{$\Pin^-(2)$-monopoles on a $4$-manifold with a tubular end}
%
%
In this section, we continue the preparation for gluing, and discuss on finite energy $\Pin^-(2)$-monopoles on  $4$-manifolds with  tubular ends.
We refer to \cite{DF} as well as \cite{KM2, FC}.
%
%
\subsection{Setting}
%
%
Let $X$ be a Riemannian $4$-manifold with a $\Spin^{c_-}$-structure containing a tubular end $ [-1,\infty)\times Y$, where $Y$ is a closed, connected, Riemannian $3$-manifold with a $\Spin^{c_-}$-structure.
More precisely, suppose we are given 
\begin{enumerate}
\item an orientation preserving isometric embedding $i\colon [-1,\infty)\times Y\to X$ such that 
$$
X^t = X\setminus i( (t,\infty)\times Y) 
$$
is compact for any $t\geq -1$,
\item an isomorphism between $\Spin^{c_-}$-structure on $[-1,\infty)\times Y$ induced from $Y$ and the one inherited from $X$ via the embedding $i$.
\end{enumerate}
\begin{Remark}
If the $\Spincm$-structure on $X$ is twisted but its restriction on the tube $[-1,\infty)\times Y$ is untwisted, then the double cover $\X$ has two tubular ends.
\end{Remark}
In order to  define weighted Sobolev norms on various sections over $X$, take a $C^\infty$-function $w\colon X\to\R$ such that 
\begin{equation}\label{eq:w}
w(t)=\left\{
\begin{aligned}
1 \ \ & \quad \text{on }X^{-1}\\
e^{\alpha t} &\quad \text{for }(t,y)\in [0,\infty)\times Y
\end{aligned}\right. 
\end{equation}
where $\alpha$ is a small positive number which will been chosen later to be suitable for our purpose.
For  a nonnegative integer $k$, we will use the weighted Sobolev norm of a section $f$ (e.g., a form or a spinor) on $X$ given by 
$$
\|f\|_{L^{2,w}_k} = \|wf \|_{L^2_k}.
$$

Let $X_1$ and $X_2$ be $4$-manifolds with tubular ends as above with isometric embeddings 
$$
i_1\colon[-1,\infty)\times Y\to X_1, \quad i_2\colon[-1,\infty)\times \bar{Y}\to X_2,
$$
where $\bar{Y}$ is $Y$ with opposite orientation.
For $T\geq 0$, let $X^{\#T}$ be the manifold obtained by gluing $X^{2T}_1$ and $X^{2T}_2$ via  the identification
$$
i_1(t,y)\sim i_2(2T-t,y).
$$
Then we naturally have an isometric embedding of a neck $i_T\colon [-T,T]\times Y\to X^{\#T}$.
(Here, the negative side is connected to $X_1^0$ and the positive side to $X_2^0$.)
When we take functions $w_1$, $w_2$ as \eqref{eq:w}, a continuous function $w_T \colon X^{\#T}\to \R$ is induced by gluing $w_1$ and $w_2$ such that  
\begin{equation}\label{eq:wT}
w_T(t) = e^{\alpha (T-|t|)}
\end{equation}
for $(t,y)\in [-T,T]\times Y$.
For the sections over $X^{\#T}$, we will use the weighted norm 
$$
\|f \|_{L^{2,w_T}_k} = \| w_T f\|_{L^2_k}.
$$

%
%
\subsection{Exponential decay}\label{subsec:expdecay}
%
%
The purpose of this subsection is to give exponential decay estimates for $\Pin^-(2)$-monopoles on a cylinder $[0,\infty) \times Y$ and a band $(-T,T)\times Y$.
Since a $\Pin^-(2)$-monopole on an untwisted $\Spincm$-structure is identified with  an ordinary Seiberg-Witten  monopole, the estimates for Seiberg-Witten monopoles on a cylinder $[0,\infty) \times Y$ hold for $\Pin^-(2)$-monopoles on an untwisted $\Spincm$-structure. 
On the other hand, we can also obtain an estimate for $\Pin^-(2)$-monopoles on a twisted $\Spincm$-structure by lifting everything to the double cover $[0,\infty)\times \tilde{Y}$ on which the corresponding $\Spincm$-structure is untwisted and applying the estimate for the Seiberg-Witten monopole.  
Thus,  invoking the results due to Froyshov \cite{FC} for the Seiberg-Witten monopoles, we  obtain the estimates for $\Pin^-(2)$-monopoles as follows.

Let $\beta$ be a nondegenerate monopole over $Y$, and $U\subset \B_Y$ is an $L^2$-closed
subset which contains no monopoles except perhaps $[\beta]$.
Define $B_t= [t-1,t+1]\times Y$. 
\begin{Theorem}[\cite{FC}, Theorem 6.3.1.]\label{thm:expdecay}
There exists a constant $\lambda_+$ which has the following significance.
For any $C>0$, there exist constants $\epsilon$ and $C_k$ for nonnegative integer $k$ such that the following holds.
Let $x=(A,\Phi)$ be a $\Pin^-(2)$-monopole in temporal gauge over $(-2,\infty)\times Y$ such that $x(t)\in U$ for some $t\geq 0$.
Set 
$$
\bar{\nu}=\|\nabla\vartheta\|_{L^2((-2,\infty)\times Y)},\quad 
\nu(t)=\|\nabla\vartheta\|_{L^2(B_t)}.
$$
If $\|\Phi\|_\infty\leq C$ and $\bar{\nu}\leq \epsilon$ then there is a smooth $\Pin^-(2)$-monopole $\alpha$ over $Y$, gauge equivalent to $\beta$, such that if $B$ is the connection part of $\pi^*\alpha$ then for every $t\geq 1$ and nonnegative integer $k$ one has 
$$
\sup_{y\in Y}|\nabla^k_B(x-\pi^*\alpha)|_{(t,y)} \leq C_k\sqrt{\nu(0)}e^{-\lambda^+ t}.
$$
\end{Theorem}
\begin{Theorem}[\cite{FC}, Theorem 6.3.2.]\label{thm:expdecay2}
There exists a constant $\lambda_+$ which has the following significance.
For any $C>0$, there exist constants $\epsilon$ and $C_k$ for nonnegative integer $k$ such that the following holds for every $T>1$.
Let $x=(A,\Phi)$ be a $\Pin^-(2)$-monopole in temporal gauge over the band $[-T-2, T+2]\times Y$ such that $x(t)\in U$ for some $t\in[-T-2,T+2]$.
Set 
$$
\bar{\nu}=\|\nabla\vartheta\|_{L^2([-T-2,T+2]\times Y)},\quad 
\nu(t)=\|\nabla\vartheta\|_{L^2(B_t)}.
$$
If $\|\Phi\|_\infty\leq C$ and $\bar{\nu}\leq \epsilon$ then there is a smooth $\Pin^-(2)$-monopole $\alpha$ over $Y$, gauge equivalent to $\beta$, such that if $B$ is the connection part of $\pi^*\alpha$ then for every $t\leq T-1$ and nonnegative integer $k$ one has 
$$
\sup_{y\in Y}|\nabla^k_B(x-\pi^*\alpha)|_{(t,y)} \leq C_k(\nu(-T)+\nu(T))^{1/2}e^{-\lambda^+ (T-|t|)}.
$$
\end{Theorem}
%
%
%
%
\subsection{Energy}
%
%
Let $Z$ be a Riemannian $\Spin^{c_-}$-$4$-manifold possibly noncompact or with boundaries, such as $X$ with a tubular end, or its compact submanifolds $X^t$ or a compact tube $[a,b]\times Y$.
Let $\mu$ be a closed $2$-form in $\Omega^2(i\lambda)$, and assume $\mu$ is the pull-back of $\eta$ on the tube. 
For configurations $(A,\Phi)$, we define the energy by
$$
\E(A,\Phi) = \frac14 \int_Z |F_A-\mu|^2 + \int_Z |\nabla_A\Phi|^2+\frac14\int_Z\left(|\Phi|^2+\frac{s}2\right)^2 - \int_Z\frac{s^2}{16} + 2\int_Z\langle\Phi,\rho(\mu)\Phi\rangle,
$$
where $s$ is the scalar curvature. 
\begin{Proposition}[\cite{KM2}, Chapter \Romnum{2} and Chapter \Romnum{8}]
{\rm (1) }If $(A,\Phi)$ is a $\Pin^-(2)$-monopole on $Z=X^T$ with a finite cylinder $(-1,T]\times Y$ near the boundary $Y$, then 
$$
\E(A,\Phi)=\frac14\int_Z(F_A-\mu)\wedge(F_A-\mu) - \int_Y\langle\Phi|_Y,D_B(\Phi|_Y)\rangle,
$$ 
where $B$ is the boundary connection induced from $A$.\\
{\rm (2)} If $(A,\Phi)$ is a $\Pin^-(2)$-monopole on $[t_0,t_1]\times Y$ in temporal gauge, then 
$$
\frac12\E(A,\Phi) = \vartheta(A(t_1),\Phi(t_1)) - \vartheta(A(t_0),\Phi(t_0)) .
$$
\end{Proposition}
%
%
%
%
\subsection{Compactness}
%
%
We invoke a compactness result due to Kronheimer and Mrowka.
\begin{Proposition}[\cite{KM2}, Theorem 5.1.1]\label{prop:loc-cpt}
Let $Z$ be a compact Riemannian $\Spin^{c_-}$-$4$-manifold  with boundary.
Suppose there exists a constant $C$ so that a sequence $(A_n,\Phi_n)$ of smooth solutions to $\Pin^-(2)$-monopole equations satisfies the bound $\E(A_n,\Phi_n)\leq C$. 
Then there exists a sequence $g_n$ of {\rm (}smooth{\rm )} gauge transformations with the following properties:
after passing to a subsequence, the transformed solutions $g_n(A_n,\Phi_n)$ converges weakly in $L^2_1$ to a $L^2_1$-configuration $(A,\Phi)$ on $Z$, and converges  strongly in $C^\infty$ on every interior domain $Z^\prime\subset Z$.
\end{Proposition}
\begin{Corollary}\label{cor:limit}
Let $x(t)=(A(t),\Phi(t))$ be a smooth monopole on $[-1,\infty)\times Y$ in temporal gauge.
If $\E(A,\Phi)$ is finite, then $[x(t)]$ converges in $\B_Y$ to some critical point as $t\to \infty$.
\end{Corollary}
\proof 
By translation, $(A_T,\Phi_T)=(A,\Phi)|_{[T-1,T+1]\times Y}$ can be considered as a monopole on $[-1,1]\times Y$.
Let $T_n$ be any sequence with $T_n\to\infty$ as $n\to\infty$.
Since $\E(A,\Phi)$ is finite, $\E(A_{T_n},\Phi_{T_n})\to 0$ as $n\to\infty$. 
Then, after some gauge transformations, we may assume $(A_{T_n},\Phi_{T_n})$ converges in $C^\infty$ on $(-1,1)\times Y$ to the pull-back of some critical point.
From this, the corollary is proved.
\endproof
\begin{Proposition}
Let $X$ be a $\Spin^{c_-}$-$4$-manifold $X$ with an end $[-1,\infty)\times Y$. 
If a smooth monopole $(A,\Phi)$ over $X$ has a finite energy $\E(A,\Phi)$, then we have either 
$$
\Phi = 0,\quad \text{ or }\quad \|\Phi\|_{C^0}\leq - \frac12 \inf_{x\in X}s(x) + 4\|\mu \|_{C^0},
$$
where $s$ is the scalar curvature of $X$. 
\end{Proposition}
\proof
By \corref{cor:limit}, we may assume $(A,\Phi)$ converges to a monopole $(B,\Psi)$ on $Y$. 
If $|\Phi|$ takes its maximum on $X$, then the argument in \cite{KM}, Lemma 2, implies the proposition.
Otherwise  we have $\|\Phi\|_{C^0}=\|\Psi\|_{C^0}$. 
Since $(B,\Psi)$ is a $3$-dimensional monopole, $\Psi$ also satisfies
$$
\Psi = 0\quad \text{ or }\quad \|\Psi\|_{C^0}\leq - \frac12 \inf_{y\in Y}s(y) + 4\|\eta \|_{C^0}.
$$
\endproof
%
%
\subsection{Weighted moduli spaces}
%
%
Throughout this subsection, we assume $X$ is a $\Spin^{c_-}$-$4$-manifold with the end $[-1,\infty)\times S^3$.
Let us fix a smooth reference connection $A^0$ which is the pull-back of $\theta$ on the tube $[0,\infty)\times S^3$.
For later purpose, we choose an integer $k$ so that $k\geq 3$.
We consider the space of configurations
$$
\CC^w=\{(A^0 +a,\Phi) \,|\, a\in L^{2,w}_k(\Lambda^1(i\lambda)), \Phi\in L^{2,w}_k(S^+)\}.
$$
Let us consider the set of gauge transformations
$$
\G^w=\{g\in L^{2}_{k+1, {\rm loc}}(\Gamma (\X\times_{\{\pm 1\}}\U(1)))\,|\, \nabla_0 g \in L_k^{2,w}\},
$$
where $\nabla_0$ denotes the covariant derivative of $A^0$.
We can prove,
%
%
\begin{Proposition}[ \cite{Taubes}, Section 7, {\it Cf.} \cite{DF}, \S4.3, \cite{FC} Chapter 2]
{\rm (1)} Let  $L\G^w$ be the set defined by
$$
L\G^w =\{\xi\in L_{k+1,{\rm loc}}^2(\Lambda^0(i\lambda)) \,|\, \nabla_0\xi\in L_{k+1}^{2,w}\}.
$$
Then each element $\xi\in L\G^w$ tends to a limit in $\Lie\Gamma_\theta\cong i\R$ at infinity, and therefore the evaluation map is defined{\rm :}
$$
r\colon L\G^w\to \Lie\Gamma_{\theta}.
$$
When we define the inner product on $L\G^w$  by
$$
\langle\xi,\eta\rangle = \langle\nabla_0\xi,\nabla_0\eta\rangle_{L^{2,w}_{k+1}} + \langle r(\xi),r(\eta)\rangle_{i\R},\quad \xi,\eta\in L\G^w,
$$
$L\G^w$ is a Hilbert space.\\
{\rm (2)} $\G^w$ is a Hilbert Lie group which is modeled on the Lie algebra  $L\G^w$.
Each element $g\in\G^w$ tends to a limit in $\Gamma_\theta$ at infinity, and the evaluation map is defined{\rm :}
$$
R\colon \G^w\to \Gamma_{\theta}.
$$

\end{Proposition}
Let $\G^w_0$ be the kernel of $R$. Then $\G^w/\G^w_0\cong\Gamma_\theta$.
Now the Lie algebra of $\G^w_0$ is given by
$$
L\G^w_0 =L_{k+1}^{2,w}(\Lambda^0(i\lambda)). 
$$

For a configuration $(A,\Phi)\in \CC^w$, the infinitesimal $\G_0^w$-action is given by the map
$$
\I_\Phi\colon L_{k+1}^{2,w}(\Lambda^0(i\lambda))\to L_k^{2,w}(\Lambda^1(i\lambda)\oplus S^+)
$$
defined by $\I_\Phi(f)=(-2df, f\Phi)$.
When $\I_\Phi^*$ is the formal adjoint of $\I_\Phi$, the adjoint of $\I_\Phi$ with respect to the {\it weighted} norm is given by
$$
\I_\Phi^{*,w}(\alpha) = w^{-2}\I_\Phi^*(w^2\alpha).
$$
This gives the decomposition({\it Cf.}~\cite{FC}):
$$
 L_k^{2,w}(\Lambda^1(i\lambda)\oplus S^+) = (\ker\I_\Phi^{*,w}\subset L_k^{2,w})\oplus \I_\Phi(L_{k+1}^{2,w}).
$$
Since the $\G^w_0$-action on $\CC^w$ is free, the quotient space $\tilde{\B}^w = \CC^w/\G^w_0$ is a Hilbert manifold, with a local model 
$$
T_{[(A,\Phi)]}\tilde{\B}^w = \ker\I_\Phi^{*,w}\cap L_k^{2,w}.
$$
The $\Pin^-(2)$-monopole map is defined as
\begin{gather*}
\Theta=\Theta_\mu\colon \CC^w\to L^{2,w}_{k-1}(\Lambda^+(i\lambda)\oplus S^-),\\
\Theta_\mu(A,\Phi)=\left(\frac12 F^+_{A}-q(\Phi)-\mu, D_{A}\Phi\right),
\end{gather*}
where $\mu$ is a (compact-supported) $i\lambda$-valued self-dual $2$-form. 
The moduli space is defined by $\M = \Theta^{-1}(0)/\G^w$.
\begin{Proposition}
The moduli space $\M$ is compact.
\end{Proposition}
\proof
Let $[(A_n,\Phi_n)]$ be any sequence in $\M$.
In general, one can prove that the sequence has a chain convergent subsequence. (\cite{DF}, Chapter 5 and \cite{FC}, Chapter 7.)
Since there is only one critical point on $Y=S^3$, the subsequence converges in   $\M$.
\endproof
The differential of $\Theta$ at $x=(A,\Phi)$ is given by
\begin{gather*}
\D_{(A,\Phi)}=D\Theta\colon L_k^{2,w}(\Lambda^1(i\lambda)\oplus S^+ )\to L^{2,w}_{k-1}(\Lambda^+(i\lambda)\oplus S^-),\\
\D_{(A,\Phi)}(a,\phi)=\left(\frac12 d^+a-Dq_\Phi(\phi), D_{A}\phi+\frac12\rho(b)\Phi\right),
\end{gather*}
where $Dq_\Phi$ is the differential of $q$.
Then 
\begin{equation}\label{eq:DI}
\D_{(A,\Phi)}\circ\I_\Phi(f) = (0, fD_A\Phi).
\end{equation}
Therefore, if $(A,\Phi)$ is a $\Pin^-(2)$-monopole solution, then $\D_{(A,\Phi)}\circ\I_\Phi(f)=0$, which forms the deformation complex:
$$
\begin{CD}
0@>>>L_{k+1}^{2,w}(\Lambda^0(i\lambda)) @>{\I_\Phi}>>L_k^{2,w}(\Lambda^1(i\lambda)\oplus S^+)@>{\D_{(A,\Phi)}}>> L^{2,w}_{k-1}(\Lambda^+(i\lambda)\oplus S^-)@>>>0.
\end{CD}
$$
The cohomology groups are denoted by $H^i_{(A,\Phi)}$.

The monopole map $\Theta$ defines a $\Gamma_\theta$-invariant section of a bundle over $\tilde{\B}^w$ whose linearization is given by $\I_\Phi^{*,w}\oplus \D_{(A,\Phi)}$.
When $Y$ is the standard $S^3$, the virtual dimension of the moduli space "framed at infinity" 
$\tilde{\M}=\Theta^{-1}(0)/\G^w_0\subset \tilde{\B}^w$ is given by
$$
\ind^+ (\I_\Phi^{*}\oplus \D_{(A,\Phi)}) +\dim\Gamma_\theta= d(c)+1,
$$
where $d(c)$ is in \eqref{eq:dc}.
The genuine moduli space is $\M=\tilde{\M}/\Gamma_\theta$ whose virtual dimension is $d(c)$.
In general, $\M$ and $\tilde{\M}$ are not smooth manifolds, and we need to perturb the equations.
Before that, we introduce a term.
\begin{Definition}
The moduli space $\M$ is said to be {\it regular} if all of elements $[(A,\Phi)]$ of  $\M$ have $H^2_{(A,\Phi)}=0$.
\end{Definition}
\begin{Remark}
If $\M$ contains no reducibles, then $H^0_{(A,\Phi)}=0$ for all $[(A,\Phi)]\in\M$. 
But the converse is not necessarily true, because the stabilizer of a $\Pin^-(2)$-monopole reducible $[(A,0)]$ on a twisted $\Spincm$-structure is $\deux$, and then $H^0_{(A,0)}=0$. 
\end{Remark}
If $b_+(X;l)\geq 1$, by perturbing the equation by adding a {\it compactly-supported} self-dual $2$-form as in \eqref{eq:monopole}, we obtain a smooth $\tilde{\M}$: 
\begin{Theorem}[\cite{FC}, Proposition 8.2.1]
Suppose $b_+(X;l)\geq 1$. 
For generic  {\it compactly-supported} self-dual $2$-forms, by perturbing the equations as in \eqref{eq:monopole}, the perturbed moduli space $\tilde{\M}$ is regular and contains no reducibles, and therefore is a smooth manifold of dimension $d(c)+1$.
Then  $\M$ is a smooth manifold of dimension $d(c)$.
\end{Theorem}
When $\M(X)$ has no reducibles, the cutting-down method described in \subsecref{subsec:cutdown} works well for $\M(X)$ in this section.
However, if  $\M(X)$ contains a reducible, we need a little care for it as follows.
Choose loops  $\gamma_1,\ldots,\gamma_b$, where $b=b_1(X;l)$, satisfying the conditions (C1) and (C2) in \subsecref{subsec:mu-map1}.
Define the map $\mathfrak{h}\colon \M(X)\to T^b$ by $\mathfrak{h}=\mathrm{Hol}_{\gamma_1}\times\cdots \times \mathrm{Hol}_{\gamma_b}$.

\begin{Theorem}\label{thm:cut-down}
Suppose $b_+(X;l)=0$, $\tilde{c}_1(E)^2=\sign(X)$ {\rm (}and hence the Dirac index is $0$ and $d(c)=b_1(X;l)${\rm )}.
For a generic choice of $\alpha\in T^b$ and a compactly-supported self-dual $2$-form, the cut-down moduli space $\M\cap \mathfrak{h}^{-1}(\alpha)$ is regular, and therefore consists of one reducible point and a finite  number of irreducible points.
\end{Theorem}
\proof 
The proof is similar to that in  \cite{Pin2}, Subsection 4.8. 
Due to the noncompactness of $X$, we need to modify the following point:
The space $L^{2,w}_{k}(\Lambda^1(i\lambda))$ is decomposed into the direct sum of $\ker d^+$ and its complement $(\ker d^+)^\perp$.
Furthermore, since $b_+(X;l)=0$, $d^+\colon (\ker d^+)^\perp \to L^{2,w}_{k-1}(\Lambda^+(i\lambda))$ is an isomorphism.
Mimicking the argument in the proof of Lemma 14.2.1 of \cite{FC}, we can prove the following.
\begin{Claim}
Fix a compact codimension-0 submanifold $K\subset X$ and let $\Omega^+_{X,K}(i\lambda)$ be the space of smooth self-dual $2$-forms on $X$ supported on $K$ with $C^\infty$-topology. 
For $(b,\mu)\in \ker d^+\oplus \Omega^+_{X,K}(i\lambda)$, let $A(b,\mu)$ be the connection $A^0 + b + (d^+)^{-1}(\mu)$. 
Let ${\cal R}$ be the set of  $(b,\mu)\in \ker d^+\oplus \Omega^+_{X,K}(i\lambda)$ such that $D_{A(b,\mu)}$ is surjective, (and hence, of course, also injective).
Then ${\cal R}$ is open-dense.  
\end{Claim}
\begin{Claim}
There exists a gauge invariant open-dense subset ${\cal R}^\prime\subset \CC^w\times \Omega^+_{X,K}(i\lambda)$ such that the restriction of the $\Pin^-(2)$-monopole map $\Theta_\mu$ to ${\cal R}^\prime$ has $0$ as regular value.
\end{Claim}
With these understood, 
$$
{\cal Z}=\{(A,\Phi,\mu)\in {\cal R}^\prime\,|\,\Theta_\mu(A,\Phi)=0\,\}
$$
is a submanifold in ${\cal R}^\prime$.
Then it suffices to apply the Sard-Smale theorem to the map
$$
\mathfrak{h}\times\pi\colon \mathcal{Z}\to T^b\times \Omega^+_{X,K}(i\lambda),
$$
where $\pi$ is the projection.
\endproof 

%
%
\section{Proofs of gluing formulae}\label{sec:gluing}
%
%
The purpose of this section is to give proofs of the gluing formulae in \subsecref{subsec:GF}.
%
%
\subsection{Gluing monopoles}\label{subsec:gluing0} 
%
%
Let $X_1$ and $X_2$ be $\Spin^{c_-}$-$4$-manifolds with ends $[-1,\infty)\times Y_1$ and $[-1,\infty)\times Y_2$, where $Y_1=\bar{Y}_2=S^3$.
Fix a reducible solution $(\theta,0)$ on $S^3$, and choose a $C^\infty$ reference connection $A_i^0$ on each $X_i$ which is the pull-back of $\theta$ on the tube.
Let $x_i=(A_i,\Phi_i)$ be finite energy monopole solutions on $X_i$ ($i=1,2$).
Furthermore, we also suppose $H^2_{x_1}=H^2_{x_2}=0$. 
We assume each $A_i$ is in temporal gauge on the tube, and if necessary, consider it as a one-parameter family of connections $\theta+a_i(t)$ on the tube.
The spinors $\Phi_i$ are also considered as one-parameter families $\Phi_i(t)$ on the tube.

Now, we construct an approximated solution on $X^{\#T}$ from $(A_1,\Phi_1)$ and $(A_2,\Phi_2)$ by splicing construction.
Choose a smooth cut-off function $\gamma$, with $\gamma(t)=1$ for $t\leq 0$ and $\gamma(t)=0$ for $t\geq 1$.
Define $x_1^\prime=(A_1^\prime,\Phi_1^\prime)$ over $X_1$ by
\begin{equation}\label{eq:prime}
\begin{aligned}
A_1^\prime =& \theta + \gamma(t-T+3)a_1(t),\\
\Phi_1^\prime=& \gamma(t-T+3)\Phi_1(t).
\end{aligned}
\end{equation}
Define $x_2^\prime=(A_2^\prime,\Phi_2^\prime)$ over $X_2$ in a similar fashion.

Fix an identification of the $\Spin^{c_-}$-structures on $[0,2T]\times Y_1$ and $[0,2T]\times Y_2$ with respect to $\theta$. 
Note that the $\Spin^{c_-}$-structures on the tubes are untwisted, which are identified with  ordinary $\Spin^{c}$-structures. 
The all possibilities of such identifications are parameterized by $\Gamma_\theta$, which are called the {\it gluing parameters}.
If we fix an identification $\sigma_0$, then the other identifications are indicated as $\sigma=\exp(v)\sigma_0$ for $v\in\Lie\Gamma_\theta\cong i\R$.
For an identification $\sigma$, we can glue $x_1^\prime$ and $x_2^\prime$ via $\sigma$ to give a configuration over $X^{\#T}$.
The glued configuration is denoted by
$$
x^\prime(\sigma) =(A^\prime(\sigma), \Phi^\prime(\sigma)).
$$
Then it is easy to see the following
\begin{Proposition}
For each $i=1,2$, let $\Gamma_i$ be the stabilizer of the monopole $x_i$.
Then $x^\prime(\sigma_1)$ and $x^\prime(\sigma_2)$ are gauge equivalent if and only if $[\sigma_1]=[\sigma_2]$ in $\Gamma_\theta/(\Gamma_1\times\Gamma_2)$, where $\Gamma_i$ are the stabilizers of $x_i$.
\end{Proposition}
Let $\mathrm{Gl}=\Gamma_\theta/(\Gamma_1\times\Gamma_2)$. 
Define the map $\mathfrak{F}^\prime\colon\mathrm{Gl}\to \B(X^{\#T})$ by the splicing construction above: $[\sigma]\mapsto [x^\prime(\sigma)]$.
If $H^2_{x_1}=H^2_{x_2}=0$ and $T$ is sufficiently large, then we can find in a unique way a monopole solution $x(\sigma)$ on $X^{\#T}$ near the spliced configuration $x^\prime(\sigma)$. 
(This is standard in gluing theory. See \cite{DK,DF,FC,Nicolaescu}.)
(The construction is explained in the Appendix.) 
Then we have a smooth map
\begin{equation}\label{eq:I}
\mathfrak{F}\colon\mathrm{Gl}\to \M(X^{\#T}),\quad [\sigma]\mapsto[x(\sigma)].
\end{equation}

Before proceeding, we give another description of the spliced family $\{[x^\prime(\sigma)]\}$ for gluing parameters $\sigma\in\Gamma_\theta$.
According to the definition of $x^\prime(\sigma)$, for different $\sigma$, $x^\prime(\sigma)$ are objects on different bundles parameterized by $\sigma$.
It is convenient if we can represent all $[x^\prime(\sigma)]$ as objects on a fixed identification, say $\sigma_0$, of  bundles.  
This is also done in \cite{DK}, \S7.2.4, in the ASD case.

Recall $X^{\#T}=X_1^0\cup([-T,T]\times Y)\cup X_2^0$, and $X_1^{2T}$ and $X_2^{2T}$ are assumed to be embedded in $X^{\#T}$.
Choose a smooth function $\lambda_1$ on $X^{\#T}$ such that 
$\lambda_1=1$ on $X_1^0$, $\lambda_1=0$ on $X_2^0$ and 
$$
\lambda_1(t,y) = \left\{
\begin{aligned}
1 & \quad-T\leq t\leq -1, y\in Y,\\
0 & \qquad 1\leq t\leq T, y\in Y,
\end{aligned}\right.
$$
and satisfies $|\nabla\lambda|=O(1)$.
Define another function $\lambda_2$ on $X^{\#T}$ by $\lambda_2 =1-\lambda_1$.
Let $v\in\Lie \Gamma_\theta = i\R$, and $\sigma=\sigma_0\exp(v)$.
Define gauge transformations $h_1$ and $h_2$ on $X^{\#T}$ by 
\begin{equation}\label{eq:hi}
\begin{aligned}
h_1 =& \exp(\lambda_2 v)\\
h_2 =& \exp(-\lambda_1 v)
\end{aligned}
\end{equation}
Note that $h_1h_2^{-1}=\exp(\lambda_1+\lambda_2)v=\exp v$.
Then $h_1x_1^\prime = h_2x_2^\prime$ over $[-2,2]\times Y$ on which $x_1^\prime$ and $x_2^\prime$ are flat, and therefore we can glue them. 
The glued configuration is denoted by $x^\prime(\sigma_0,v)$.
Then, by definition, it can be seen that $x^\prime(\sigma)$ and  $x^\prime(\sigma_0,v)$ are gauge equivalent.
Often, we will not distinguish these two, and use the same symbol  $x^\prime(\sigma)$.
%
%
\subsection{Gluing maps between the moduli spaces}\label{subsec:gluingformulae}
%
%
The gluing construction \eqref{eq:I} can be globalized to whole moduli spaces.
In fact, we can define the map 
$$
\Xi\colon \tilde{\M}(X_1)\times_{\Gamma_\theta}\tilde{\M}(X_2)\to \M(X^{\#T}),
$$
for sufficiently large $T$.
\begin{Theorem}\label{thm:gluing}
Let $X_1$ and $X_2$ be $\Spin^{c_-}$-$4$-manifolds with ends $[-1,\infty)\times Y_1$ and $[-1,\infty)\times Y_2$, where $Y_1=\bar{Y}_2=S^3$.
Suppose the following.
\begin{itemize}
\item The $\Spincm$-structure on $X_1$ may be twisted or untwisted, and $\M(X_1)$ contains no reducibles.
\item The $\Spincm$-structure on $X_2$ is twisted, and $\M(X_2)$ may contain a reducibles.
\item Both of $\M(X_1)$ and $\M(X_2)$ are regular, and $\dim\M(X_1)=\dim\M(X_2)=0$.
\end{itemize}
Then $\Xi$ is a diffeomorphism between $1$-dimensional compact manifolds.
\end{Theorem}
%
%
%
%
\begin{Theorem}\label{thm:open}
Suppose $X_1$ is a $\Spincm$-$4$-manifold with the end $[-1,\infty)\times S^3$ whose moduli space $\M(X_1)$ is regular and contains no reducibles.
Suppose $X_2$ is a $\Spinc${\rm (}untwisted $\Spincm${\rm )}-$4$-manifold with the end $[-1,\infty)\times S^3$ such that $b_1(X_2)=b_+(X_2)=0$ and $\dim\M(X_2)=-1$.
Then $\Xi$ induces a diffeomorphism 
$$
\M(X_1)\to\M(X^{\#T}).
$$
\end{Theorem}
With the results in the previous subsections understood, we can prove these theorems  by  a similar way to those of the corresponding theorems in the Seiberg-Witten and Donaldson theory. (See \cite{DK,DF,FC,Nicolaescu}).
A proof based on \cite{DK,DF} will be explained in the Appendix.
%
%
\subsection{The images of the map $\mathfrak{F}$}\label{subsec:I}
%
%
To prove the gluing formulae, we want to know what is the homology class of the image of $\mathfrak{F}$ in $H_*(\B)$. 
The homology class depends on whether each of the $\Spincm$-structures on $X_1$ and $X_2$ is twisted or untwisted, and whether each of monopoles $x_1$ and $x_2$ is irreducible or not.
We call an irreducible/reducible monopole on a twisted $\Spincm$-structure {\it $\Pin^-(2)$-irreducible/reducible}, and an irreducible/reducible monopole on an  untwisted $\Spincm$-structure {\it $\U(1)$-irreducible/reducible}.
We assume that the $\Spincm$-structures of $x_2$ is twisted. 
Then $\B(X^{\#T})$ is homotopy equivalent to $\RP^\infty\times T^{b_1(X^{\#T};l)}$. 
Let $\nu_2^*=\mu_{\mathcal{F}}(\nu_2)$ and $\hat{\delta}^*=\hat{\mu}_{\mathcal{E}}(\delta)$ for the loops $\nu_2$ and  $\delta$ in $X^{\#T}$ chosen as in {\rm \subsecref{subsec:GF}}. 
For monopoles $x_1$ and $x_2$ on $X_1$ and $X_2$, let $C$ be the image of $\mathfrak{F}$.
Suppose $x_1$ and $x_2$ are not $\U(1)$-reducible. 
Then $C$ is a circle.
\begin{Theorem}\label{thm:C}
 For the homology classes  $[C]\in H_1(\B;\Z)$ and $[C]_2\in H_1(\B;\Z_2)$ of $C$, we have the following{\rm :} 
\begin{enumerate}
\item If $x_1$ is a $\U(1)$-irreducible and $x_2$ is a $\Pin^-(2)$-reducible, then $\langle \nu_2^*, [C]_2\rangle\neq 0$.
\item If $x_1$ is a $\U(1)$-irreducible and $x_2$ is a $\Pin^-(2)$-irreducible, then $[C]=[C]_2= 0$.
\item  If $x_1$ is a $\Pin^-(2)$-irreducible and $x_2$ is a $\Pin^-(2)$-reducible, $\langle \hat{\delta}^*, [C]\rangle= \pm1$.
\item If both of $x_1$ and $x_2$ are $\Pin^-(2)$-irreducibles, then $\langle \hat{\delta}^*, [C]\rangle=\pm 2$.
\end{enumerate}
\end{Theorem}
Before proving the theorem, we give some preliminaries.
In the following, we simplify the notation as $\G=\G^w$, $\G_0=\G_0^w$ and $\K=\K_\gamma$ which is in \remref{rem:Kgamma}. 
Let $\K_0=\K\cap\G_0$.
For each $i=1, 2$, let $\Sc_i$ be the set of solutions which are $\G$-equivalent to $x_i$. 
Now, we prove the assertions (1) and (2).
\proof[Proof of {\rm (1)} and {\rm (2)}]
We have a commutative diagram whose vertical and horizontal arrows are exact:
$$
\begin{CD}
@. 1 @.  1 \\
@. @AAA @AAA\\
@. \deux @= \deux \\
@. @AAA @AAA\\
1@>>> \G_0 @>>> \G @>>> \Gamma_\theta@>>> 1,\\
@. @AAA @AAA @|\\
1@>>> \K_0 @>>> \K @>>> \Gamma_\theta@>>> 1,\\
@. @AAA @AAA \\
@. 1 @.  1 
\end{CD}
$$
We also have the following diagrams of various quotient maps:
\[\xymatrix{
& \Sc_2/\K_0 \ar[ld]_{\deux} \ar[rdd]^{\Gamma_\theta} &\\
\Sc_2/\G_0 \ar[rdd]_{\Gamma_\theta} &&& \Sc_1/\G_0 \ar[dd]_{\Gamma_\theta}\\
&& \Sc_2/\K\ar[ld]^{\deux} & \\
&\Sc_2/\G& & \Sc_1/\G
}\]
By definition, $\Sc_1/\G$ and $\Sc_2/\G$ are  one-point sets.
Then $\Sc_1/\G_0$ is a circle on which $\Gamma_\theta$ acts freely.
Hence, $C=\Ima \mathfrak{F}$ can be written as 
$$
C=(\Sc_1/\G_0)\times_{\Gamma_\theta}(\Sc_2/\G_0) = \frac{(\Sc_1/\G_0)\times_{\Gamma_\theta}(\Sc_2/\K_0)}{\deux} = (\Sc_2/\K_0)/\deux.
$$
First, let us consider the case of (2).
In this case, $\G$ acts on $\Sc_2$ freely. 
Therefore $\Sc_2/\K_0\cong \Gamma_\theta\times \deux$, and we can see that the homology class of $C$ is zero.
In the case of (1), each element of $\Sc_2$ has the stabilizer $\deux\subset\G$.
Since $\G_0\cap\deux=\{1\}$, we see that $\Sc_2/\G_0\cong \Gamma_\theta/\deux$ and $[C]$ is the generator of $H_1(\RP^\infty;\Z_2)\cong\Z_2$.  
\endproof
We give an alternative proof which gives more intuitive understanding of the homology class of $[C]$.
\proof[Alternative proof of {\rm (1)} and {\rm (2)}]
For a section $\Phi_i$ of the spinor bundle $S^+_i$ of $c_i$ ($i=1,2$), let $\Phi_i^\prime$ be the cut-off section as in \eqref{eq:prime}.  
Define
$$
\Gamma(S^+_i)^\prime := \{\Phi_i^\prime\,|\, \Phi_i\in \Gamma(S^+_i)\,\}.
$$
Let $S^+_{\sigma_0}=S^+_1\#_{\sigma_0}S^+_2$ be the glued spinor bundle over $X^{\#T}$ via the gluing parameter $\sigma_0$.
Then we can assume $\Gamma(S^+_1)^\prime\oplus\Gamma(S^+_2)^\prime$ is a subspace of $\Gamma(S^+_{\sigma_0})$ via the splicing construction.
For the monopoles $x_i=(A_i,\Phi_i)$ and $\sigma=\sigma_0\exp(v)$ ($v\in \Lie\Gamma_\theta$), define the configuration $y(\sigma)$ on $X^{\#T}$ by
$$
y(\sigma) = (A^\prime(\sigma_0),(\exp(v)\Phi_1^\prime,\Phi_2^\prime))
$$
where $(\exp(v)\Phi_1^\prime,\Phi_2^\prime)\in \Gamma(S^+_1)^\prime\oplus\Gamma(S^+_2)^\prime$ is assumed to be an element of  $\Gamma(S^+_{\sigma_0})$ as above.
Then the homology class $[C]$ is represented by
$$
\{y(\sigma)\}_{\sigma\in\Gamma_\theta} = \{A^\prime(\sigma_0)\}\times C_1\times \{\Phi_2^\prime\},
$$
where $C_1=\{\exp(v)\Phi_1^\prime\}_v$.
Note that $C_1$ is a circle in the complex line generated by $\Phi_1^\prime$.
Let $\mathcal{P}=(\Gamma(S^+_{\sigma_0})\setminus\{0\})/\deux$. 
Then $\mathcal{P}$ is homotopy equivalent to $\RP^\infty$.
Consider the following map
$$
q\colon  \mathcal{P}\to \B^*=\CC^*/\G,\quad [\Phi]\mapsto [(A^\prime(\sigma_0),\Phi)].
$$
Then the map $q$ induces an injective homomorphism 
$$
q_*\colon H_*(\mathcal{P})\to H_*(\B^*).
$$
If $x_2$ is a $\Pin^-(2)$-reducible, then $\Phi_2^\prime\equiv 0$ and $[C_1\times\{\Phi_2^\prime\}]$ is a generator of $ H_1(\mathcal{P} ;\Z_2)$.
On the other hand, if $x_2$ is a $\Pin^-(2)$-irreducible, then $\Phi_2^\prime\neq 0$ and $[C_1\times\{\Phi_2^\prime\}]$ is null-homologous.
\endproof
%
%
%
%
In order to prove the assertions (3) and (4), we first consider the gluing of connections.
For each $i=1,2$, let $A_i$ be a connection on the characteristic bundle $E_i$ for $c_i$.
For $\sigma\in\Gamma_\theta$, let  $A_1\#_\sigma A_2$ be the spliced connection on $E=E_1\#_\sigma E_2$ as in \subsecref{subsec:gluing0}.  
Note that $A_1\#_\sigma A_2$ is gauge equivalent to $A_1\#_{-\sigma} A_2$, where $-\sigma = \sigma\exp\pi i$.
\begin{Lemma}\label{lem:S}
Let $S=\{[A_1\#_\sigma A_2]\}_{\sigma\in\Gamma_\theta/\deux}\subset\A(E)/\G$ be the set of gauge equivalence classes of the family $\{A_1\#_\sigma A_2\}_{\sigma\in\Gamma_\theta}$.
Then its homology class $[S]\in H_1(\A(E)/\G;\Z)$ satisfies the following:
\begin{enumerate}
\item $\langle\bar{\alpha}_i,[S]\rangle = \langle\bar{\beta}_j,[S]\rangle =0$ for  $i=1,\ldots,b_1(l_1)$, $j=1,\ldots,b_1(l_2))$,
\item $\langle\bar{\delta},[S]\rangle =\pm 1$,
\end{enumerate}
where $\bar{\alpha}_i$, $\bar{\beta}_j$, $\delta\in H^1(\A/\G;\Z)$ as in \remref{rem:jacobian}.
\end{Lemma}
\proof
The assertion (1) is obvious.
We prove the assertion (2).
Fix $\sigma_0\in\Gamma_\theta$ as based point, and the spliced connections $A_1\#_\sigma A_2$ for other $\sigma$ are constructed by using \eqref{eq:hi} as in \subsecref{subsec:gluing0}.
For $\sigma\in\Gamma_\theta$, $A_1\#_{\sigma} A_2$ and $A_1\#_{-\sigma} A_2$ are gauge equivalent by the gauge transformation $\check{g}$ such that
$$
\check{g}=\left\{\begin{aligned}
1 \qquad& \text{ on $X_1^0$}\\
-1 \qquad& \text{ on $X_2^0$}\\
\exp(\lambda_2\pi i) &\text{  on $[-T,T]\times Y$}
\end{aligned}\right.
$$
where $\lambda_2$ is the function defined around \eqref{eq:hi}.
On the other hand, for any $w$ with $0<w<\pi$, if we put $\sigma_w=\sigma\exp(iw)$, then $A_1\#_\sigma A_2$ and $A_1\#_{\sigma_w}A_2$ are not gauge equivalent.
Therefore $S$ is a circle embedded in $\A(E)/\G$.
By taking homotopy class and projection, we have a surjection $\rho\colon \G\to H^1(X;l)/\mathrm{Tor}$ (see \cite{Pin2}, Lemma 4.22).
Then it suffices to prove $\langle\rho(\check{g}),[\delta]\rangle=\pm 1$.
To see this,  consider the following commutative diagram:
$$
\begin{CD}
\tilde{\G}=\Map(\X;\U(1))@>{\tilde{\rho}}>> H^1(\X;\Z)\\
@A{\varpi}AA @A{\varpi^\prime}AA\\
\G=\Gamma(\X\times_{\deux}\U(1)) @>{\rho}>> H^1(X;l)/\mathrm{Tor},
\end{CD}
$$
where the maps $\varpi$ and $\varpi^\prime$ are the pull-back maps to the double covering $\X$.
Note the following:
\begin{itemize}
\item The image of $\varpi$ is the fixed point set $\tilde{\G}^I$, where the $I$-action is given by $I\tilde{g}=\overline{\iota^*\tilde{g}}$. 
\item Let $\X_i$ ($i=1,2$) be the double coverings of $X_i$. 
Then $\X$ is the connected sum ``at two points" of $\X_1$ and $\X_2$. 
That is, this is obtained as follows: 
For each $i=1,2$, removing two $4$-balls from each of $\X_i$, we obtain a manifold $\X_i^\prime$ whose boundary $\tilde{Y}_i$ is a disjoint union of two $S^3$.
Then $\X = \X_1^\prime\cup_{\tilde{Y}_1 = \tilde{Y}_2}\X_2^\prime$.
\end{itemize}
Consider a circle $\tilde{\delta}$ embedded in $\X$ starting from a point $x_1$ in $\X_1^\prime$ and entering $\X_2^\prime$ via a component of $\tilde{Y}_1 = \tilde{Y}_2$ and returning to $x_1$ via another component of $\tilde{Y}_1 = \tilde{Y}_2$.
Then the restriction of $\varpi(\check{g})$ to $\tilde{\delta}$ gives a degree one map from $\tilde{\delta}$ to $\U(1)$.
\endproof
\proof[Proof of {\rm (3)} and {\rm (4)}]
Let us consider the projection 
$$
\pi\colon C=(\Sc_1/\G_0)\times_{\Gamma_\theta}(\Sc_2/\G_0)\to S,
$$
which is defined by $\pi([x_1],[x_2]) = [A_1\#_\sigma A_2]$, where each $A_i$ is the connection part of $x_i$.
Note that $\pi$ is a map between two $S^1$.
Then, $\pi$ has degree $1$ in the case of (3), and degree $2$ in the case (4).
\endproof
%
%
\subsection{Proofs of the gluing theorems}\label{subsec:proofGF}
%
%
\proof[Proof of \thmref{thm:GF}]
First assume $d(c_i)=\dim\M(X_i)=0$ for $i=1,2$.
For each $i=1,2$, let $X_i^\prime$ be the manifold with cylindrical end obtained from removing a $4$-ball from $X_i$. 
By perturbing the equations with a compactly-supported $2$-form, \thmref{thm:open} implies that $\tilde{\M}(X_i)\cong \tilde{\M}(X_i^\prime)$ for a metric on $X_i$ with long neck.
By the assumption, $\M(X_1^\prime)$ consists of odd numbers, say $k$, of $\U(1)$-irreducible points.
The assumption that $\dim \M(X_2)=0$ implies $b_1^l= b_1(X_2;l_2)=0$, and then $\M(X_2)$ consists of one $\Pin^-(2)$-reducible point and maybe several $\Pin^-(2)$-irreducible points.
By \thmref{thm:gluing}, $\M(X^{\#T})$ is a disjoint union of several circles:
$$
\M(X^{\#T})=\bigcup_{i=1}^k C_i\cup \bigcup_j C_j^\prime,
$$
where $C_i$ are obtained by gluing $\U(1)$-irreducibles and a $\Pin^-(2)$-reducible, and $C_j^\prime$  are made from  $\U(1)$-irreducibles and  $\Pin^-(2)$-irreducibles.
Then \thmref{thm:C}(1)(2) implies that $\langle h,[\M(X^{\#T})]\rangle=k$ mod $2$, and this implies the theorem. 

In the case when $d(c_1)$ or $d(c_2)$ is positive,  \thmref{thm:gluing} can be generalized to give the diffeomorphism between $1$-dimensional cut-down moduli spaces:
$$
\Xi\colon \tilde{M}_1\times_{\Gamma_\theta} \tilde{M}_2 \to M_T,
$$
where
\begin{equation}\label{eq:MT} 
\begin{aligned}
\tilde{M}_1 &= \tilde{\M}(X_1)\cap \bigcap _{i\in I}V_{\alpha_{i}}\cap \bigcap_{k=1}^a V_{x_k},\\
\tilde{M}_2 &= \tilde{\M}(X_2)\cap\mathfrak{h}^{-1}(\alpha),\\
M_T &= \M(X^{\#T})\cap\bigcap _{i\in I}V_{\alpha_{i}}\cap \bigcap_{k=1}^a V_{x_k}\cap\mathfrak{h}^{-1}(\alpha),
\end{aligned}
\end{equation}
and $\tilde{M}_1$, $\tilde{M}_2$ and $M_T$ are assumed to be smooth and  $1$-dimensional.
When $N$ is a closed submanifold of $\M(X^{\#T})$, as a homology class, 
$$
[N\cap V_{x_0}] = \mu_{\cal F}(x_0)\cap [N] = (\mu_{\cal F}(\nu)\cup\mu_{\cal F}(\nu))\cap [N].
$$
From these, the theorem follows.
\endproof
\proof[Proof of \thmref{thm:blowup}]
This is a corollary of  \thmref{thm:open}.
\endproof
\proof[Proof of \thmref{thm:twistedGF}]
The proof is similar to that of \thmref{thm:GF}, by using  \thmref{thm:C} (3)(4).
\endproof

\proof[Proof of \thmref{thm:vanishing}]
For each $i=1,2$, $\M(X_i)$ is perturbed to have no reducibles since $b_+(X_i;l_i)\geq 1$.
The cut-down moduli space $M_T$ as in \eqref{eq:MT}  is a disjoint union of circles $C_i$.
In the case (\romnum1), each $C_i$ is null-homologous by \thmref{thm:C}(2).
In the case (\romnum2),  $\langle \gamma_0^*,[C_i]\rangle=\pm 2$ by \thmref{thm:C}(4).
Therefore the $\Z_2$-valued invariant is zero.  
\endproof
By the proof of the case (\romnum2) of \thmref{thm:vanishing}, \thmref{thm:ZGF} is true if the glued moduli space is orientable.
The orientability of the glued moduli space follows from the next lemma:
\begin{Lemma}
For $i=1,2$, let $X_i$ be a closed oriented connected $4$-manifold with a twisted $\Spincm$-structure $c_i$ whose Dirac index is positive and even, and $A_i$ be a connection on the characteristic bundle $E_i$.
Then for $S$ in \lemref{lem:S}, the restriction of $\ind\delta_{Dirac}$ to $S$, $\ind (\delta_{Dirac}|_S)$, is orientable.  
\end{Lemma}
\proof
We construct a framing of the index bundle $\ind (\delta_{Dirac}|_S)$.
For simplicity, we assume $\ind D_{A_2}=2$, and the general case will be clear.
Let us consider the family $\{D_{A_1\#_\sigma A_2}\}_{\sigma\in\Gamma_\theta}$.
By Proposition 2.2 in \cite{AS4}, we may assume $\Coker D_{A_1\#_\sigma A_2}=0$ for any $\sigma$.
Since $\ker D_\theta =0$ on $S^3$, we can construct an isomorphism for each $\sigma$ (\cite{DF}, \S3.3):
$$
\alpha_\sigma\colon\ker D_{A_1}\oplus\ker D_{A_2}\to \ker D_{A_1\#_\sigma A_2}.
$$
In the proof of \lemref{lem:S}, we have seen that $A_1\#_\sigma A_2$ is gauge equivalent to $A_1\#_{-\sigma}A_2$ by a gauge transformation $g$.
Now we can see that, for $\psi\in\Ker D_{A_1}$ and $\phi\in\Ker D_{A_2}$, 
$$
\alpha_\sigma(\psi,\phi) = g\alpha_{-\sigma}(\psi, -\phi).
$$
Let $\{\psi^{j}\}$ be a basis for $\ker D_{A_1}$, and $\{\phi^{1},\phi^{2}\}$ be a basis for $\ker D_{A_2}$.
Fix $\sigma_0\in\Gamma_\theta$ and let $\sigma_w=\sigma_0\exp(i w)$ for $0\leq w \leq \pi$, and 
$$
\begin{pmatrix}
\phi_w^{1}\\
\phi_{w}^{2}
\end{pmatrix} =
\begin{pmatrix}
\cos w & -\sin w\\
\sin w &\cos w
\end{pmatrix}
\begin{pmatrix}
\phi^{1}\\
\phi^{2}
\end{pmatrix}.
$$
Then the following gives a framing for $\ind(\delta_{Dirac}|_S)$:
$$
\{\alpha_{\sigma_w}(\psi^{j}, \phi_w^1), \alpha_{\sigma_w}(\psi^{j}, \phi_w^2) \}.
$$
\endproof
\begin{Corollary}\label{cor:ori}
For each $i=1,2$, let $X_i$ be a closed oriented connected $4$-manifold with a twisted $\Spincm$-structure which has the following properties{\rm :}
\begin{itemize}
\item the index of the Dirac operator is positive and even, and 
\item the moduli space $\M(X_i)$ is orientable.
\end{itemize}
Then the glued moduli space $\M(X_1\#X_2)$ is also orientable.
\end{Corollary}
\proof[Proof of \thmref{thm:ZGF}]
Since each of $\M(X_i)$ is orientable, \corref{cor:ori} implies the moduli space of $X_0\#\cdots\# X_n$ is also orientable.
The statement for the invariant is proved by \thmref{thm:C}.
\endproof
%
%
\section{Proofs of \thmref{thm:genus} and \corref{cor:genus-cover}}
%
%
We begin with the proof of \thmref{thm:genus}.
Our proof of \thmref{thm:genus} is similar to the proof of Thom conjecture due to Kronheimer and Mrowka \cite{KM}.
({\it Cf.} \cite{Nicolaescu}.)
%
%
\subsection{Reduction to the case when $\alpha\cdot\alpha=0$}
%
%
Suppose $n:=\alpha\cdot\alpha>0$.
Let $X^\prime = X\# n\overline{\CP}^2$, and $E_i$ $(i=1,\ldots,n)$ be the $(-1)$-sphere in the $i$-th $\overline{\CP}^2$.
Take the connected sum in $X^\prime$, 
$$
\Sigma^\prime=\Sigma\# E_1\# \cdots\#E_n.
$$
Then $[\Sigma^\prime]\cdot[\Sigma^\prime]=0$.

Even if we replace $X$ by $X^\prime$, the $\Pin^-(2)$-monopole invariant is unchanged by \thmref{thm:blowup}.
Furthermore, even if we replace $\X$ by $\X^\prime$, the Seiberg-Witten invariant is also unchanged by the ordinary blow-up formulae \cite{FS, Nicolaescu}.
The quantity $-\chi(\Sigma)$ and $\alpha\cdot\alpha + |\tilde{c}_1(E)\cdot\alpha|$ are also unchanged.
Thus, we may assume $\alpha\cdot\alpha=0$. 

In the remainder of this section, we suppose $(X,\alpha,\Sigma)$ satisfies the assumption of the beginning of \subsecref{subsec:genus}, and 
\begin{itemize}
\item $\alpha=[\Sigma] \in H_2(X;l)$ has infinite order, and
\item  $\alpha \cdot \alpha=0$.
\end{itemize}
%
%
\subsection{The case when $\chi(\Sigma)>0$ }\label{subsec:genus0}
%
%
Here, we prove that, under the assumption of \thmref{thm:genus},  the Euler characteristic of $\Sigma$ cannot be positive: 
\begin{Proposition}\label{prop:genus0}
If $\chi(\Sigma)>0$, then the $\Pin^-(2)$-monopole invariants of $(X,c)$ and the Seiberg-Witten invariants of $(\X,\tilde{c})$  are trivial.
\end{Proposition}
\proof
The Seiberg-Witten case is proved by Theorem 1.1.1 in \cite{FC} or Proposition 4.6.5 in \cite{Nicolaescu}.
The $\Pin^-(2)$-monopole case is similar.
Take a tubular neighborhood $N$ of $\Sigma$, and let $Y=\partial N$ and $X_0=\overline{X\setminus N}$.
Then $Y$ admits a positive scalar curvature metric $g_Y$. 
Decompose $X$ as $X=X_0\cup_YN$.
For a positive real number $T$, let us insert a cylinder between $X_0$ and $N$ as:
$$
X_T = X_0\cup([-T,T]\times Y)\cup N.
$$ 
Fix a metric on $X_T$ which is product on the cylinder: $dt^2 + g_Y$.
Let $\alpha_\R$ be the class in $H_2(X_T;\lambda)=H_2(X_T;l\otimes\R)$ corresponding to $\alpha\in H_2(X_T;l)$.
Since $\alpha$ is suuposed to have infinite order, $\alpha_\R$ is a nonzero class in  $H_2(X_T;\lambda)$. 
Let $a\in H^2(X_T;\lambda)$ be the Kronecker dual of $\alpha_\R$ such that $\langle a,\alpha_\R\rangle=1$.    
Then the image of $a$ by the restriction map $r\colon H^2(X_T;\lambda)\to H^2(Y;i^*\lambda)$ is also a nonzero class.
Choose a $2$-form $\eta\in \Omega^2(Y;i^*\lambda)$ representing $r(a)$.
Let us perturb the $\Pin^-(2)$-monopole equations on $Y$ by $\eta$ as in \eqref{eq:Pin2eq3}. 
Since every  $\Pin^-(2)$-monopole solution for a positive scalar curvature metric $g_Y$ is reducible, a generic small choice of $\eta$ makes the perturbed Chern-Simons-Dirac functional \eqref{eq:CSD} have no critical point. 
Choose a $2$-form $\mu\in\ i\Omega^2(X;\lambda)$ whose restriction to the cylinder is the pull-back of $i\eta$.  

Now suppose the $\Pin^-(2)$-monopole invariants of $(X,c)$ is nontrivial. 
Then the moduli space $\M(X_T)$ is nonempty for all $T$.
Taking the limit $T\to \infty$, we can obtain a finite energy solution on the manifold with cylindrical end, $X_0\cup[-1,\infty)\times Y$.
Since a finite energy solution should converge to a critical point at infinity (\corref{cor:limit}), this is a contradiction.
\endproof
%
%
\subsection{The case when $\Sigma$ is nonorientable}\label{subsec:nonori}
%
%
Take a tubular neighborhood $N$ of $\Sigma$, and let $Y=\partial N$ and $X_0=\overline{X\setminus N}$.
Decompose $X$ as $X=X_0\cup_YN$.
For a large $T>0$, insert a long cylinder between $X_0$ and $N$ as:
$$
X_T = X_0\cup([-T,T]\times Y)\cup N.
$$ 
Fix a metric on $X_T$ which is product on the cylinder: $dt^2 + g_Y$.
(Below, we will take a special metric $g_Y$ on $Y$.)
Let $\tilde{X}_T$ be the associated double covering.
Then 
$$
\tilde{X}_T=\tilde{X}_0\cup([-T,T]\times\tilde{Y})\cup \tilde{N},
$$
where $\tilde{Y}=S^1\times \tilde{\Sigma}$ and $\tilde{N}=D^2\times\tilde{\Sigma}$.
(The object with $\tilde{\ }$ is the associated double covering.)
Take the metric $g_Y$ on $Y$ so that its pull-back metric on $\tilde{Y}=S^1\times \tilde{\Sigma}$ is of the form 
$$
d\theta^2 + g_{\tilde{\Sigma}},
$$
where $g_{\tilde{\Sigma}}$ is the metric with constant scalar curvature $-2\pi(4g(\tilde{\Sigma})-4)$.
Then the volume of $\tilde{\Sigma}$ is $1$.

Now, consider the limit $T\to \infty$.
For $\X_T$, the following is known.
\begin{Proposition}[\cite{KM}, Proposition 8]
If the Seiberg-Witten invariant of $(\X,\tilde{c})$ is nontrivial, then 
there is a translation invariant Seiberg-Witten solution on $\R\times\tilde{Y}$.
\end{Proposition}
The same method of proof as in \cite{KM} yields the following:
\begin{Proposition}\label{prop:translation}
If the $\Pin^-(2)$-monopole invariant of $(X,c)$ is nontrivial, then 
there exists a translation invariant $\Pin^-(2)$-monopole solution on $\R\times Y$. 
\end{Proposition}
Under the situation of \propref{prop:translation}, by pulling back the $\Pin^-(2)$-monopole solution on $\R\times Y$ to $\R\times\tilde{Y}$, we also have a translation invariant Seiberg-Witten solution on $\R\times\tilde{Y}$.

By the argument in \cite{KM}, the existence of a translation invariant solution on $\R\times\tilde{Y}$ implies
$$
-\chi(\tilde{\Sigma})\geq |c_1(L)[\tilde{\Sigma}]|,
$$
where $L$ is the determinant line bundle of the $\Spinc$-structure $\tilde{c}$. 
This immediately implies 
$$
-\chi(\Sigma) \geq |\tilde{c}_1(E)[{\Sigma}]|.
$$

%
%
\subsection{The case when $\Sigma$ is orientable}\label{subsec:nori}
%
%
Since the restriction of the local system $l$ to $\Sigma$ is trivial for orientable $\Sigma$, the restrictions of the $\Spincm$-structure to $Y$ and $N$ are untwisted. 
This reduces the argument to the Seiberg-Witten case \cite{KM}.
Let us consider the case when the Seiberg-Witten invariant of $(\X,\tilde{c})$ is nontrivial. 
Since $\Sigma$ is orientable, $\tilde{\Sigma}$ has two components: $\tilde{\Sigma}=\tilde{\Sigma}_1\cup\tilde{\Sigma}_2$.
Then take a tubular neighborhood $\tilde{N}_1$ of $\tilde{\Sigma}_1$, and let $\tilde{Y}_1=\partial\tilde{N}_1$ and $\tilde{X}_0=\overline{\X\setminus\tilde{N}_1}$.
Let us consider 
$$
\X^\prime_T = \tilde{X}_0\cup ([-T,T]\times \tilde{Y}_1)\cup \tilde{N}_1,
$$
for large $T$.
This also reduces the argument to the Seiberg-Witten case \cite{KM}.

\proof[Proof of \corref{cor:genus-cover}]
Since $(\iota_*)^2=\id$,  $H_2(\X;\Q)$ splits into $(\pm 1)$-eigenspaces.
Then $(-1)$-eigenspace is identified with $H_2(X;l\otimes\Q)$. 
Let $\pi\colon\X\to X$ be the projection. 
Then $\pi_*\colon H_2(\X;\Q)\to H_2(X;l\otimes\Q)$ can be identified with $\alpha\mapsto \frac12(\alpha -\iota_*\alpha)$.
It follows from these and the assumption that $\Sigma\cap\iota\Sigma=\emptyset$ that $\pi(\Sigma)$ satisfies the conditions in \thmref{thm:genus-ap}.
\endproof

\renewcommand{\thesection}{\Alph{section}}
\setcounter{section}{0}
%
%
\section{Appendix}\label{sec:appendix}
%
%
The purpose of this appendix is to give a proof of \thmref{thm:gluing}.
The proof is based on \cite{DK}, \S7.2 and \cite{DF}, Chapter 4.
%
%
\subsection{The construction of the map $\mathfrak{F}$}\label{subsec:gluing}
%
%
First, we give the construction of the map $\mathfrak{F}$ of \eqref{eq:I}.
Let $x_i=(A_i,\Phi_i)$ be finite energy monopole solutions on $X_i$ ($i=1,2$).
Fix a gluing parameter $\sigma_0\in\Gamma_\theta$.
Let $x_0^\prime=x^\prime(\sigma_0)$ be the spliced configuration as in \subsecref{subsec:gluing0}.
The goal is to find a true solution $x(\sigma_0)$ near $x_0^\prime$ under the assumptions $H^2_{x_1}=H^2_{x_2}=0$. 

The monopole map on $X^{\#T}$ is defined as a map between weighted spaces:
$$
\Theta\colon L_{k}^{2,w_T}\to L_{k-1}^{2,w_T}.
$$
Since the monopole solutions $x_1$ and $x_2$ decay exponentially(\subsecref{subsec:expdecay}), we have an estimate
$$
\|\Theta(x^\prime_0)\|_{L_{k-1}^2}=O(e^{-\delta_0 T}).
$$
Therefore we also have
\begin{equation}\label{eq:estTheta}
\|\Theta(x^\prime_0)\|_{L_{k-1}^{2,w_T}}=O(e^{(\alpha-\delta_0) T}).
\end{equation}
Now we assume $\alpha<\delta_0$ so that the quantity above will be small for large $T$, and set
\begin{equation}\label{eq:delta}
\delta = \delta_0-\alpha.
\end{equation}
We want to solve the equation for $y=(a,\phi)\in L_k^{2,w_T}(\Lambda^1(i\lambda)\oplus S^+)$ 
\begin{equation*}
\Theta(x^\prime_0+y)=0.
\end{equation*}
This equation is equivalent to 
\begin{equation}\label{eq:tosolve}
(\D_{x_0^\prime} + n)(y) = -\Theta(x_0^\prime),
\end{equation}
where $\D_{x_0^\prime}$ is the linearization of $\Theta$, and $n$ is the quadratic term:
$$
n(y)=(q(\phi),\rho(a)\phi).
$$
To solve \eqref{eq:tosolve}, we first solve the linear version of it.
For this purpose, we construct the right inverse $Q$ of the linear operator $\D_{x_0^\prime}$.
The operator $Q$ will be constructed by splicing the right inverses $Q_1$ and $Q_2$ for the linearizations of $\Theta$, $\D_{x_1}$ and $\D{x_2}$, over $X_1$ and $X_2$.
First, we have the following:
\begin{Proposition}[\cite{DF}, \S3.3] \label{prop:Qi}
For each monopole solutions $x_i$ $(i=1,2)$, if $H^2_{x_i}=0$, then there exists the right inverse $Q_i$ for $D_{x_i}${\rm :}
That is, there exists a map $Q_i\colon L_k^{2,w_i}\to L^{2,w_i}_{k-1}$ and a constant $C_i$ which satisfy{\rm :}
$$
D_{x_i}\circ Q_i(u)=u,\quad \|Q_iu\|_{L^{2,w_i}_k}\leq C_i\|u\|_{L^{2,w_i}_{k-1}},
$$
for every $u \in L^{2,w_i}_{k-1}$.
\end{Proposition}
The proof is a simple adaptation of the argument due to Donaldson \cite{DF}, \S3.3.

Now, an approximate inverse $Q^\prime$ for $D_{x_0^\prime}$ is constructed by splicing as follows:
Recall $X^{\#T}$ is considered as the union $X^{\#T}= X_1^T\cup X_2^T$.
Let $\chi_1:X^{\#T}\to\R$ be the characteristic function of $X_1^T$, that is 
$$
\chi_1(x)=\left\{
\begin{aligned}
1 & \quad x\in X_1^T,\\ 
0&  \quad x\in X^{\#T}\setminus X_1^T.
\end{aligned}\right.
$$
Choose the function $\gamma_1$ such that
\begin{itemize}
\item $\gamma_1=1$ over the support of $\chi_1$,
\item the support of $\gamma_1$ is in $X^0_1\cup [-T,T]\times Y$, and
\item $|\nabla\gamma_1|=O(T^{-1})$.
\end{itemize} 
Take $\chi_2$ and $\gamma_2$ symmetrically.
Then we have $\gamma_1\chi_1 + \gamma_2\chi_2=1$ everywhere.
Now define 
$$
Q^\prime(u)=\gamma_1Q_1(\chi_1u)+\gamma_2Q_2(\chi_2u).
$$ 
Note that the $w_T$-norm of $\chi_1u$ is equal to the $w_1$-norm of that since the weight functions are equal on its support.
Thus we have
$$
\|Q_1(\chi_1u)\|_{L^{2,w_T}_k}\leq C_1\|u\|_{L^{2,w_T}_{k-1}},
$$
where $C_1$ is the constant in \propref{prop:Qi}.

Recall that $x_1$ and $x_2$ are  monopole solutions on $X_1$ and $X_2$, $x_1^\prime$ and $x_2^\prime$ are configurations flattened on the ends and $x_0^\prime$ is the spliced configuration on $X^{\#T}$.
Therefore the linearization $\D_{x_0^\prime}$ is equal to $\D_{x_i^\prime}$ on the support of $\gamma_i$.
Then
\begin{align*}
D_{x_0^\prime}Q^\prime u=&D_{x_0^\prime}(\gamma_1Q_1(\chi_1u)+\gamma_2Q_2(\chi_2u))\\
=& \gamma_1\D_{x_1^\prime}Q_1(\chi_1 u)+ \gamma_2\D_{x_2^\prime}Q_2(\chi_2 u) + \nabla\gamma_1*Q_1(\chi_1 u)+ \nabla\gamma_2*Q_2(\chi_2 u),
\end{align*}
where $*$ means an algebraic multiplication.
Let us estimate each term of the last equation.
The $w_T$-norm of $\nabla\gamma_1*Q_1(\chi_1 u)$ is less than $w_1$-norm of it since $w_T$ is smaller than $w_1$. 
Therefore
$$
\|\nabla\gamma_1*Q_1(\chi_1 u)\|_{L^{2,w_T}_k}\leq \text{const.}T^{-1}\|u\|_{L^{2,w_T}_{k-1}}.
$$
Next we want to estimate $\gamma_1\D_{x_1^\prime}Q_1(\chi_1 u)$.
The operator $Q_1$ is not the right inverse for $\D_{x_1^\prime}$, but is that of   $\D_{x_1}$.
Since $x_1$ decay exponentially, the operator norm of the difference of these two is estimated as
$$
\|\D_{x_1}-\D_{x_1^\prime}\|_{\text{OP}}=O(e^{-\delta T}).
$$ 
Then 
$$
\|(\id-\D_{x_1^\prime}\circ Q_1)u\|_{L^{2,w_1}_{k-1}}\leq \|\D_{x_1}-\D_{x_1^\prime}\|_{\text{OP}}\cdot\|Q_1u\|_{L^{2,w_1}_k}\leq \text{const.}e^{-\delta T}\|u\|_{L^{2,w_1}_{k-1}}.
$$
Summing up these, we obtain
$$
\|(\id - \D_{x_0^\prime}\circ Q)u \|_{L^{2,w_T}_{k-1}}\leq \text{const.}(e^{-\delta T}+T^{-1})\|u\|_{L^{2,w_T}_{k-1}}\leq CT^{-1}\|u\|_{L^{2,w_T}_{k-1}},
$$
for some constant $C$.
If we take a large $T$ so that $CT^{-1}<1$, we obtain the inverse $(\D_{x_0^\prime}Q^\prime)^{-1}$ by iteration. 
Then the true right inverse $Q$ for $\D_{x_0^\prime}$ is given by
$$
Q=Q^\prime(\D_{x_0^\prime}Q^\prime)^{-1}.
$$
For summary,
\begin{Proposition}
There exists the operator $Q\colon L^{2,w_T}_{k-1}\to L^{2,w_T}_k$ which satisfies
\begin{equation}\label{eq:C}
(\D_{x_0^\prime}\circ Q) u=u, \quad \|Qu\|_{L^{2,w_T}_k}\leq C\|u\|_{L^{2,w_T}_{k-1}},
\end{equation}
for every $u\in L^{2,w_T}_{k-1}$.
\end{Proposition}
Now we begin to seek the solution for \eqref{eq:tosolve}.
The main tool for this is the contraction mapping principle.
We seek the solution of the form $y=Qu$. 
So to solve is 
\begin{equation}\label{eq:tosolve2}
u+n(Qu) = -\Theta(x_0^\prime).
\end{equation}
Set $U=L^{2,w_T}_{k-1}$ and $V=L^{2,w_T}_k$. 
Since $n$ is a quadratic map, we have an estimate that there is a constant $M$ such that
\begin{equation}\label{eq:M}
\|n(y_1)-n(y_2)\|_U\leq M\|y_1-y_2\|_V(\|y_1\|_V+\|y_2\|_V),
\end{equation}
for any $y_1$, $y_2$ in $V$.
Then we have
$$
\|n(Qu_1)-n(Qu_2)\|_U\leq MC^2 \|u_1-u_2\|_U(\|u_1\|_U+\|u_2\|_U),
$$
where the constants $C$ and $M$ are those in \eqref{eq:C} and \eqref{eq:M}.
Now if, for instance, $\|\Theta(x_0^\prime)\|_U\leq (100MC^2)^{-1}$, then there exists a unique solution to \eqref{eq:tosolve2}.
By \eqref{eq:estTheta}, $\|\Theta(x_0^\prime)\|_U$ can be arbitrary small if we take a sufficient large $T$.
Thus for large $T$, we have a unique solution $u$. 
Let $y=Qu$. 
Then $x_0^\prime + y$ is a required monopole solution which is in $L^{2,w_T}_k$, and hence in $C^\infty$.

Thus for each $\sigma\in\Gamma$, we can find a monopole solution $x(\sigma)$ in a unique way near the spliced configuration $x^\prime(\sigma)$.
The correspondence $\sigma\mapsto x(\sigma)$ descends to the map $\mathfrak{F}$. 
%
%
\subsection{The image of $\mathfrak{F}$}\label{subsec:gluingpara}
%
%
%
%
We would like to characterize the image of $\mathfrak{F}$.
%
Let $d$ be the metric on $\B(X^{\#T})$ given by  
$$
d([x],[y])=\inf_{g\in\G}\|x-gy\|_V,
$$
where $V=L^{2,w_T}_{k-1}$.
For $\varepsilon>0$, let $\UU(\varepsilon)\subset\B(X^{\#T})$ be the open set
\begin{equation}\label{eq:Ueps}
\UU(\varepsilon)=\{[x]\,|\,d([x],\mathfrak{F}^\prime(\Gamma_\theta))<\varepsilon, \|\Theta(x)\|_U<\varepsilon\}.
\end{equation}
\begin{Proposition}\label{prop:continuity}
If $H_{x_i}^0=H_{x_i}^1=H_{x_i}^2=0$ for $i=1,2$, then for small enough $\varepsilon$ there exists $T(\varepsilon)$ so that for $T>T(\varepsilon)$ any point in $\UU(\varepsilon)$ can be represented by a configuration of the form $x^\prime(\sigma)+Q_\sigma u$ with $\|u\|_U<\text{\rm const.}\varepsilon$, where $Q_\sigma$ is the right inverse for $D_{x^\prime(\sigma)}$.
\end{Proposition}
Assuming the proposition, we have
\begin{Corollary}\label{cor:continuity}
If $H_{x_i}^0=H_{x_i}^1=H_{x_i}^2=0$ for $i=1,2$, then for $\varepsilon$ and $T(\varepsilon)$ in \propref{prop:continuity}, and for every $T>T(\varepsilon)$, the intersection $\UU(\varepsilon)\cap \M(X^{\#T})$ is equal to the image of $\mathfrak{F}\colon\Gamma_\theta\to\M$.
\end{Corollary}
The corollary follows from the argument in \subsecref{subsec:gluing} since under the given assumptions there is a unique small solution $u$ to the equation $\Theta(x^\prime(\sigma)+Q_\sigma u)=0$.
%
%
\subsection{Proof of \propref{prop:continuity}: Closedness}
%
%
We will prove \propref{prop:continuity} by continuity method.
Let $[y]$ be an element of $\UU(\varepsilon)$. 
Then there exists $x^\prime\in \mathfrak{F}^\prime(\Gamma_\theta)$ with $\|x^\prime - y\|_V<\varepsilon$.
Let us write $y=x^\prime + b$ and consider the path for $t\in[0,1]$,
$$
y_t=x^\prime + t b.
$$
By gauge transformation, we may assume $y$ and $b$ are smooth, and so is $y_t$ for all $t$. 
It can be seen that, for given $\varepsilon$, if we take $T$ large enough, the class $[y_t]$ is in $\UU(\varepsilon)$ for every $t\in[0,1]$.
Let us define the subset $S\subset [0,1]$ as the set of $t$ which has the property that there exist $g_t\in \G$, $x_{\sigma_t}^\prime\in \mathfrak{F}^\prime(\Gamma_\theta)$ and $u_t\in U=L^{2,w_T}_{k-1}$ such that
\begin{equation}\label{eq:gtyt}
g_ty_t = x_{\sigma_t}^\prime + Q_{\sigma_t}(u_t),
\end{equation}
with $\|u_t\|_U<\nu$, where $\nu$ will be chosen below.
Obviously $0\in S$.
We prove $S$ is open and closed.

Let us prove the closedness. 
Suppose $t\in S$. 
Then there exist $g_t$, $x_{\sigma_t}^\prime$  and $u_t$ so that \eqref{eq:gtyt} holds.
Applying $\Theta$ on both sides of \eqref{eq:gtyt}, we have
\begin{equation}\label{eq:Thetagtyt}
\Theta(g_ty_t) = \Theta(x_{\sigma_t}^\prime) + \D_{\sigma_t}Q_{\sigma_t}u_t + n(Q_{\sigma_t}u_t) = \Theta(x_{\sigma_t}^\prime) + u_t + n(Q_{\sigma_t}u_t). 
\end{equation}
Then we have an estimate
$$
\|u_t\|_U\leq\|\Theta(y_t)\|_U + \|\Theta(x^\prime({\sigma_t}))\|_U + \|n(Q_{\sigma_t}u_t)\|_U \leq \varepsilon + \text{const.}e^{-\delta T} + (C_{\sigma_t})^2\|u_t\|^2_U,
$$
where $C_{\sigma_t}$ is the constant for $Q_{\sigma_t}$ so that $\|Q_{\sigma_t}u\|\leq C_{\sigma_t}\|u\|$.
Since $\Gamma_\theta$ is compact, $C_{\sigma}$ $(\sigma\in\Gamma_\theta)$ is bounded above by some constant $N$ as
\begin{equation}\label{eq:N}
C_{\sigma}\leq N.
\end{equation}
Rearranging this and taking $\nu$ so that $\|u\|\leq\nu\leq (2N^2)^{-1}$, we have
\begin{equation}\label{eq:ut}
\frac12\|u_t\|_U\leq(1-N^2\|u_t\|)\|u_t\|\leq \varepsilon + \text{const.}e^{-\delta T}.
\end{equation}
This estimate implies the following:
\begin{Lemma}\label{lem:open-closed}
Suppose $\nu\leq(2N^2)^{-1}$ so that the estimate \eqref{eq:ut} holds.
Then we can choose small $\varepsilon$ and large $T$ so that $\|u\|<\nu$ implies $\|u\|\leq\frac12\nu$. 
\end{Lemma}
Thus the open condition $\|u\|<\nu$ is also closed.

Suppose we have $t_i\in S$ with $t_i\to t_\infty$.
By definition, for each $t_i$, there exist $u_i=u_{t_i}$, $\sigma_i=\sigma_{t_i}$ and $y_i=y_{t_i}$, and if we set $x_i=x^\prime(\sigma_{t_i}) + Q_{\sigma_{t_i}} u_i$, then $g_iy_i=x_i$ holds.
Then obviously $y_{i}=x^\prime +t_ib$ converge to $y_\infty=x^\prime +t_\infty b$ in $C^\infty$.
Since $\Gamma_\theta$ is compact, $\sigma_i$ converge to some $\sigma_{\infty}$. 
Then the spliced configurations $x^\prime(\sigma_{t_i})$ converge to $x^\prime(\sigma_{t_\infty})$ in $C^\infty$.
By the uniform bound $\|u_i\|_U<\nu$, taking a subsequence, we have a weak limit $u_\infty$ so that $u_i\to u_\infty$ in $U=L^{2,w_T}_{k-1}$.
Then $x_i$ converge weakly in $L^{2,w_T}_k$, and we may assume $u_i$ converge weakly in $L^{2,w_T}_{k+1}$ and strongly in $L^{2,w_T}_k$.
Now we would like to see that $u_i$ converge to $u_\infty$ strongly.
By \eqref{eq:Thetagtyt}, 
\begin{equation}\label{eq:uiuj}
\|u_i-u_j\|_U\leq \|\Theta(g_iy_i)-\Theta(g_jy_j)\|_U + \|\Theta(x^\prime(\sigma_i))- \Theta(x^\prime(\sigma_j))\|_U +\|n(Q_{\sigma_i}u_i)- n(Q_{\sigma_j}u_j)\|_U.
\end{equation}
If $i,j\to \infty$, then the second term of the right hand side $\|\Theta(x^\prime(\sigma_i))- \Theta(x^\prime(\sigma_j))\|_U$ tends to $0$, because $x^\prime(\sigma_i)$ converge in $C^\infty$.
The first term is estimated, for instance, as
$$
\|\Theta(g_iy_i)-\Theta(g_jy_j)\|_U\leq \|g_i-g_j\|_{L^{2,w_T}_{k-1}}\cdot\|\Theta(y_j)\|_{C^0} + \|g_i\|_{L^{2,w_T}_{k-1}}\cdot\|\Theta(y_i)-\Theta(y_j)\|_{C^0},
$$
where the right hand side tends to $0$ if $i,j\to \infty$.
For the third term,
\begin{align*}
\|n(Q_{\sigma_i}u_i)- n(Q_{\sigma_j}u_j)\|_U\leq& M\|Q_{\sigma_i}u_i - Q_{\sigma_j}u_j\|_U(\|Q_{\sigma_i}u_i\|_U+\|Q_{\sigma_j}u_j\|_U)\\
\leq & M(\|Q_{\sigma_i} - Q_{\sigma_j}\|_{\text{OP}}\|u_i\|_U + C_{\sigma_j}\|u_i-u_j\|_U)(C_{\sigma_i}\|u_i\|_U+C_{\sigma_j}\|u_j\|_U)\\
\leq & M(\|Q_{\sigma_i} - Q_{\sigma_j}\|_{\text{OP}}\|u_i\|_U + N\|u_i-u_j\|_U)\cdot N(\|u_i\|_U+\|u_j\|_U).
\end{align*}
where $N$ is the constant in \eqref{eq:N}.
If we assume $\|u_i\|_U<(4MN^2)^{-1}$, then, by rearranging \eqref{eq:uiuj}, we can see that the sequence $\{u_i\}$ is a Cauchy sequence in $U$, and $u_\infty$ is the strong limit.

Now we choose $\nu$ so that $\nu\leq \min\{(2N)^{-1}, (4MN^2)^{-1} \}$, and choose $\varepsilon$ and $T$ as in \lemref{lem:open-closed}.
Then $\{u_i\}$ converge strongly to $u_\infty$, and by \lemref{lem:open-closed}, the limit $u_\infty$ satisfies $\|u_\infty\|_U<\nu$. 
This means $t_\infty\in S$, and the closedness is proved.

%
%
\subsection{Proof of \propref{prop:continuity}: Openness}
%
%
Let us prove the openness.
To prove the openness, we use the implicit function theorem.
Suppose $t_0\in S$ with $0\leq t_0<1$ so that there exist $g_0$, $\sigma_0$ and $u_0$ so that $g_0y_{t_0} = x^\prime_0 + Q_0u_0$. 
To prove is $[t_0,t_0+\epsilon)\subset S$ for small $\epsilon$. 
In fact, we will prove any configuration $z$ close to $y_{t_0}$ is gauge equivalent to some $x^\prime_v + Q_v(u_0+w)$ for some $v\in\Lie\Gamma_\theta$ and $w\in U$, where $x^\prime_v=x^\prime(\sigma_0,v)$ and $Q_v=Q_{(\sigma_0\exp v)}$.
Define a map
$$
{\cal F}\colon \Omega^0(i\lambda)\times\Lie\Gamma_\theta\times(\Omega^+(i\lambda)\oplus\Gamma(S^-))\to \Omega^1(i\lambda)\oplus\Gamma(S^+)
$$
by
$$
{\cal F}(f,v,w) = \exp(f)(x_v^\prime + Q_v(u_0+w)) - (x^\prime_0 + Q_0(u_0)).
$$
We need to show that ${\cal F}$ is surjective onto a neighborhood of $0$. This follows from the implicit function theorem if the derivative $D{\cal F}$ of ${\cal F}$ at $(0,0,0)$ is surjective. 
If $x^\prime_0=(A^\prime_0,\Phi^\prime_0)$, the derivative $D{\cal F}$ is 
$$
D{\cal F}_{(0,0,0)}(f,v,w) = \I_{\Phi^\prime_0}(f) + \partial_vx_v^\prime + \partial_vQ_v(u_0) + Q_0(w), 
$$
where $\I_{\Phi^\prime_0}(f)=(-df,f\Phi^\prime_0)$ and $\partial_v$ means the derivative with respect to $v$.

More precisely, $\partial_vx^\prime_v$ can be written as follows:
For the connection part $A^\prime(\sigma_0,v)$ of $x^\prime_v$, set
$$
j(v) = \left.\frac{\partial}{\partial s}A^\prime(\sigma_0,sv)\right|_{s=0}\in \Omega^1(i\lambda).
$$
Then $j(v)=d(\lambda_2 v)=-d(\lambda_1 v)$ on $[-1,1]\times Y$, and $j(v)=0$ outside of $[-1,1]\times Y$.
Now we have
$$
\partial_vx^\prime_v = (j(v),0)\in \Omega^1(i\lambda)\oplus\Gamma(S^+).
$$ 
The term $\partial_vQ_v$ will be discussed below.

In order to prove the surjectivity of $D{\cal F}$, we define a map 
$$
{\cal T} \colon \Omega^0(i\lambda)\times\Lie\Gamma_\theta\times(\Omega^+(i\lambda)\oplus\Gamma(S^-))\to \Omega^1(i\lambda)\oplus\Gamma(S^+)
$$
by 
$$
{\cal T}(f,v,w) = \I_{\Phi^\prime_0}(f) + (j(v),0) + Q_0(w) =D{\cal F} - \partial_vQ_v(u_0).
$$
Let $B_1$ be the completion of the domain of ${\cal T}$ in the norm:
$$
\|(f,v,w)\|_{B_1}=\|\I_{\Phi^\prime_0}(f) + j(v)\|_V + \|w\|_U,
$$
where $U=L^{2,w_T}_{k-1}$ and $V=L^{2,w_T}_k$.
This is a norm by \lemref{lem:C0} below.
Let $B_2$ be the completion of the range in the norm:
$$
\|(a,\phi)\|_{B_2} = \|(a,\phi)\|_V + \|\D_{x^\prime_0}(a,\phi)\|_U.
$$

Now, the fact that $\|\cdot\|_{B_1}$ is a norm follows from the following:
\begin{Lemma}\label{lem:C0}
If $H_{x_1}^0=H_{x_2}^0=0$, then there exists a constant $L$ independent of $T$ such that, for any $f\in\Omega^0(i\lambda)$ and any $v\in\Lie\Gamma_\theta$, we have
$$
\|f\|_{C^0} + |v| \leq L\|\I_{x^\prime_0}(f) + j(v)\|_V.
$$ 
\end{Lemma}
\proof
Let $f_1 = f + (1-\lambda_1)$ over $\supp(\lambda_1)\subset X_1$, and 
$f_2 = f - (1-\lambda_2)$ over $\supp(\lambda_2)\subset X_2$.
Then, for $i=1,2$,
$$
df+j(v) = df_i 
$$
over $X_i^T$, and $f_1-f_2=v$ over $[-1,1]\times Y$.
Then each of $\|f\|_{C^0}$ and $|v|$ is bounded above by $\|f_1\|_{C^0}+\|f_2\|_{C^0}$.
On the other hand, if $H_{x_i}^0=0$, then there exists a constant $L_i$ for each $i$ so that 
$$
\|f_i\|_{C^0}\leq L_i\|\I_{\Phi_i}(f_i)\|_{L^{2,w_i}_k}.
$$
Then we have
$$
\|f_i\|_{C^0}\leq L_i\|\I_{\Phi^\prime_i}(f_i)\|_{L^{2,w_i}_k} + \|\Phi_i - \Phi^\prime_i\|_{L^{2,w_i}_k}\|f_i\|_{C^0}.
$$
By the exponential decay, $\|\Phi_i - \Phi^\prime_i\|_{L^{2,w_i}_k}=O(e^{-\delta T})$. 
So we can choose large $T$ so that $\|\Phi_i - \Phi^\prime_i\|_{L^{2,w_i}_k}<\frac12$, say.
Rearranging this, we obtain a bound for $\|f\|_{C^0}$ by $\|\I_{\Phi^\prime_i}(f)\|_{L^{2,w_i}_k}$.
Since $\I_{\Phi^\prime_i}(f_i)$ is supported on $\supp\lambda_i$, the $L^{2,w_i}_k$ and $V=L^{2,w_T}_k$ norms of it are uniformly equivalent, and the lemma is proved.
\endproof
Thus ${\cal T}$ is a bounded map from $B_1$ to $B_2$.
In fact, the following holds:
\begin{Lemma}
There exists a constant $K$ independent of $T$ so that
\begin{equation}\label{eq:K}
\|(f,v,w)\|_{B_1}\leq K\|{\cal T}(f,v,w)\|_{B_2}.
\end{equation}
\end{Lemma}
\proof
Let $\alpha={\cal T}(f,v,w)=\I_{\Phi^\prime_0}(f)+j(v)+Q_0(w)$.
We consider $\D_{x^\prime_0}\alpha$. 
By \eqref{eq:DI}, we have
$$
\D_{x^\prime_0}\I_{\Phi^\prime_0}(f)=(0,fD_{A^\prime_0}\Phi^\prime_0).
$$
On the other hand, 
$$
\D_{x^\prime_0}(j(v))=(-d(j(v)), j(v)\Phi^\prime_0)=0,
$$ because $\supp j(v)\cap \supp\Phi^\prime_0=\emptyset$.
Thus we have $$
\D_{x^\prime_0}\alpha=(0,fD_{A^\prime_0}\Phi^\prime_0) + w.
$$
Since $\|f\|_{C^0}\leq L\|\I_{x^\prime_0}(f) + j(v)\|_V$ by \lemref{lem:C0}, and $\|D_{A^\prime_0}\Phi^\prime_0\|_U=O(e^{-\delta T})$, we obtain
\begin{align*}
\|w\|_U\leq & \|\D_{x^\prime_0}\alpha\|_U + \|f\|_{C^0}\|D_{A^\prime_0}\Phi^\prime_0\|_U \\
\leq & \|\D_{x^\prime_0}\alpha\|_U + \text{const.}e^{-\delta T}\|\I_{x^\prime_0}(f) + j(v)\|_V\\
=& \|\D_{x^\prime_0}\alpha\|_U + \text{const.}e^{-\delta T} \|\alpha-Q_0(w)\|_V\\
\leq & \|\D_{x^\prime_0}\alpha\|_U + \text{const.}e^{-\delta T} (\|\alpha\|_V + C\|w\|_U ).
\end{align*}
Thus, when $T$ is sufficiently large, we obtain a bound $\|w\|_U\leq K_1\|\alpha\|_{B_2}$ for some constant $K_1$.
Therefore we have
$$
\|\I_{x^\prime_0}(f) + j(v)\|_V = \|\alpha-Q_0(w)\|_V \leq (1+CK_1)\|\alpha\|_{B_2}.
$$
Combining the last two inequalities, we can find a constant $K$ so that \eqref{eq:K} holds.
\endproof
\begin{Corollary}\label{cor:T}
The kernel of ${\cal T}$ is zero, and the image of ${\cal T}$ is closed in $B_2$.
\end{Corollary}
Now we use the index theorem to prove $\cal T$ is the isomorphism.
\begin{Proposition}
If $H_{x_i}^0=H_{x_i}^1=H_{x_i}^2=0$ for $i=1,2$, then the operator ${\cal T}$ is an isomorphism from $B_1$ to $B_2$ with operator norm $\|{\cal T}^{-1}\|_\text{\rm OP}\leq K$.
\end{Proposition}
\proof
The operator $Q_0$ is a pseudo-differential operator whose symbol is homotopic to $(D_{x^\prime_0})^*(1+(D_{x^\prime_0})^*D_{x^\prime_0})^{-1}.$
Thus $\I_{\Phi^\prime_0}\oplus Q_0$ is Fredholm, and the index is calculated as
$$
\ind \left[\I_{\Phi^\prime_0}\oplus Q_0\right] = \ind \left[(\I_{\Phi^\prime_0})^*\oplus \D_{x^\prime_0}\right]^* = -\ind \left[(\I_{\Phi^\prime_0})^*\oplus \D_{x^\prime_0}\right].
$$
Then
\begin{align*}
\ind{\cal T} =& \dim\Lie\Gamma_\theta-\ind \left[(\I_{\Phi^\prime_0})^*\oplus \D_{x^\prime_0}\right]\\
=& \dim\Lie\Gamma_\theta-\left\{\ind \left[(\I_{\Phi_1})^*\oplus \D_{x_1}\right] + \ind \left[(\I_{\Phi_2})^*\oplus \D_{x_2}\right] + \dim\Lie\Gamma_\theta \right\} =0.
\end{align*}
Now the proposition immediately follows from \corref{cor:T}.
\endproof
Recall $D{\cal F}={\cal T} + \partial Q_v(u_0)$, and we have seen that ${\cal T}$ is an isomorphism from $B_1$ to $B_2$ which satisfies \eqref{eq:K}.
If we see the operator norm of the map $v\mapsto \partial Q_v(u_0)$ is less than $K^{-1}$ in \eqref{eq:K}, then ${\cal F}$ is also invertible, and the proof of \propref{prop:continuity} is completed.

Let us evaluate the norm of $ \partial Q_v(u_0)$.
Recall $Q_v$ is constructed as 
\begin{equation}\label{eq:Qv}
Q_v = Q^\prime_v(\D_{x^\prime_v}Q^\prime_v)^{-1},
\end{equation}
where $Q^\prime_v$ is the spliced operator which can be written as
$$
Q^\prime_v = Q^\prime_{v,1} + Q^\prime_{v,2},
$$
with 
$$
Q^\prime_{v,i}=h_iQ^\prime_ih_i^{-1},\quad\text{ and }\quad Q^\prime_i(u) = \gamma_iQ_i(\chi_iu),
$$
where $h_i$ are the gauge transformations in \eqref{eq:hi}, and $Q_i$, $\gamma_i$ and $\chi_i$ are defined around \propref{prop:Qi}.
Then the differential of $Q_v^\prime$ with respect to $v$ at $v=0$ is given by
$$
\partial_v Q_v^\prime(u) =  [(1-\lambda_1)v_1,Q_1u] + [(1-\lambda_2)v_2,Q_2u].
$$
Then we have
$$
\|\partial_v Q_v^\prime(u)\|_V\leq\text{const.}|v|\|u\|_U.
$$
Similarly, the differential $\partial_v(\D_{x^\prime_v}Q_v^\prime(u))$ is bounded as
$$
\|\partial_v(\D_{x^\prime_v}Q_v^\prime(u))\|_V\leq\text{const.}|v|\|u\|_U.
$$
By differentiating \eqref{eq:Qv}, we obtain
$$
\partial_v Q_v = \{\partial Q_v^\prime - Q_0\partial(\D_{x_v^\prime}Q_v^\prime)\}(\D_{x^\prime_0}Q^\prime_0)^{-1}.
$$
Hence we obtain the estimate
\begin{equation}\label{eq:dQv}
\|\partial_v Q_v(u)\|_V\leq\text{const.}|v|\|u\|_U.
\end{equation}
Differentiating the identity $\D_{x_v^\prime}Q_v=1$, we have
$$
\D_{x^\prime_0}(\partial_vQ_v(u)) = -(\partial_v\D_{x_v^\prime})Q_0(u) = (0,-j(v)\phi),
$$
where $\phi$ is the spinor component of $Q_0(u)$.
Therefore we have the estimate
$$
\|\D_{x^\prime_0}(\partial_vQ_v(u))\|_U\leq\text{const.}\|j(v)\|_{L^4}\|Q_0(u)\|_V\leq \text{const.}|v|\|u\|_U,
$$
because of the following facts:
\begin{itemize}
\item H\"{o}lder's inequality  implies $\|ab\|_U\leq \|a\|_V\|b\|_{L^4}$, and
\item we may assume the $L^4$ norm of $\nabla\lambda_i$ is independent of $T$, and therefore $\|j(v)\|_{L^4}\leq \text{const.}|v|$.
\end{itemize}

Summing up these, we obtain
$$
\|\partial_vQ_v(u)\|_{B_2}\leq \text{const.}|v|\|u\|_U\leq \text{const.}\|(f,v,u)\|_{B_2}\|u\|_U.
$$
Now if $\|u\|_U$ is small (i.e., $\nu$ is small), then $D{\cal F}$ is invertible and the proof of \propref{prop:continuity} is completed.
%
%
\subsection{The injectivity of the map $\mathfrak{F}$}
%
%
Now, we prove that the map $\mathfrak{F}$ is injective:
\begin{Proposition}[\cite{DK}, \S7.2.6]\label{prop:injective}
For monopoles $x_i$ on $X_i$ with $H_{x_i}^2=0$ {\rm (}$i=1,2${\rm )}, and  for sufficiently small $\varepsilon$,  the map $\mathfrak{F}$ of \eqref{eq:I} is injective.
\end{Proposition}
\proof
If $H^0_{x_i}\neq 0$ for some $i$, then $\mathrm{Gl}$ is one point, and therefore $\mathfrak{F}$ is obviously injective.
Suppose  $H^0_{x_i} =  0$  for $i=1,2$.
For the fixed identification $\sigma_0$ and any $v$, suppose the following:
\begin{itemize}
\item $\mathfrak{F}(\sigma_0)$ is represented by $x^\prime(\sigma_0)+y_0$.
\item $\mathfrak{F}(\sigma_0\exp v)$ is represented by $x^\prime(\sigma_0,v)+y_v$.
\item $x^\prime(\sigma_0)+y_0 $ and $x^\prime(\sigma_0,v)+y_v$ are gauge equivalent.
\item $x^\prime(\sigma_0)$ and $x^\prime(\sigma_0,v)$ are not gauge equivalent. 
\end{itemize}
First, we claim that, if $x^\prime(\sigma_0)+y_0 $ and $x^\prime(\sigma_0,v)+y_v$ are gauge equivalent, then we may assume they are equivalent by a gauge transformation in the identity component. 
Recall that $\pi_0\G\cong H^1(X;l)\cong \Z^{b_1(X;l)}\oplus\Z_2$.
Let $\rho\colon \G\to \pi_0\G$ be the projection.
Let us consider the case when both of $x_1$ and $x_2$ are monopoles on twisted $\Spincm$-structures. 
(The untwisted $\Spincm$-cases are easier.)
Let us write the connection terms of $x^\prime(\sigma_0)$ and $x^\prime(\sigma_0,v)$ as $A(\sigma_0)$ and $A(\sigma_0,v)$, and let $a_0$ and $a_v$ be the $1$-form components of $y_0$ and $y_v$. 
By \lemref{lem:S}, there exists $t\in H^1(X;l)$ which is represented by a gauge transformation $\check{g}$ such that $A(\sigma_0,\pi i) = \check{g} A(\sigma_0)$.
Hence, as de Rham classes, $n[t]=[A(\sigma_0,n\pi i) - A(\sigma_0)]$ for $n\in\Z$.
Suppose $x^\prime(\sigma_0)+y_0= g(x^\prime(\sigma_0,v)+y_v)$ for some $g\in\G$. 
Since the de Rham classes of $a_0$ and $a_v$ are very small for large $T$, we see that $[A(\sigma_0)+a_0-(A(\sigma_0,v)+a_v)] $ should be $n[t]$ for some $n\in\Z$, and therefore $\rho(g)$ is in $\Z\langle t \rangle\oplus\Z_2$.
Then by replacing $x^\prime(\sigma_0,v)$ by $(\pm 1)\cdot \check{g}^{-n} \cdot x^\prime(\sigma_0,v)$, we may assume   $x^\prime(\sigma_0)+y_0 $ and $x^\prime(\sigma_0,v)+y_v$ are gauge equivalent by a gauge transformation of the form $g=\exp(\chi)$ for some $\chi\in\Omega^0(i\lambda)$.
By restricting on $X_i^{2T}$, we obtain gauge transformations $g_i$ over $X_i^{2T}$ so that  $x^\prime(\sigma_0)+y_0 = g_i(x^\prime(\sigma_0,v)+y_v)$. 
Then, for the connection parts, we have
$$
A_i^\prime(\sigma_0)+a_0 = g_i(A_i^\prime(\sigma_0,v) + y_v)= 
g_ih_i (A_i^\prime(\sigma_0)+a_v),
$$
where $h_i$ are the gauge transformations in \eqref{eq:hi}.
Set $g_i^\prime = g_ih_i$. 
We may assume $g_i^\prime = \exp(\chi_i)$ for some $\chi_i$.
Then we have $-2d\chi_i=a_0-a_v$, and therefore
$$
\|\chi_i\|_{C^0}\leq \text{const.}\|d\chi_i\|_{L^{2,w_i}_k(X_i^{2T})}\leq \text{const.}\|a_0-a_v\|_{L^{2,w_i}_k(X_i^{2T})}.
$$
On the overlapping region, the compatibility condition for $g_i$ implies $|\chi_1-\chi_2|=|v|$.
Thus we have 
\begin{equation}\label{eq:av1}
|v|\leq \text{const.}\|a_0-a_v\|_V.
\end{equation}

On the other hand, $y_v$ is given as $y_v= Q_v(u_v)$ for a $u_v$ such that $u_v + n(Q_v(u_v)) = -\Theta(x^\prime(\sigma_o,v))$.
Since $ \Theta(x^\prime(\sigma_o,v))$ is supported on the region where $h_i=1$, the $v$-derivative of $u_v$ is given by $\partial_vu_v =-\partial_v(n(Q_v(u_v)))$.
By calculating the derivative (by using \eqref{eq:M}), we have
$$
\|\partial_v(n(Q_v(u_v))) \|_U\leq\text{const.} \|\partial_v(Q_v(u_v))\|_V\cdot\|Q(u)\|_V\leq\text{const.}\|(\partial_vQ_v)(u)+Q(\partial_vu_v)\|_V\cdot\|u\|_U
$$ 
Since $\|u \|_U\leq \text{const.}\varepsilon$, the estimate \eqref{eq:dQv} implies
$$
\|\partial_v u_v\|_U\leq\text{const.}(|v|\varepsilon+\|\partial_v u_v\|_U)\varepsilon
$$
Rearranging this, we have $\|\partial_v u_v\|_U\leq\text{const.}|v|\varepsilon^2$, and hence
$$
|v|\leq \text{const.}\|a_0-a_v\|_V\leq \text{const.}|v|\varepsilon^2.
$$
Thus for small $\varepsilon$, we obtain $v=0$.
\endproof
%
%
\subsection{Proof of \thmref{thm:gluing}}
%
%
Now we prove \thmref{thm:gluing}. 
Suppose the assumptions in \thmref{thm:gluing} are satisfied.
Since $\M_i$ are compact,  we can define for sufficiently large $T$ the global gluing map,
\begin{equation}\label{eq:gluing}
\Xi\colon \tilde{\M}_1\times_{\Gamma_\theta}\tilde{\M}_2\to \M(X^{\#T}).
\end{equation}
We need some more things.
Let $\varepsilon$ and $T(\varepsilon)$ be the constants in \propref{prop:continuity}, and take $T>T(\varepsilon)$.
 For $\tau$ such that $T>\tau>T(\varepsilon)$, let $K_1^\tau=X_1^\tau$, $K_2^\tau=X_2^\tau$ and $K^\tau =K_1^\tau\cup K_2^\tau$.
We can assume $K^\tau$ as a submanifold of $X^{\#T}$.  
So by restricting to $K^\tau$, we can compare  configurations on the different manifolds $X_1\cup X_2$ and $X^{\#T}$.
Let $\B(K^\tau)$ be the space of the configurations modulo gauge over $K^\tau$. 
We may identify $\B(K^\tau)=\B(K_1^\tau)\times\B(K_2^\tau)$.
For $a=[x_1]\times[x_2]$ and $b=[y_1]\times[y_2]$ in $\B(K^\tau)$, we define the metric 
$$
d_{K^\tau}(a,b)=\inf_{g_1\in\G(K_1^\tau)}\|g_1x_1-y_1\|_V + \inf_{g_1\in\G(K_2^\tau)}\|g_2x_2-y_2\|_V.
$$
For monopoles $x_i$ ($i=1,2$) on $X_i$, let $x_i^\prime$ be the flattened configuration, and $\mathfrak{F}^\prime\colon \mathrm{Gl}\to \B(X^{\#T})$  the map splicing $x_1^\prime$ and $x_2^\prime$ with a gluing parameter $\sigma$. 
If $w$ is a monopole on $X^{\#T}$, then there exists a constant $C$ such that,
$$
d_{K^\tau}([w]|_{K^\tau},[x_1^\prime]\times[x_2^\prime]|_{K^\tau})< C d([w],\mathfrak{F}^\prime(\mathrm{Gl})).
$$
Conversely, we have
\begin{Proposition}\label{prop:tau}
There exists a constant $\tau$ with $\tau>T(\varepsilon)$ such that if 
$$
d_{K^\tau}([w]|_{K^\tau},[x_1^\prime]\times[x_2^\prime]|_{K^\tau})< \frac{\varepsilon}2, 
$$
then
$$
d([w],\mathfrak{F}^\prime(\mathrm{Gl}))<\varepsilon.
$$
\end{Proposition}
\proof
Let us consider the disjoint union $K^T = K_1^T\cup K_2^T$, where $K_i^T=X_i^T$.
Note that $X^{\#T}$ is made by gluing $K_1^T$ and $K_2^T$.
For a monopole $x$ on $X^{\#T}$, let us consider the restriction
$[x]|_{K^T} = [x|_{K_1^T}]\times[x|_{K_2^T}]\in\B(K_1^T)\times\B(K_2^T)$.
Then $d_{K^T}([x]|_{K^T},[x_1^\prime]\times[x_2^\prime]|_{K^T})<\varepsilon$ implies $d([x],\mathfrak{F}^\prime(\mathrm{Gl}))<\varepsilon$.
Let $K^{T-\tau}$ be the disjoint union of $\overline{(X_1^T\setminus X_1^\tau)}$ and $\overline{(X_2^T\setminus X_2^\tau)}$.
The exponential decay estimate implies that there exists a constant $C$ such that, for every monopole $w$ on $X^{\#T}$ and every monopoles $x_1$ and $x_2$ on $X_1$ and $X_2$, 
$$
d_{K^{T-\tau}}([w],[x_1^\prime]\times[x_2^\prime])<Ce^{-\delta\tau}.
$$
Hence if $\tau$ is large enough, then, say, 
$d_{K^{T-\tau}}([w],[x_1^\prime]\times[x_2^\prime])<\varepsilon/10$,
and the proposition holds.
\endproof
By \corref{cor:continuity}, we obtain the following:
\begin{Corollary}
For $\tau$ in \propref{prop:tau}, if $w$ is a monopole on $X^{\#T}$ with $d_{K^\tau}([w]|_{K^\tau},[x_1^\prime]\times[x_2^\prime]|_{K^\tau})< {\varepsilon}/2$ as above, then $[w]$ is in the image $\mathfrak{F}(\mathrm{Gl})$. 
\end{Corollary}
In order to define the inverse of the gluing map, we need to make monopoles on $X_i$ from a monopole on $X^{\#T}$. 

Suppose $x$ is a monopole on $X^{\#T}$ with $H^2_x=0$.
Let us consider the configuration $x^\prime$ obtained by making $x$ flattened on the neck.
More precisely, using the function $\gamma$ in \subsecref{subsec:gluing}, we define the function $\bar{\gamma}$ by
$$
\bar{\gamma} (t) =\left\{
\begin{aligned}
\gamma(-t-3),&\quad t\geq 0,\\
\gamma(t+3),\ &\quad t<0,
\end{aligned}\right.
$$  
and let 
$$
x^\prime = \bar{\gamma} x + (1-\bar{\gamma})(\theta,0).
$$
For each $i$, restricting $x^\prime$ to $X_i^{T}$, assuming $X_i^T\subset X_i$ and extending $x^\prime$ over $X_i$ obviously, we obtain an approximate monopole $x^\prime_i$ on each $X_i$.
Taking a large $T$ and arguing as in \subsecref{subsec:gluing}, we can construct, for each $i=1,2$, a genuine monopole $y_i$ on $X_i$ which is close to $x_i^\prime$. 
To do this, first we need to construct a right inverse $Q_i$ for the operator $\D_{x_i^\prime}$ for each $i$.
The operator $Q_1$, say, is constructed by splicing the right inverse $Q_{x^\prime}$ for $\D_{x^\prime}$ over $X^{\#T}$ with the right inverse for the operator $\D_{(\theta,0)}$ over the cylinder $(-2T,\infty)\times Y$  as in \subsecref{subsec:gluing}.
Then, by the contraction mapping principle, we can find a genuine monopole $y_i$ near $x_i^\prime$ for each $i$. 
Taking a large $T$, we may assume 
$$
d_{K^\tau}([x]|_{K^\tau}, [y_1^\prime]\times[y_2^\prime]|_{K^\tau})<\frac{\varepsilon}2.
$$
So the monopole class $[x]$ is in the image $\mathfrak{F}(\mathrm{Gl})$ for gluing $y_1$ and $y_2$.
By \propref{prop:injective}, we find the inverse image of $[x]$ for the gluing map $\Xi$, and we can see that $\Xi$ is a diffeomorphism.
Thus \thmref{thm:gluing}  is proved.

\end{document}